\documentclass{article}
\usepackage{amsmath}
\usepackage{amssymb}

\usepackage{theorem}
\numberwithin{equation}{section}

\newtheorem{theorem}{Theorem}[section]
\newtheorem{conjecture}{Conjecture}[section]
\newtheorem{proposition}[theorem]{Proposition}
\newtheorem{corollary}[theorem]{Corollary}
\newtheorem{lemma}[theorem]{Lemma}
\newtheorem{definition}{Definition}[section]

\begin{document}
\title{Global well-posedness and scattering for the radial, defocusing, cubic wave equation with almost sharp initial data}
\date{\today}
\author{Benjamin Dodson}
\maketitle

\noindent \textbf{Abstract:} In this paper we prove that the cubic wave equation is globally well - posed and scattering for radial initial data lying in a weighted, slightly supercritical space. This space of functions is slightly smaller than the general critical space.
%satisfying

%\begin{equation}
%\| |x|^{2 \epsilon} u_{0} \|_{\dot{H}^{1/2 + \epsilon}(\mathbf{R}^{3})} + \| u_{0} \|_{\dot{H}^{1/2 + \epsilon}(\mathbf{R}^{3})} + \| |x|^{2 \epsilon} u_{1} \|_{\dot{H}^{-1/2 + \epsilon}(\mathbf{R}^{3})} + \| u_{1} \|_{\dot{H}^{-1/2 + \epsilon}(\mathbf{R}^{3})} < \infty.
%\end{equation}

\section{Introduction}
The defocusing, cubic nonlinear wave equation

\begin{equation}\label{1.1}
\partial_{t}^{2} u - \Delta u = -u^{3}, \hspace{5mm} u(0, \cdot) = u_{0}, \hspace{5mm} u_{t}(0, \cdot) = u_{1},
\end{equation}
is invariant under the scaling

\begin{equation}\label{1.1.1}
u(t,x) \mapsto \lambda u(\lambda t, \lambda x).
\end{equation}
That is, for any $\lambda > 0$, $(\ref{1.1.1})$ solves $(\ref{1.1})$ with initial data $(\lambda u_{0}(\lambda x)$, $\lambda^{2} u_{1}(\lambda x))$. This scaling preserves the $\dot{H}^{1/2} \times \dot{H}^{-1/2}$ norm of the initial data:

\begin{equation}\label{1.2}
\| \lambda u_{0}(\lambda x) \|_{\dot{H}^{1/2}(\mathbf{R}^{3})} = \| u_{0} \|_{\dot{H}^{1/2}(\mathbf{R}^{3})}, \hspace{5mm} \text{and} \hspace{5mm} \| \lambda^{2} u_{1}(\lambda x) \|_{\dot{H}^{-1/2}(\mathbf{R}^{3})} = \| u_{1} \|_{\dot{H}^{-1/2}(\mathbf{R}^{3})}.
\end{equation}

\noindent Because of this fact, $\dot{H}^{1/2} \times \dot{H}^{-1/2}$ is called the critical Sobolev space for $(\ref{1.1})$.\vspace{5mm}

\noindent A complete dichotomy has been proved for the local theory of $(\ref{1.1})$. On the negative side, \cite{CCT} proved that $(\ref{1.1})$ is locally ill - posed for $u_{0} \in \dot{H}^{s}$, $u_{1} \in \dot{H}^{s - 1}$ for any $s < \frac{1}{2}$. The proof exploited the scaling in $(\ref{1.1.1})$ in order to prove the existence of a solution that was not continuous in time at $t = 0$.\vspace{5mm}

\noindent On the positive side, \cite{LS} proved that $(\ref{1.1})$ is locally well - posed for initial data lying in the Sobolev space $\dot{H}^{1/2}(\mathbf{R}^{3}) \times \dot{H}^{-1/2}(\mathbf{R}^{3})$. Moreover, for data lying in $\dot{H}^{s} \times \dot{H}^{s - 1}$, $\frac{1}{2} < s < \frac{3}{2}$, \cite{LS} proved that $(\ref{1.1})$ is locally well - posed, with the time of existence depending only on the size of the initial data.\vspace{5mm}

\noindent \textbf{Remark:} In this paper, well - posedness refers to the standard definition that a unique solution that is continuous in time exists, and that the solution depends continuously on the initial data.\vspace{5mm}

\noindent A natural question in light of the result of \cite{LS} is when can local well - posedness results be extended to global well - posedness. It is a well - known fact from ordinary differential equations theory that this does not always hold. For example, consider the equation

\begin{equation}\label{1.2.1}
u_{tt} = u^{3}.
\end{equation}

\noindent The solution $u(t) = \sqrt{2} t^{-1}$ to $(\ref{1.2.1})$ certainly blows up in finite time. The paper \cite{AM} proved that solutions to

\begin{equation}\label{1.2.2}
u_{tt} - \Delta u = |u|^{p - 1} u,
\end{equation}

\noindent for $1 < p < 5$, $d = 3$, also can exhibit such ordinary differential equations type growth. On the other hand, for the ordinary differential equation

\begin{equation}\label{1.2.3}
u_{tt} = -u^{3},
\end{equation}

\noindent $u(t)$ is concave down in time when $u$ is positive, and is concave up in time when $u$ is negative, so the solution to $(\ref{1.2.3})$ should not blow up, but rather be global. It is for this reason, as well as due to the fact that

\begin{equation}\label{1.2.4}
\Delta u - u^{3} = 0,
\end{equation}

\noindent has no nonzero solutions, that it is conjectured that

\begin{conjecture}\label{c1.1}
$(\ref{1.1})$ is globally well - posed and scattering in $\dot{H}^{1/2}(\mathbf{R}^{3}) \times \dot{H}^{-1/2}(\mathbf{R}^{3})$.
\end{conjecture}

\noindent Scattering also has the usual definition.

\begin{definition}[Scattering]\label{d1.2}
A solution to $(\ref{1.1})$ is said to scatter forward in time if the solution exists for all $t > 0$ and there exists $(u_{0}^{+}, u_{1}^{+}) \in \dot{H}^{1/2}(\mathbf{R}^{3}) \times \dot{H}^{-1/2}(\mathbf{R}^{3})$ such that

\begin{equation}\label{1.3}
\aligned
\lim_{t \rightarrow \infty} \| u(t) - S(t)(u_{0}^{+}, u_{1}^{+}) \|_{\dot{H}^{1/2}(\mathbf{R}^{3})} = 0, \\
\lim_{t \rightarrow \infty} \| \partial_{t} u(t) - \partial_{t} S(t)(u_{0}^{+}, u_{1}^{+}) \|_{\dot{H}^{-1/2}(\mathbf{R}^{3})} = 0,
\endaligned
\end{equation}
where $S(t)$ is the solution operator to the linear wave equation $u_{tt} - \Delta u = 0$. Similarly, $u$ is said to scatter backward in time if there exists $(u_{0}^{-}, u_{1}^{-}) \in \dot{H}^{1/2} \times \dot{H}^{-1/2}$ such that $(\ref{1.3})$ holds for $t \rightarrow -\infty$.
\end{definition}
The counterpart to conjecture $\ref{c1.1}$ holds for the quintic problem.

\begin{theorem}\label{t1.2}
The defocusing, quintic wave equation problem

\begin{equation}\label{1.4}
\partial_{t}^{2} u - \Delta u = -u^{5}, \hspace{5mm} u(0, \cdot) = u_{0}, \hspace{5mm} u_{t}(0, \cdot) = u_{1},
\end{equation}
is globally well - posed and scatters both forward and backward in time for any $(u_{0}, u_{1}) \in \dot{H}^{1} \times L^{2}$.
\end{theorem}
\emph{Proof:} See \cite{Gril}. $\Box$\vspace{5mm}

\noindent \textbf{Remark:} In this case $(\ref{1.4})$ is invariant under the scaling $\lambda^{1/2} u(\lambda t, \lambda x)$, which preserves the $\dot{H}^{1}(\mathbf{R}^{3}) \times L^{2}(\mathbf{R}^{3})$.\vspace{5mm}

\noindent $(\ref{1.4})$ is called an energy - critical problem because a solution to $(\ref{1.4})$ conserves the energy

\begin{equation}\label{1.5}
 E(u(t)) = \frac{1}{2} \int |\nabla u(t,x)|^{2} dx + \frac{1}{2} \int (u_{t}(t,x))^{2} dx + \frac{1}{6} \int (u(t,x))^{6} dx = E(u(0))
 \end{equation}
 for the entire time of its existence. This quantity is unchanged by the scaling $\lambda^{1/2} u(\lambda t, \lambda x)$, controls the $\dot{H}^{1}(\mathbf{R}^{3}) \times L^{2}(\mathbf{R}^{3})$, and by the Sobolev embedding theorem 
 
 \begin{equation}\label{1.6}
 E(u(0)) \lesssim \| u_{0} \|_{\dot{H}^{1}(\mathbf{R}^{3})}^{2} + \| u_{0} \|_{\dot{H}^{1}(\mathbf{R}^{3})}^{6} + \| u_{1} \|_{L^{2}(\mathbf{R}^{3})}^{2}.
 \end{equation}
 However, for the cubic initial value problem there does not exist a conserved quantity which controls the $\dot{H}^{1/2} \times \dot{H}^{-1/2}$ norm of $(u(t), u_{t}(t))$. In fact, for radial data, the existence of such a quantity would prove conjecture $\ref{c1.1}$.

 \begin{theorem}\label{t1.3}
Suppose $u_{0} \in \dot{H}^{1/2}(\mathbf{R}^{3})$, $u_{1} \in \dot{H}^{-1/2}(\mathbf{R}^{3})$ are radial functions, and $u$ solves $(\ref{1.1})$ on a maximal interval $0 \in I \subset \mathbf{R}$ with

\begin{equation}\label{1.7}
\sup_{t \in I} \| u(t) \|_{\dot{H}^{1/2}(\mathbf{R}^{3})} + \| u_{t}(t) \|_{\dot{H}^{-1/2}(\mathbf{R}^{3})} < \infty.
\end{equation}
Then $I = \mathbf{R}$ and the solution $u$ scatters both forward and backward in time.
\end{theorem}
\emph{Proof:} See \cite{DL}. $\Box$\vspace{5mm}

\noindent In this paper we prove that $(\ref{1.1})$ is globally well - posed and scattering for $u_{0}$, $u_{1}$ contained in a subset of $\dot{H}^{1/2}(\mathbf{R}^{3}) \times \dot{H}^{-1/2}(\mathbf{R}^{3})$.

\begin{theorem}\label{t1.4}
Suppose there exists a positive constant $\epsilon > 0$ such that

\begin{equation}\label{1.8}
\| u_{0} \|_{\dot{H}^{1/2 + \epsilon}(\mathbf{R}^{3})} + \| |x|^{2 \epsilon} u_{0} \|_{\dot{H}^{1/2 + \epsilon}(\mathbf{R}^{3})} \leq A < \infty,
\end{equation}
and

\begin{equation}\label{1.9}
\| u_{1} \|_{\dot{H}^{-1/2 + \epsilon}(\mathbf{R}^{3})} + \| |x|^{2 \epsilon} u_{1} \|_{\dot{H}^{-1/2 + \epsilon}(\mathbf{R}^{3})} \leq A < \infty.
\end{equation}
Then $(\ref{1.1})$ has a global solution and there exists some $C(A, \epsilon) < \infty$ such that

\begin{equation}\label{1.10}
\int_{\mathbf{R}} \int (u(t,x))^{4} dx dt \leq C(A, \epsilon),
\end{equation}
which proves that $u$ scatters both forward and backward in time.
\end{theorem}

\noindent \textbf{Remark:} Theorem $\ref{t1.4}$ with $\epsilon = 0$ would imply conjecture $\ref{c1.1}$.\vspace{5mm}

\noindent The proof of theorem $\ref{t1.4}$ is based on two previous results.

 \begin{theorem}[Global well - posedness]\label{t1.5}
 For any $\epsilon > 0$, if $u_{0}$ and $u_{1}$ are radial functions and $u_{0} \in \dot{H}^{1/2}(\mathbf{R}^{3}) \cap \dot{H}^{1/2 + \epsilon}(\mathbf{R}^{3})$, $u_{1} \in \dot{H}^{-1/2}(\mathbf{R}^{3}) \cap \dot{H}^{-1/2 + \epsilon}(\mathbf{R}^{3})$, then $(\ref{1.1})$ is globally well - posed.
 \end{theorem}
 \emph{Proof:} See \cite{D1}. $\Box$\vspace{5mm}

\noindent \textbf{Remark:} The initial data in theorem $\ref{t1.4}$ satisfies $u_{0} \in \dot{H}^{1/2 + \epsilon}(\mathbf{R}^{3}) \cap \dot{H}^{1/2 - \epsilon}(\mathbf{R}^{3})$ and $u_{1} \in \dot{H}^{-1/2 + \epsilon}(\mathbf{R}^{3}) \cap \dot{H}^{-1/2 - \epsilon}(\mathbf{R}^{3})$, and thus $(\ref{1.1})$ has a global solution under such initial data.\vspace{5mm}

\noindent The proof of theorem $\ref{t1.5}$ used the I - method, which is an improvement over the Fourier truncation method. For example, using the I - method, \cite{CKSTT} improved the results of \cite{B} for the nonlinear Schr{\"o}dinger equation. On the wave equation side, \cite{Roy} extended the results of \cite{KPV} for the cubic wave equation to $s > \frac{13}{18}$ and to $s > \frac{7}{10}$ if $u$ has radial symmetry. Perhaps more importantly, \cite{KTwm} proved a well - posedness result which was technically unattainable via the Fourier truncation method. See \cite{D1} for a more detailed discussion of the history of the I - method.\vspace{5mm}

 \noindent The second result utilized the conformal transformation and energy methods.
 
 \begin{theorem}\label{t1.6}
 Assume that $A$, $\epsilon$ are positive constants. Let $(u_{0}, u_{1}) \in \dot{H}^{1} \times L^{2}$ be radial initial data so that
 
 \begin{equation}\label{1.11}
 \int |\nabla u_{0}(x)|^{2} (1 + |x|)^{1 + 2 \epsilon} dx + \int |u_{1}(x)|^{2} (1 + |x|)^{1 + 2 \epsilon} dx \leq A^{2},
 \end{equation}
 then the solution to $(\ref{1.1})$ scatters in both time directions with
 
 \begin{equation}\label{1.12}
 \| u \|_{L_{t,x}^{4}(\mathbf{R} \times \mathbf{R}^{3})} \leq C(A, \epsilon) < \infty.
 \end{equation}
 \end{theorem}
 \emph{Proof:} See \cite{Shen}. $\Box$\vspace{5mm}

 \noindent The proof in \cite{Shen} follows four steps.
 
 \begin{enumerate}
 \item Since $u_{0} \in \dot{H}^{1}$ and $u_{1} \in L^{2}$, then by the conservation of energy $(\ref{1.1})$ has a global solution under the initial data in theorem $\ref{t1.6}$.
 
 \item Define the conformal transformation of $u$,
 
 \begin{equation}\label{1.13}
 v(\tau, y) = \frac{\sinh |y|}{|y|} e^{\tau} u(e^{\tau} \frac{\sinh |y|}{|y|} \cdot y, t_{0} + e^{\tau} \cosh |y|), \hspace{5mm} (\tau, y) \in \mathbf{R} \times \mathbf{R}^{3},
 \end{equation}
 $t_{0} < 0$ is a fixed constant. Then \cite{Shen} shows that $v(0, y)$ has finite energy.
 
 \item Prove that if $v$ solves the conformal wave equation,
 
 \begin{equation}\label{1.14}
 \| v \|_{L_{\tau, y}^{4}(\mathbf{R} \times \mathbf{R}^{3})} \lesssim E(v(0)).
 \end{equation}
 
 \item Show that $(\ref{1.14})$ implies that $\| u \|_{L_{t,x}^{4}(\mathbf{R} \times \mathbf{R}^{3})} < \infty$, which proves scattering.
 \end{enumerate}
 
 \noindent The conformal energy is invariant under the scaling $(\ref{1.1.1})$. Because the conformal energy is conserved, to prove theorem $\ref{t1.4}$ it suffices to prove that ``most" (in a sense that will be more fully defined later) of a solution $u$ to $(\ref{1.1})$ with initial data lying in $(\ref{1.8})$ and $(\ref{1.9})$ has uniformly bounded conformal energy. This implies scattering.\vspace{5mm}

 \noindent By theorem $\ref{t1.5}$, step one also holds for the initial data in theorem $\ref{t1.4}$. In place of step two, we will prove that $v(0, y) \in \dot{H}^{1/2 + \epsilon}(\mathbf{R}^{3})$ and $v_{\tau}(0, y) \in \dot{H}^{-1/2 + \epsilon}(\mathbf{R}^{3})$. Next, as in \cite{D1}, we will utilize the I - method to prove that $\| v \|_{L_{\tau, y}^{4}(\mathbf{R} \times \mathbf{R}^{3})} < \infty$. Step four may be copied directly from \cite{Shen}, proving theorem $\ref{t1.4}$.\vspace{5mm}

 \noindent \textbf{Acknowledgements:} The author of this paper was supported by NSF grant DMS - 1500424.

 \section{Linear Estimates and harmonic analysis}
 In this section we will collect several estimates concerning the linear wave equation and harmonic analysis. These estimates will be utilized throughout the paper. None of the results in this section are new, in fact all are very well - known.

 \begin{definition}[Littlewood - Paley partition of unity]\label{d7.1}
Suppose $\psi \in C_{0}^{\infty}(\mathbf{R}^{3})$ is a radial, decreasing function supported on $|x| \leq 2$, $\psi = 1$ on $|x| \leq 1$. Then for any $j \in \mathbf{Z}$ we define the Littlewood - Paley projection

\begin{equation}\label{7.1}
(P_{j} f)(x) = \mathcal F^{-1} ((\psi(\frac{\xi}{2^{j}}) - \psi(\frac{2 \xi}{2^{j}})) \hat{f}(\xi))(x),
\end{equation}

\noindent where $\mathcal F$ is the Fourier transform

\begin{equation}
\mathcal F f(\xi) = (2 \pi)^{-3/2} \int e^{-ix \cdot \xi} f(x) dx,
\end{equation}

\noindent and

\begin{equation}\label{7.2}
\mathcal F^{-1}(\hat{f}(\xi))(x) = (2 \pi)^{-3/2} \int e^{ix \cdot \xi} \hat{f}(\xi) d\xi.
\end{equation}

\noindent Also define the operators

\begin{equation}\label{7.3}
(P_{\leq j} f)(x) = \sum_{l \leq j} P_{l} f,
\end{equation}

\noindent and $P_{> j} = 1 - P_{\leq j}$.
\end{definition}

\noindent \textbf{Remark:} Since $\psi$ is a $C_{0}^{\infty}(\mathbf{R}^{3})$ function, $P_{j} f$ is the convolution of $f$ with a Schwartz function $\mathcal F^{-1} P_{j}(x)$ that satisfies 

\begin{equation}\label{7.4}
(\mathcal F^{-1} P_{j})(x) = \check{\psi}_{j}(x) \lesssim_{l} 2^{3j} (1 + 2^{j} |x|)^{-l}, \hspace{5mm} \text{for any} \hspace{5mm} l \in \mathbf{Z}.
\end{equation}

\noindent By direct computation this gives the Sobolev embedding theorem

\begin{equation}
\| P_{j} f \|_{L^{q}(\mathbf{R}^{3})} \lesssim 2^{j(\frac{3}{p} - \frac{3}{q})} \| P_{j} f \|_{L^{p}(\mathbf{R}^{3})},
\end{equation}
when $q \geq p$. There is also the radial Sobolev embedding theorem

\begin{equation}
\| |x| (P_{j} f) \|_{L^{\infty}} \lesssim 2^{j/2} \| P_{j} f \|_{L^{2}}.
\end{equation}

\noindent Next recall the energy estimate

\begin{theorem}[Energy estimate]\label{t7.2}
If $u$ solves $u_{tt} - \Delta u = 0$ on an interval $I$, with $t_{0} \in I$, then

\begin{equation}\label{7.5}
\| \nabla u(t) \|_{L_{t}^{\infty} L_{x}^{2}(I \times \mathbf{R}^{3})} + \| u_{t}(t) \|_{L_{t}^{\infty} L_{x}^{2}(I \times \mathbf{R}^{3})} = \| \nabla u(t_{0}) \|_{L^{2}(\mathbf{R}^{3})} + \| u_{t}(t_{0}) \|_{L^{2}(\mathbf{R}^{3})}.
\end{equation}
\end{theorem}
Also recall the Strichartz estimates of \cite{Stri}.

\begin{theorem}[Strichartz estimate]\label{t7.3}
If $u$ solves $u_{tt} - \Delta u = 0$ on an interval $I$, with $t_{0} \in I$, then

\begin{equation}\label{7.6}
\| u(t) \|_{L_{t,x}^{4}(I \times \mathbf{R}^{3})} \lesssim \| u(t_{0}) \|_{\dot{H}^{1/2}(\mathbf{R}^{3})} + \| u_{t}(t_{0}) \|_{\dot{H}^{-1/2}(\mathbf{R}^{3})}.
\end{equation}
\end{theorem}
Also recall the endpoint Strichartz estimate of \cite{KM} for the radial wave equation.

 \begin{theorem}[Endpoint Strichartz estimate]\label{t7.4}
 If $u$ solves $u_{tt} - \Delta u = 0$ on an interval $I$, with $t_{0} \in I$, then

\begin{equation}\label{7.6}
\| u(t) \|_{L_{t}^{2} L_{x}^{\infty}(I \times \mathbf{R}^{3})} \lesssim \| u(t_{0}) \|_{\dot{H}^{1}(\mathbf{R}^{3})} + \| u_{t}(t_{0}) \|_{L^{2}(\mathbf{R}^{3})}.
\end{equation}
\end{theorem}
\textbf{Remark:} Then by interpolation, for any $0 < \sigma < 1$, if $u$ is a radial solution to $u_{tt} - \Delta u = 0$, $0 < \sigma < 1$,

\begin{equation}\label{7.6.0}
\frac{1}{p} + \frac{3}{q} = \frac{3}{2} - \sigma,
\end{equation}
then

\begin{equation}\label{7.6.1}
\| u \|_{L_{t}^{p} L_{x}^{q}(I \times \mathbf{R}^{3})} \lesssim \| u(t_{0}) \|_{\dot{H}^{\sigma}(\mathbf{R}^{3})} + \| u_{t}(t_{0}) \|_{\dot{H}^{\sigma - 1}(\mathbf{R}^{3})}.
\end{equation}
\textbf{Remark:} This fact is also true when $u$ is nonradial. See \cite{GV}.

 \begin{theorem}[Local energy estimate]\label{t7.5}
 Let $B_{R} = \{ x : |x| \leq R \}$ and let $A_{R} = B_{2R} \setminus B_{R} = \{ x : R \leq |x| \leq 2R \}$. Also let $B_{l} = B_{2^{l}}$ and $A_{l} = A_{2^{l}}$. Then define the norms
 
 \begin{equation}\label{7.7}
 \| f \|_{l^{\infty} L_{t,x}^{2}(J \times \mathbf{R}^{3})} = \sup_{l \in \mathbf{Z}} 2^{-l/2} \| f \|_{L_{t,x}^{2}(J \times B_{l})}.
 \end{equation}
 Also let
 
 \begin{equation}\label{7.8}
 \| g \|_{l^{1} L_{t,x}^{2}(J \times \mathbf{R}^{3})} = \sum_{l \in \mathbf{Z}} 2^{l/2} \| g \|_{L_{t,x}^{2}(J \times A_{l})}.
 \end{equation}
 Also, for any $j \in \mathbf{Z}$ let
 
 \begin{equation}\label{7.8.1}
 \| f \|_{l_{j}^{\infty} L_{t,x}^{2}(J \times \mathbf{R}^{3})} = \sup_{l \in \mathbf{Z}, l \geq -j} 2^{-l/2} \| f \|_{L_{t,x}^{2}(J \times B_{l})}.
 \end{equation}
 Then if $u$ is a radial solution to $u_{tt} - \Delta u = 0$ on an interval $I$, with $t_{0} \in I$,
 
 \begin{equation}\label{7.9}
 \| \nabla u \|_{l^{\infty} L_{t,x}^{2}(J \times \mathbf{R}^{3})} + \| u_{t} \|_{l^{\infty} L_{t,x}^{2}(J \times \mathbf{R}^{3})} \lesssim \| \nabla u(t_{0}) \|_{L^{2}(\mathbf{R}^{3})} + \| u_{t}(t_{0}) \|_{L^{2}(\mathbf{R}^{3})}.
 \end{equation}
 
 \end{theorem}
 \emph{Proof:} The proof follows by the sharp Huygens principle. Without loss of generality suppose $t_{0} = 0$, $u(0) = f$, and $u_{t}(0) = g$. By time reversal symmetry we can assume $t \geq 0$. If $t \leq R$, then by the energy estimate (theorem $\ref{t7.2}$),
 
 \begin{equation}\label{7.10}
 \| \nabla u(t) \|_{L_{t,x}^{2}([-R, R] \times \mathbf{R}^{3})} + \| u_{t}(t) \|_{L_{t,x}^{2}([-R, R] \times \mathbf{R}^{3})} \lesssim R^{1/2} \| f \|_{\dot{H}^{1}(\mathbf{R}^{3})} + R^{1/2} \| g \|_{L_{x}^{2}(\mathbf{R}^{3})}.
 \end{equation}
 Next, since $u$ is radial, if $r \leq R$ and $t > R$,
 
 \begin{equation}\label{7.11}
 r u(t,r) = \frac{1}{2} (r + t) f(r + t) - \frac{1}{2} (t - r) f(t - r) + \frac{1}{2} \int_{t - r}^{t + r} s g(s) ds.
 \end{equation}
 Now then,
 
 \begin{equation}\label{7.12}
 \aligned
 \partial_{r} (r u(t,r)) = \frac{1}{2} f(r + t) + \frac{1}{2} (r + t) f'(r + t) + \frac{1}{2} f(t - r) + \frac{1}{2} (t - r) f'(t - r) \\ + \frac{1}{2} (t + r) g(t + r) + \frac{1}{2} (t - r) g(t - r), \\
  \partial_{t} (r u(t,r)) = \frac{1}{2} f(r + t) + \frac{1}{2} (r + t) f'(r + t) - \frac{1}{2} f(t - r) - \frac{1}{2} (t - r) f'(t - r) \\ + \frac{1}{2} (t + r) g(t + r) - \frac{1}{2} (t - r) g(t - r).
  \endaligned
  \end{equation}
 Now by Fubini's theorem, H{\"o}lder's inequality, and Hardy's inequality,
 
 \begin{equation}\label{7.13}
 \aligned
 \int_{R}^{\infty} \int_{0}^{R} f(r + t)^{2} dr dt + \int_{R}^{\infty} \int_{0}^{R} f(t - r)^{2} dr dt \\ \lesssim R \int_{0}^{\infty} f(s)^{2} ds \lesssim R \| \frac{1}{|x|} f \|_{L^{2}(\mathbf{R}^{3})}^{2} \lesssim R \| \nabla f \|_{L^{2}(\mathbf{R}^{3})}^{2},
 \endaligned
 \end{equation}
 while
 
 \begin{equation}\label{7.14}
 \aligned
 \int_{R}^{\infty} \int_{0}^{R} (r + t)^{2} f'(r + t)^{2} dr dt + \int_{R}^{\infty} \int_{0}^{R} f'(t - r)^{2} (t - r)^{2} dr dt \\ \lesssim R \int_{0}^{\infty} f'(s)^{2} s^{2} ds \lesssim R \| \nabla f \|_{L^{2}(\mathbf{R}^{3})}^{2},
 \endaligned
 \end{equation}
and

\begin{equation}\label{7.15}
\int_{R}^{\infty} \int_{0}^{R} (t + r)^{2} g(t + r)^{2} dr dt + \int_{R}^{\infty} \int_{0}^{R} (t - r)^{2} g(t - r)^{2} dr dt \lesssim R \| g \|_{L^{2}(\mathbf{R}^{3})}^{2}.
\end{equation}
Since $r u_{t} = \partial_{t} (r u)$ and $r u_{r} = \partial_{r} (ru) - u$, to complete the proof of theorem $\ref{t7.5}$ it only remains to show that

\begin{equation}\label{7.16}
\int_{R}^{\infty} \int_{0}^{R} u(t, r)^{2} dr dt \lesssim R \| \nabla f \|_{L^{2}(\mathbf{R}^{3})}^{2} + R \| g \|_{L^{2}(\mathbf{R}^{3})}^{2}.
\end{equation}
However this follows directly from the endpoint Strichartz estimate in theorem $\ref{t7.4}$ and H{\"o}lder's inequality. $\Box$

\begin{theorem}\label{t7.6}
If $u$ is a radial solution to

\begin{equation}\label{7.17}
u_{tt} - \Delta u = F_{1} + F_{2} + F_{3} + F_{4}, \hspace{5mm} u(0) = u_{0}, \hspace{5mm} u_{t}(0) = u_{1},
\end{equation}
then for any $j$,

\begin{equation}\label{7.18}
\aligned
\| P_{j} u \|_{L_{t}^{2} L_{x}^{\infty}(J \times \mathbf{R}^{3})} + 2^{j/2} \| P_{j} u \|_{L_{t,x}^{4}(J \times \mathbf{R}^{3})} + 2^{j} \| P_{j} u \|_{L_{t}^{\infty} L_{x}^{2}(J \times \mathbf{R}^{3})} + \| P_{j} u_{t} \|_{L_{t}^{\infty} L_{x}^{2}(J \times \mathbf{R}^{3})} \\
+ \| P_{j} u \|_{l_{j}^{\infty} L_{t,x}^{2}(J \times \mathbf{R}^{3})} \lesssim 2^{j} \| P_{j} u_{0} \|_{L^{2}(\mathbf{R}^{3})} + \| P_{j} u_{1} \|_{L^{2}(\mathbf{R}^{3})} + \| P_{j} F_{1} \|_{L_{t}^{1} L_{x}^{2}(J \times \mathbf{R}^{3})} \\
+ \| P_{j} F_{2} \|_{L_{t,x}^{4/3}(J \times \mathbf{R}^{3})} + \| P_{j} F_{3} \|_{l^{1} L_{t,x}^{2}(J \times \mathbf{R}^{3})} + 2^{j} \| P_{j} F_{4} \|_{L_{t}^{2} L_{x}^{1}(J \times \mathbf{R}^{3})}.
\endaligned
\end{equation}
To simplify notation let

\begin{equation}\label{7.18.1}
\aligned
\mathcal S(j, u) = \| P_{j} u \|_{L_{t}^{2} L_{x}^{\infty}(J \times \mathbf{R}^{3})} + 2^{j/2} \| P_{j} u \|_{L_{t,x}^{4}(J \times \mathbf{R}^{3})} + 2^{j} \| P_{j} u \|_{L_{t}^{\infty} L_{x}^{2}(J \times \mathbf{R}^{3})} \\ + \| P_{j} u_{t} \|_{L_{t}^{\infty} L_{x}^{2}(J \times \mathbf{R}^{3})} + \| P_{j} u \|_{l_{j}^{\infty} L_{t,x}^{2}(J \times \mathbf{R}^{3})}.
\endaligned
\end{equation}
\end{theorem}
\emph{Proof:} Let $\tilde{P}_{j} = P_{j - 1} + P_{j} + P_{j + 1}$, so that, by definition $\ref{d7.1}$, $\tilde{P}_{j} P_{j} = P_{j}$. By Duhamel's principle

\begin{equation}\label{7.20}
u(t) = S(t)(u_{0}, u_{1}) + \int_{0}^{t} S(t - \tau)(0, F_{1} + F_{2} + F_{3} + F_{4}) d\tau
\end{equation}
solves $(\ref{7.17})$, and

\begin{equation}\label{7.21}
P_{j} u(t) = S(t)(P_{j} u_{0}, P_{j} u_{1}) + \int_{0}^{t} S(t - \tau)(0, P_{j} F_{1} + P_{j} F_{2} + P_{j} F_{3} + P_{j} F_{4}) d\tau
\end{equation}
By theorems $\ref{t7.2}$, $\ref{t7.3}$, and $\ref{t7.4}$,

\begin{equation}\label{7.22}
\aligned
\| P_{j} S(t) (u_{0}, u_{1}) \|_{L_{t}^{2} L_{x}^{\infty}(J \times \mathbf{R}^{3})} + 2^{j/2} \| P_{j} S(t)(u_{0}, u_{1}) \|_{L_{t,x}^{4}(J \times \mathbf{R}^{3})} \\ + 2^{j} \| P_{j} S(t)(u_{0}, u_{1}) \|_{L_{t}^{\infty} L_{x}^{2}(J \times \mathbf{R}^{3})} + \| P_{j} \partial_{t} S(t)(u_{0}, u_{1}) \|_{L_{t}^{\infty} L_{x}^{2}(J \times \mathbf{R}^{3})} \\ \lesssim \| \nabla P_{j} u_{0} \|_{L^{2}(\mathbf{R}^{3})} + \| P_{j} u_{1} \|_{L^{2}(\mathbf{R}^{3})} \lesssim 2^{j} \| P_{j} u_{0} \|_{L^{2}(\mathbf{R}^{3})} + \| P_{j} u_{1} \|_{L^{2}(\mathbf{R}^{3})}.
\endaligned
\end{equation}
Meanwhile, by theorem $\ref{t7.5}$,

\begin{equation}\label{7.23}
\| \nabla P_{j} S(t)(u_{0}, u_{1}) \|_{l^{\infty} L_{t,x}^{2}(J \times \mathbf{R}^{3})} \lesssim 2^{j} \| P_{j} u_{0} \|_{L^{2}(\mathbf{R}^{3})} + \| P_{j} u_{1} \|_{L^{2}(\mathbf{R}^{3})}.
\end{equation}
To simplify notation let $v_{j} = P_{j} S(t)(u_{0}, u_{1})$. Now, by $(\ref{7.4})$, if $R \geq 2^{-j}$, $|x| \leq R$, by Taylor's formula and $\int \check{\psi}_{j}(x) dx = 0$,

\begin{equation}\label{7.24}
\aligned
P_{j} v_{j}(t,x) = 2^{3j} \int \check{\psi}_{j}(2^{j}(x - y)) v_{j}(t,y) dy = 2^{3j} \int \check{\psi}_{j}(2^{j}(x - y)) [v_{j}(t,y) - v_{j}(t,x)] dy \\ = 2^{3j} \int_{0}^{1} \int \check{\psi}_{j}(2^{j}(x - y)) (y - x) \cdot \nabla v_{j}(t, x + \tau(y - x)) dy d\tau,
\endaligned
\end{equation}
and so by $(\ref{7.4})$, if $R = 2^{-j + m}$, $m \geq 0$,

\begin{equation}\label{7.25}
\aligned
R^{-1/2} \| \tilde{P}_{j} P_{j} u \|_{L_{t,x}^{2}(J \times B_{R})} \lesssim R^{-1/2} \sum_{l \geq 0} 2^{-10 (l + m)} 2^{-j} R^{-1/2} \| \nabla v_{j} \|_{L_{t,x}^{2}(J \times B_{-j + m + l})} \\ \lesssim \| P_{j} u_{0} \|_{L^{2}(\mathbf{R}^{3})} + 2^{-j} \| P_{j} u_{1} \|_{L^{2}(\mathbf{R}^{3})}.
\endaligned
\end{equation}
Thus,

\begin{equation}\label{7.26}
\mathcal S(j, S(t)(u_{0}, u_{1})) \lesssim 2^{j} \| P_{j} u_{0} \|_{L^{2}(\mathbf{R}^{3})} + \| P_{j} u_{1} \|_{L^{2}(\mathbf{R}^{3})}.
\end{equation}
Next, by the principle of superposition combined with $(\ref{7.26})$, if

\begin{equation}\label{7.27}
u^{1}(t) = \int_{0}^{t} S(t - \tau)(0, F_{1}) d\tau,
\end{equation}
then

\begin{equation}\label{7.28}
\mathcal S(j, u^{1}(t)) \lesssim \| P_{j} F_{1} \|_{L_{t}^{1} L_{x}^{2}(J \times \mathbf{R}^{3})}.
\end{equation}
Next, by the Christ - Kiselev lemma, duality, theorem $\ref{t7.3}$, and the sine addition formulas (see \cite{GV}), if

\begin{equation}\label{7.29}
u^{2}(t) = \int_{0}^{t} S(t - \tau)(0, P_{j} F_{2}) d\tau, \hspace{5mm} \mathcal S(j, u^{2}(t)) \lesssim 2^{j/2} \| P_{j} F_{2} \|_{L_{t,x}^{4/3}(J \times \mathbf{R}^{3})}.
\end{equation}
Next, let

\begin{equation}\label{7.30}
u^{3}(t) = \int_{0}^{t} S(t - \tau)(0, P_{j} F_{3}) d\tau.
\end{equation}
First suppose that $F_{3}$ is supported on $r \leq R$. By the fundamental solution of the wave equation,

\begin{equation}\label{7.31}
r u^{3}(t, r) = \int_{\sup(0, t - r)}^{t} \int_{r - (t - \tau)}^{r + (t - \tau)} s F_{3}(s, \tau) ds d\tau + \int_{0}^{\sup(0, t - r)} \int_{(t - \tau) - r}^{(t - \tau) + r} s F_{3}(s, \tau) ds d\tau.
\end{equation}
By the support properties of $F_{3}$, the integrals 

\begin{equation}\label{7.32}
\int_{r - (t - \tau)}^{r + (t - \tau)} s F_{3}(s, \tau) ds, \hspace{5mm} \text{and} \hspace{5mm} \int_{(t - \tau) - r}^{(t - \tau) + r} s F_{3}(s, \tau) ds,
\end{equation}
 are only nonzero if $|(t - \tau) - r| \leq R$. Thus, if for some $k \in \mathbf{Z}$, $kR \leq (t - r) < (k + 1) R$, $(\ref{7.32})$ is zero unless $(k - 1) R \leq \tau \leq (k + 1) R$. Thus, the supports of
 
 \begin{equation}\label{7.33}
 \int_{kR}^{\inf((k + 1) R, t)} S(t - \tau)(0, F_{3}) d\tau
 \end{equation}
 are finitely intersecting. Now, by H{\"o}lder's inequality in time,
 
 \begin{equation}\label{7.34}
 \| F_{3} \|_{L_{t}^{1} L_{x}^{2}([kR, (k + 1) R] \times \mathbf{R}^{3})} \lesssim R^{1/2} \| F_{3} \|_{L_{t,x}^{2}([kR, (k + 1) R] \times \mathbf{R}^{3})},
 \end{equation}
so by theorems $\ref{t7.2}$ - $\ref{t7.5}$,

\begin{equation}\label{7.35}
\aligned
\| u^{3}(t) \|_{L_{t}^{2} L_{x}^{\infty}(J \times \mathbf{R}^{3})} + \| \nabla u^{3}(t) \|_{L_{t}^{\infty} L_{x}^{2}(J \times \mathbf{R}^{3})} + \| \partial_{t} u^{3}(t) \|_{L_{t}^{\infty} L_{x}^{2}(J \times \mathbf{R}^{3})} \\ + \| \nabla u^{3}(t) \|_{l^{\infty} L_{t,x}^{2}(J \times \mathbf{R}^{3})} + \| \partial_{t} u^{3}(t) \|_{l^{\infty} L_{t,x}^{2}(J \times \mathbf{R}^{3})} \lesssim R^{1/2} \| F_{3} \|_{l^{1} L_{t,x}^{2}(J \times \mathbf{R}^{3})}.
\endaligned
\end{equation}
If $F_{3}$ is compactly supported then $P_{j} F_{3}$ need not be, but making a calculation similar to $(\ref{7.23})$ - $(\ref{7.26})$ proves that

\begin{equation}\label{7.36}
\mathcal S(j, u^{3}(t)) \lesssim \| P_{j} F_{3} \|_{l^{1} L_{t,x}^{2}(J \times \mathbf{R}^{3})}.
\end{equation}
Finally let

\begin{equation}\label{7.37}
u^{4}(t) = \int_{0}^{t} S(t - \tau)(0, F_{4}) d\tau.
\end{equation}
Since $P_{j} = \tilde{P}_{j} P_{j}$,

\begin{equation}\label{7.38}
P_{j} u^{4}(t) = P_{j} \int_{0}^{t} S(t - \tau)(0, \tilde{P}_{j} F_{4}) d\tau.
\end{equation}
Now, by the Sobolev embedding theorem,

\begin{equation}\label{7.39}
\| \tilde{P}_{j} F_{4} \|_{L_{t,x}^{2}(J \times \mathbf{R}^{3})} \lesssim 2^{3j/2} \| \tilde{P}_{j} F_{4} \|_{L_{t}^{2} L_{x}^{1}(J \times \mathbf{R}^{3})},
\end{equation}
so in particular,

\begin{equation}\label{7.40}
2^{-j/2} \| \tilde{P}_{j} F_{4} \|_{L_{t,x}^{2}(J \times B_{-j})} \lesssim 2^{j} \| \tilde{P}_{j} F_{4} \|_{L_{t}^{2} L_{x}^{1}(J \times \mathbf{R}^{3})}.
\end{equation}
Also, by the radial Sobolev embedding theorem,

\begin{equation}\label{7.41}
\| |x| \tilde{P}_{j} F_{4} \|_{L_{t,x}^{2}(J \times \mathbf{R}^{3})} \lesssim 2^{j/2} \| \tilde{P}_{j} F_{4} \|_{L_{t}^{2} L_{x}^{1}(J \times \mathbf{R}^{3})}.
\end{equation}
Therefore,

\begin{equation}\label{7.42}
\sum_{l > -j} 2^{l/2} \| \tilde{P}_{j} F_{4} \|_{L_{t,x}^{2}(J \times A_{l})} \lesssim 2^{j/2} \sum_{l > -j} 2^{-l/2} \| \tilde{P}_{j} F_{4} \|_{L_{t}^{2} L_{x}^{1}(J \times \mathbf{R}^{3})} \lesssim 2^{j}  \| \tilde{P}_{j} F_{4} \|_{L_{t}^{2} L_{x}^{1}(J \times \mathbf{R}^{3})}.
\end{equation}
Thus, by $(\ref{7.36})$,

\begin{equation}\label{7.43}
\mathcal S(j, u^{4}(t)) \lesssim 2^{j} \| P_{j} \tilde{P}_{j} F_{4} \|_{L_{t}^{2} L_{x}^{1}(J \times \mathbf{R}^{3})} = 2^{j} \| P_{j} F_{4} \|_{L_{t}^{2} L_{x}^{1}(J \times \mathbf{R}^{3})}.
\end{equation}
This completes the proof of the theorem. $\Box$

\begin{theorem}\label{t7.7}
If $u$ is a radial solution to the wave equation

\begin{equation}\label{7.44}
u_{tt} - \Delta u = F_{1} + F_{2} + F_{3}, \hspace{5mm} u(t_{0}) = u_{0}, \hspace{5mm} u_{t}(t_{0}) = u_{1},
\end{equation}
on the interval $J$ with $t_{0} \in I$, and $0 < \sigma < 1$ satisfies

\begin{equation}\label{7.44.1}
\frac{1}{p} = 1 - \frac{\sigma}{2}, \hspace{5mm} \frac{1}{q} = \frac{1}{2} + \frac{\sigma}{2},
\end{equation}

\begin{equation}\label{7.45}
\aligned
\| \nabla u \|_{l^{\infty} L_{t,x}^{2}(J \times \mathbf{R}^{3})} + \| u_{t} \|_{l^{\infty} L_{t,x}^{2}(J \times \mathbf{R}^{3})} + \| u \|_{L_{t}^{2} L_{x}^{\infty}(J \times \mathbf{R}^{3})} \lesssim \| u_{0} \|_{\dot{H}^{1}(\mathbf{R}^{3})} + \| u_{1} \|_{L^{2}(\mathbf{R}^{3})} \\ + \| F_{1} \|_{L_{t}^{1} L_{x}^{2}(J \times \mathbf{R}^{3})} + \| |\nabla|^{\sigma} F_{2} \|_{L_{t}^{p} L_{x}^{q}(J \times \mathbf{R}^{3})} + \| F_{3} \|_{l^{1} L_{t,x}^{2}(J \times \mathbf{R}^{3})}.
\endaligned
\end{equation}
\end{theorem}
\emph{Proof:} By the endpoint Strichartz estimate (theorem $\ref{t7.4}$) and the local energy - estimate (theorem $\ref{t7.5}$),

\begin{equation}\label{7.46}
\aligned
\| \nabla S(t)(u_{0}, u_{1}) \|_{l^{\infty} L_{t,x}^{2}(J \times \mathbf{R}^{3})} + \| \partial_{t} S(t)(u_{0}, u_{1}) \|_{l^{\infty} L_{t,x}^{2}(J \times \mathbf{R}^{3})} \\ + \| S(t)(u_{0}, u_{1}) \|_{l^{\infty} L_{t,x}^{2}(J \times \mathbf{R}^{3})} \|_{L_{t}^{2} L_{x}^{\infty}(J \times \mathbf{R}^{3})} \lesssim \| u_{0} \|_{\dot{H}^{1}(\mathbf{R}^{3})} + \| u_{1} \|_{L^{2}(\mathbf{R}^{3})}.
\endaligned
\end{equation}
Also by the principle of superposition, if $u^{1}(t) = \int_{0}^{t} S(t - \tau) F_{1}(\tau) d\tau$,

\begin{equation}\label{7.47}
\| \nabla u^{1} \|_{l^{\infty} L_{t,x}^{2}(J \times \mathbf{R}^{3})} + \| \partial_{t} u^{1} \|_{l^{\infty} L_{t,x}^{2}(J \times \mathbf{R}^{3})} + \| u^{1} \|_{L_{t}^{2} L_{x}^{\infty}(J \times \mathbf{R}^{3})} \lesssim \| F_{1} \|_{L_{t}^{1} L_{x}^{2}(J \times \mathbf{R}^{3})}.
\end{equation}
Next, by the Christ - Kiselev lemma (again see \cite{GV}), when $\sigma < 1$, if $u^{2}(t) = \int_{0}^{t} S(t - \tau) F_{2}(\tau) d\tau$,

\begin{equation}\label{7.48}
\| u^{2} \|_{l^{\infty} L_{t,x}^{2}(J \times \mathbf{R}^{3})} + \| \partial_{t} u^{2} \|_{l^{\infty} L_{t,x}^{2}(J \times \mathbf{R}^{3})} + \| u^{2} \|_{L_{t}^{2} L_{x}^{\infty}(J \times \mathbf{R}^{3})} \lesssim \| |\nabla|^{\sigma} F_{2} \|_{L_{t}^{p} L_{x}^{q}(J \times \mathbf{R}^{3})}.
\end{equation}
Finally, as in $(\ref{7.33})$, if $F_{3}$ is supported on $|x| \leq R$, then the supports of

\begin{equation}\label{7.49}
\int_{kR}^{\inf(t, (k + 1) R} S(t - \tau) (0, F_{3})(\tau) d\tau,
\end{equation}
are finitely intersecting, so since by H{\"o}lder's inequality in time,

\begin{equation}\label{7.50}
\| F_{3} \|_{L_{t}^{1} L_{x}^{2}([kR, (k + 1)R] \times \mathbf{R}^{3})} \lesssim R^{1/2} \| F_{3} \|_{L_{t,x}^{2}([kR, (k + 1) R] \times \mathbf{R}^{3})},
\end{equation}
which proves that if $u^{3}(t) = \int_{0}^{t} S(t - \tau)(0, F_{3}) d\tau$,

\begin{equation}\label{7.51}
\| u^{3} \|_{l^{\infty} L_{t,x}^{2}(J \times \mathbf{R}^{3})} + \| \partial_{t} u^{3} \|_{l^{\infty} L_{t,x}^{2}(J \times \mathbf{R}^{3})} + \| u^{3} \|_{L_{t}^{2} L_{x}^{\infty}(J \times \mathbf{R}^{3})} \lesssim \| F_{3} \|_{L_{t}^{1} L_{x}^{2}(J \times \mathbf{R}^{3})}.
\end{equation}
$\Box$

\section{The conformal symmetry for the wave equation}
In this section we outline the proof of the main theorem, theorem $\ref{t1.4}$. The proof of the theorem depends on several propositions, which will then be proved in subsequent sections. First, observe that by theorem $\ref{t1.5}$, if $u_{0} \in \dot{H}^{1/2 + \epsilon}(\mathbf{R}^{3}) \cap \dot{H}^{1/2}(\mathbf{R}^{3})$ and $u_{1} \in \dot{H}^{-1/2 + \epsilon}(\mathbf{R}^{3}) \cap \dot{H}^{-1/2}(\mathbf{R}^{3})$ are radial functions, then $(\ref{1.1})$ has a global solution. We will show in lemma $\ref{l8.1}$ that this is true.\vspace{5mm}

\noindent Now let $v = \mathbf{T} u$, where $\mathbf{T}$ is the conformal transformation

\begin{equation}\label{6.1}
(\mathbf{T} u)(y, \tau) = \frac{\sinh |y|}{|y|} e^{\tau} u(\tilde{\mathbf{T}} (y, \tau)) = \frac{\sinh |y|}{|y|} e^{\tau} u(e^{\tau} \frac{\sinh |y|}{|y|} y, t_{0} + e^{\tau} \cosh(|y|)), \hspace{5mm} (y, \tau) \in \mathbf{R}^{3} \times \mathbf{R},
\end{equation}
where $t_{0} < 0$ will be defined later. Since $u$ is radially symmetric, taking $s \in [0, \infty)$, $\tau \in \mathbf{R}$,

\begin{equation}\label{6.2}
v(s, \tau) = \frac{\sinh s}{s} e^{\tau} u(e^{\tau} \sinh s, t_{0} + e^{\tau} \cosh s).
\end{equation}
Then by direct computation, taking $w(r,t) = r u(r,t)$,

\begin{equation}\label{6.3}
v_{\tau \tau} - \Delta v = v_{\tau \tau} - v_{ss} - \frac{2}{s} v_{s} = \frac{1}{s} [(s v)_{\tau \tau} - (s v)_{ss}] = -(\frac{s}{\sinh s})^{2} v^{3} = -\phi(s) v^{3}.
\end{equation}
\textbf{Remark:} See section five of \cite{Shen} for more details concerning the calculation.\vspace{5mm}

\noindent We prove

\begin{equation}\label{6.4}
\int_{\mathbf{R}} \int_{0}^{\infty} (v(s, \tau))^{4} \phi(s) s^{2} ds d\tau < \infty
\end{equation}
in four steps.

\begin{proposition}[Initial data]\label{p6.1}
There exists some constant $C_{0}(A, \epsilon) < \infty$ such that

\begin{equation}\label{6.5}
\| v(s, 0) \|_{\dot{H}^{1/2 + \epsilon}(\mathbf{R}^{3})} + \| \partial_{\tau} v(s, 0) \|_{\dot{H}^{-1/2 + \epsilon}(\mathbf{R}^{3})} \lesssim_{A, \epsilon} 1.
\end{equation}
\end{proposition}
%\textbf{Remark:} For technical reasons we will also require that $C_{0}(A, \epsilon)$ be larger than the constant in Hardy's inequality, that is,

%\begin{equation}
%\| \frac{1}{|x|} f \|_{L^{2}(\mathbf{R}^{3})} \leq C_{0}(A, \epsilon) \| \nabla f \|_{L^{2}(\mathbf{R}^{3})}.
%\end{equation}
\noindent Now let $N = 2^{k_{0}}$ and define the I - operator

\begin{equation}\label{6.6}
I = P_{\leq k_{0}} + \sum_{j > 0} 2^{-j(\frac{1}{2} - \epsilon)} P_{j + k_{0}}.
\end{equation}
Then

\begin{equation}
\frac{1}{2} \| \nabla Iv \|_{L^{2}(\mathbf{R}^{3})}^{2} + \frac{1}{2} \| Iv_{\tau} \|_{L^{2}(\mathbf{R}^{3})}^{2} \lesssim_{A, \epsilon} N^{1 - 2 \epsilon}.
\end{equation}
Also, by the Sobolev embedding theorem and H{\"o}lder's inequality,

\begin{equation}
\frac{1}{4} \int \phi(s) |Iv(s, 0)|^{4} s^{2} ds \lesssim_{A, \epsilon} N^{1 - 2 \epsilon}.
\end{equation}
Therefore, there exists some $C_{0}(A, \epsilon)$ such that

\begin{equation}
E(Iv(0)) \leq C_{0}(A, \epsilon) N^{1 - 2 \epsilon}.
\end{equation}

\begin{proposition}[Long time Strichartz estimate]\label{p6.2}
Suppose $\mathcal J$ is an interval on which

\begin{equation}\label{6.7}
\int_{\mathcal J} \int \phi(x) (Iv(x,t))^{4} dx dt \leq C_{1} N^{1 - 2 \epsilon},
\end{equation}
and $\sup_{t \in \mathcal J} E(Iv(t)) \leq 2 C_{0}(A, \epsilon) N^{1 - 2 \epsilon}$. Then there exists $k_{0}(C_{1}, A, \epsilon)$ sufficiently large such that,

\begin{equation}\label{6.8}
\sup_{j \geq k_{0} - 7} 2^{j(\epsilon - 1/2)} \mathcal S(j, v) \lesssim C_{0}(A, \epsilon).
\end{equation}
\end{proposition}
\textbf{Remark:} It is important to observe that the implicit constant in $(\ref{6.8})$ crucially does not depend on $C_{1}$.

\begin{proposition}[Almost conservation of energy]\label{p6.3}
Suppose $\mathcal J$ is an interval on which

\begin{equation}\label{6.9}
\int_{\mathcal J} \int_{0}^{\infty} \phi(s) (\frac{\cosh s}{\sinh s}) (Iv(s,\tau))^{4} s^{2} ds d\tau \leq C_{1} N^{1 - 2 \epsilon},
\end{equation}
 $\sup_{t \in \mathcal J} E(Iv(t)) \leq 2 C_{0}(A, \epsilon) N^{1 - 2 \epsilon}$, and $E(Iv(0)) \leq C_{0}(A, \epsilon) N^{1 - 2 \epsilon}$. Then for $k_{0}(C_{1}, A, \epsilon)$ sufficiently large,

\begin{equation}\label{6.10}
\sup_{\tau \in \mathcal J} E(Iv(\tau)) \leq \frac{3}{2} C_{0}(A, \epsilon) N^{1 - 2 \epsilon}.
\end{equation}
\end{proposition}

\begin{proposition}[Almost Morawetz estimate]\label{p6.4}
If $\mathcal J$ is an interval on which $E(Iv(\tau)) \leq 2 C_{0}(A, \epsilon) N^{1 - 2 \epsilon}$, and

\begin{equation}
\int_{\mathcal J} \int \phi(s) (\frac{\cosh s}{\sinh s}) (Iv(s, \tau))^{4} s^{2} ds d\tau,
\end{equation}
then for $k_{0}(C_{1}, A, \epsilon)$ sufficiently large,

\begin{equation}\label{6.11}
\int_{\mathcal J} \int \phi(s)(\frac{\cosh s}{\sinh s}) (Iv(s, \tau))^{4} s^{2} ds d\tau \lesssim C_{0}(A, \epsilon) N^{1 - 2\epsilon}.
\end{equation}
\end{proposition}
Armed with propositions $\ref{p6.1}$ - $\ref{p6.4}$,

\begin{equation}\label{6.12}
\int_{\mathbf{R}} \int_{0}^{\infty} \phi(s) v(s, \tau)^{4} s^{2} ds d\tau < \infty
\end{equation}
may be proved by a bootstrap argument. Choose some $C_{1}(A, \epsilon) >> C_{0}(A, \epsilon)$ and define the set

\begin{equation}\label{6.13}
\aligned
J_{0} = \{ T \in \mathbf{R} : \int_{-T}^{T} \int_{0}^{\infty} \phi(s) (\frac{\cosh s}{\sinh s}) (Iv(s, \tau))^{4} s^{2} ds d\tau \leq C_{1} N^{1 - 2 \epsilon}, \\ \hspace{5mm} \text{and} \hspace{5mm} \sup_{\tau \in [-T, T]} E(Iv(\tau)) \leq 2 N^{1 - 2 \epsilon} \}.
\endaligned
\end{equation}
%Now by Hardy's inequality and the Sobolev embedding theorem,

%\begin{equation}\label{6.14}
%\aligned
%\int \phi(s) (\frac{\cosh s}{\sinh s}) (Iv(s, \tau))^{4} s^{2} ds \lesssim \| \frac{1}{|x|} Iv(\tau) \|_{L_{x}^{2}(\mathbf{R}^{3})} \| Iv(\tau) \|_{L_{x}^{6}(\mathbf{R}^{3})}^{3} + \| Iv(\tau) \|_{L_{x}^{4}(\mathbf{R}^{3})}^{4} \\ \lesssim E(Iv(\tau))^{2} + E(Iv(\tau)) \lesssim N^{2 - 4 \epsilon} \| v(\tau) \|_{\dot{H}^{1/2 + \epsilon}(\mathbf{R}^{3})}^{4}(1 + \| v(\tau) \|_{\dot{H}^{1/2}(\mathbf{R}^{3})}^{4}) \\ + N^{1 - 2 \epsilon} \| v(\tau) \|_{\dot{H}^{1/2 + \epsilon}(\mathbf{R}^{3})}^{2}(1 + \| v(\tau) \|_{\dot{H}^{1/2}(\mathbf{R}^{3})}^{2}).
%\endaligned
%\end{equation}
%Then by theorem $\ref{t1.5}$ this quantity is bounded for any $\tau \in \mathbf{R}$, this implies that $J_{0}$ contains some open interval around $0$. 

\noindent By the dominated convergence theorem, $J_{0}$ is a closed set. Also, clearly $0 \in J_{0}$. Therefore, it suffices to prove that $J_{0}$ is open, which would then imply that $J_{0} = \mathbf{R}$.\vspace{5mm}

\noindent Now by proposition $\ref{p6.3}$,

\begin{equation}\label{6.15}
E(Iv(\tau)) \leq \frac{3}{2} C_{0}(A, \epsilon) N^{1 - 2 \epsilon}
\end{equation}
for all $\tau \in J_{0}$. Also, proposition $\ref{p6.4}$ implies

\begin{equation}\label{6.16}
\int_{J_{0}} \int \phi(s) (\frac{\cosh s}{\sinh s}) (Iv(s, \tau))^{4} s^{2} ds d\tau \lesssim C_{0}(A, \epsilon) N^{1 - 2 \epsilon}.
\end{equation}
Therefore, by local well - posedness (lemma $\ref{l9.0}$) there exists some open interval $J_{1}$ that contains $J_{0}$, such that

\begin{equation}\label{6.17}
\int_{J_{1}} \int_{0}^{\infty} \phi(s) (\frac{\cosh s}{\sinh s}) (Iv(s, \tau))^{4} s^{2} ds d\tau \leq C_{1} N^{1 - 2 \epsilon},
\end{equation}
and

\begin{equation}
\sup_{t \in J_{1}} E(Iv(t)) \leq 2 C_{0}(A, \epsilon) N^{1 - 2 \epsilon}.
\end{equation}
Therefore $J_{0}$ is both open and closed in $\mathbf{R}$, and since $J_{0}$ is non - empty, $J_{0} = \mathbf{R}$, and

\begin{equation}\label{6.18}
\int_{\mathbf{R}} \int_{0}^{\infty} \phi(s) (\frac{\cosh s}{\sinh s}) (Iv(s, \tau))^{4} s^{2} ds d\tau \lesssim C_{0}(A, \epsilon) N^{1 - 2 \epsilon}.
\end{equation}
Meanwhile, by proposition $\ref{p6.2}$, $(\ref{6.18})$, and the fact that $\frac{\cosh s}{\sinh s} \gtrsim 1$,

\begin{equation}\label{6.19}
\sum_{j \geq k_{0}} \| P_{j} v \|_{L_{\tau, x}^{4}(\mathbf{R} \times \mathbf{R}^{3})} \lesssim_{\epsilon} C_{0}(A, \epsilon),
\end{equation}
so

\begin{equation}\label{6.20}
\int_{\mathbf{R}} \int_{0}^{\infty} \phi(s) ((1 - I) v(s,\tau))^{4} s^{2} ds d\tau \lesssim_{\epsilon} C_{0}(A, \epsilon)^{4}.
\end{equation}
Therefore, we have proved

\begin{equation}\label{6.21}
\int_{\mathbf{R}} \int_{\mathbf{R}^{3}} \phi(s) v(s, \tau)^{4} s^{2} ds d\tau \lesssim_{A, \epsilon} 1.
\end{equation}
Now we can follow the argument in section $6.3$ of \cite{Shen} to complete the proof of theorem $\ref{t1.4}$. Let $\chi(x)$ be a smooth function,

\begin{equation}\label{6.22}
\chi(x) = \left\{
	\begin{array}{ll}
		1  & \mbox{if } |x| \geq 1 \\
		0 & \mbox{if } |x| \leq \frac{1}{2}
	\end{array}
\right.
\end{equation}
Fixing $\delta > 0$ to be a small, fixed constant, $R(A, \epsilon)$ may be chosen to be sufficiently large so that

\begin{equation}\label{6.23}
\| \chi(\frac{x}{R(A, \epsilon)}) u_{0} \|_{\dot{H}^{1/2}(\mathbf{R}^{3})} + \| \chi(\frac{x}{R(A, \epsilon)}) u_{1} \|_{\dot{H}^{-1/2}(\mathbf{R}^{3})} \leq \delta.
\end{equation}
Then by small data arguments,

\begin{equation}\label{6.24}
w_{tt} - \Delta w = -w^{3}, \hspace{5mm} w(0,x) = \chi u_{0}, \hspace{5mm} w_{t}(0,x) = \chi u_{1},
\end{equation}
has a solution

\begin{equation}\label{6.25}
\| w \|_{L_{t,x}^{4}(\mathbf{R} \times \mathbf{R}^{3})} \lesssim \delta.
\end{equation}
Also, by finite propagation speed this implies $u = w$ for $|x| \geq R + |t|$, and thus

\begin{equation}\label{6.26}
\int_{0}^{\infty} \int_{|x| \geq R + t} u(x,t)^{4} dx dt \lesssim \delta.
\end{equation}
Choose $t_{0} < 0$ in $(\ref{6.1})$ to satisfy $t_{0}^{2} > R^{2} + 1$. Then

\begin{equation}\label{6.27}
\{ (t, r) : t \geq 0, 0 \leq r \leq R + t \} \subset \Omega = \{ (x, t) : |x|^{2} < (t - t_{0})^{2} - 1, t > t_{0} \} = \tilde{\mathbf{T}}(\{ (y, \tau) : \tau > 0 \}).
\end{equation}
Now by the change of variables formula,

\begin{equation}\label{6.28}
dx dt = 4 \pi r^{2} dr dt = e^{4 \tau} (\frac{\sinh |y|}{|y|})^{2} dy d\tau,
\end{equation}
so

\begin{equation}\label{6.29}
\aligned
\int \int_{\Omega} u(x,t)^{4} dx dt = \int_{0}^{\infty} \int e^{4 \tau} (\frac{\sinh s}{s})^{2} (u(e^{\tau} \sinh s, t_{0} + e^{\tau} \cosh s))^{4} s^{2} ds d\tau \\
= \int_{0}^{\infty} \int_{0}^{\infty} (\frac{s}{\sinh s})^{2} (e^{\tau} \frac{\sinh s}{s} u(e^{\tau} \sinh s, t_{0} + e^{\tau} \cosh s))^{4} s^{2} ds d\tau \\
= \int_{0}^{\infty} \int_{0}^{\infty} \phi(s) v(s, \tau)^{4} s^{2} ds d\tau \lesssim_{A, \epsilon} 1.
\endaligned
\end{equation}
This finally proves

\begin{equation}\label{6.30}
\int_{0}^{\infty} \int_{\mathbf{R}^{3}} u(x,t)^{4} dx dt < \infty.
\end{equation}
Theorem $\ref{t1.4}$ then follows by time reversal symmetry. $\Box$\vspace{5mm}

\noindent It only remains to pay off the debt incurred in propositions $\ref{p6.1}$ - $\ref{p6.4}$ as well as a local well - posedness result. The proof of local well - posedness is relatively straightforward. Assume throughout that the $N$ in $(\ref{6.6})$ is large.

\begin{lemma}[First local result]\label{l9.0}
Suppose $E(Iv(0)) \lesssim N^{1 - 2 \epsilon}$. Then the conformal wave equation

\begin{equation}
v_{tt} - \Delta v = -\phi(x) v^{3},
\end{equation}
is locally well - posed on some interval.
\end{lemma}
\emph{Proof:} By Strichartz estimates, H{\"o}lder's inequality, and Bernstein's inequality, since $P_{> k_{0}} |\nabla|^{1/2} I \sim P_{> k_{0}} N^{1/2 - \epsilon} |\nabla|^{\epsilon}$,

\begin{equation}
\aligned
\| |\nabla|^{1/2} Iv \|_{L_{t,x}^{4}(J \times \mathbf{R}^{3})} + \| Iv \|_{L_{t}^{2} L_{x}^{\infty}(J \times \mathbf{R}^{3})} + \| Iv \|_{L_{t}^{\infty} L_{x}^{6}(J \times \mathbf{R}^{3})} \\ \lesssim \| \nabla Iv(0) \|_{L^{2}(\mathbf{R}^{3})} + \| Iv_{t}(0) \|_{L^{2}(\mathbf{R}^{3})}
+ |J|^{1/2} \| Iv \|_{L_{t}^{2} L_{x}^{\infty}(J \times \mathbf{R}^{3})} \| Iv \|_{L_{t}^{\infty} L_{x}^{6}(J \times \mathbf{R}^{3})}^{2} \\
+ N^{1/2 - \epsilon} \| |\nabla|^{\epsilon} (1 - I) v \|_{L_{t,x}^{4}(J \times \mathbf{R}^{3})} \| (1 - I) v \|_{L_{t,x}^{4}(J \times \mathbf{R}^{3})}^{2} \\
\lesssim N^{1/2 - \epsilon} + |J|^{1/2}  \| Iv \|_{L_{t}^{2} L_{x}^{\infty}(J \times \mathbf{R}^{3})} \| Iv \|_{L_{t}^{\infty} L_{x}^{6}(J \times \mathbf{R}^{3})}^{2} + N^{-1} \| |\nabla|^{1/2} Iv \|_{L_{t,x}^{4}(J \times \mathbf{R}^{3})}^{3}.
\endaligned
\end{equation}
Taking $|J| \leq \frac{1}{N^{2}}$ proves

\begin{equation}
\| |\nabla|^{1/2} Iv \|_{L_{t,x}^{4}(J \times \mathbf{R}^{3})} + \| Iv \|_{L_{t}^{2} L_{x}^{\infty}(J \times \mathbf{R}^{3})} + \| Iv \|_{L_{t}^{\infty} L_{x}^{6}(J \times \mathbf{R}^{3})} \lesssim N^{1/2 - \epsilon}.
\end{equation}
This gives local well - posedness. $\Box$

\begin{lemma}[Second local result]\label{l9.1}
Suppose that $\phi(x) = (\frac{|x|}{\sinh |x|})^{2}$, and $v$ solves the equation

\begin{equation}\label{9.1}
v_{tt} - \Delta v = -\phi(x) v^{3}, \hspace{5mm} E(Iv(0)) \lesssim N^{1 - 2 \epsilon},
\end{equation}
on the interval $J_{k}$ and

\begin{equation}\label{9.2}
\int_{J_{k}} \int \phi(x) (Iv(x,t))^{4} dx dt \leq \delta^{4},
\end{equation}
for some small $\delta > 0$. ($\delta$ may be independent of $N$). Then

\begin{equation}\label{9.3}
\| Iv \|_{L_{t}^{2} L_{x}^{\infty}(J_{k} \times \mathbf{R}^{3})} + \sup_{j} \mathcal S(j, Iv) \lesssim N^{1/2 - \epsilon}.
\end{equation}
\end{lemma}
\emph{Proof:} Suppose $J_{k} = [a_{k}, b_{k}]$. Following an argument similar to lemma $\ref{l9.0}$,

\begin{equation}\label{9.4}
\aligned
\| Iv \|_{L_{t}^{2} L_{x}^{\infty}(J_{k} \times \mathbf{R}^{3})} + \| |\nabla|^{1/2} Iv \|_{L_{t,x}^{4}(J_{k} \times \mathbf{R}^{3})} \lesssim \| \nabla Iv(a_{k}) \|_{L^{2}(\mathbf{R}^{3})} + \| Iv_{t}(a_{k}) \|_{L^{2}(\mathbf{R}^{3})} \\
+ \| Iv \|_{L_{t}^{2} L_{x}^{\infty}(J_{k} \times \mathbf{R}^{3})} \| \phi(x)^{1/2} Iv \|_{L_{t,x}^{4}(J_{k} \times \mathbf{R}^{3})}^{2} \\ + N^{1/2 - \epsilon} \| |\nabla|^{\epsilon} (1 - I) v \|_{L_{t,x}^{4}(J_{k} \times \mathbf{R}^{3})} \| (1 - I) v \|_{L_{t,x}^{4}(J_{k} \times \mathbf{R}^{3})}^{2} \\
\lesssim N^{1/2 - \epsilon} + \delta^{2} \| Iv \|_{L_{t}^{2} L_{x}^{\infty}(J_{k} \times \mathbf{R}^{3})} + N^{-1} \| |\nabla|^{1/2} Iv \|_{L_{t,x}^{4}(J_{k} \times \mathbf{R}^{3})}^{3}.
\endaligned
\end{equation}
The theorem then follows by the contraction mapping principle and theorem $\ref{t7.6}$. $\Box$

\begin{corollary}\label{c9.2}
Suppose $J$ is an interval on which

\begin{equation}\label{9.5}
\int_{J} \int \phi(x) (Iv(x, t))^{4} dx dt \leq C N^{2(1 - s)},
\end{equation}
and $E(Iv)(t) \leq C N^{2(1 - s)}$ for all $t \in J$. Then

\begin{equation}\label{9.6}
\| Iv \|_{L_{t}^{2} L_{x}^{\infty}(J \times \mathbf{R}^{3})} + \sup_{j} \mathcal S(j, Iv) \lesssim C^{1/2} N^{1 - 2 \epsilon}.
\end{equation}
\end{corollary}
\emph{Proof:} Partition $J$ into $\lesssim N^{1 - 2 \epsilon}$ subintervals $J_{k}$ such that

\begin{equation}\label{8.7}
\int_{J_{k}} \int \phi(x) (Iv(x,t))^{4} dx dt \leq \delta^{4}.
\end{equation}
Then apply the previous lemma. $\Box$

\section{Initial data}
Now it only remains to prove propositions $\ref{p6.1}$ - $\ref{p6.4}$. In this section we prove proposition $\ref{p6.1}$. In order to utilize $(\ref{6.1})$, it is first necessary to prove that the solution to $(\ref{1.1})$ is global, which follows directly from proving that the initial data $u_{0}$ and $u_{1}$ lie in the set of data prescribed in \cite{D1}. Throughout this section, $f \lesssim g $ denotes $f \leq C(A, \epsilon) g$.

\begin{lemma}\label{l8.1}
Suppose that for some $\epsilon > 0$, $u_{0}$, $u_{1}$ are radial functions with $u_{0} \in \dot{H}^{\frac{1}{2} + \epsilon}(\mathbf{R}^{3})$ and $u_{1} \in \dot{H}^{\epsilon - \frac{1}{2}}(\mathbf{R}^{3})$. Also suppose that $|x|^{2 \epsilon} u_{0} \in \dot{H}^{\frac{1}{2} + \epsilon}(\mathbf{R}^{3})$ and $|x|^{2 \epsilon} u_{1} \in \dot{H}^{\epsilon - \frac{1}{2}}(\mathbf{R}^{3})$. Then this implies

\begin{equation}\label{8.1}
u_{0} \in \dot{H}^{1/2 + \epsilon}(\mathbf{R}^{3}) \cap \dot{H}^{1/2 - \epsilon}(\mathbf{R}^{3}), \hspace{5mm} \text{and} \hspace{5mm} u_{1} \in \dot{H}^{-1/2 + \epsilon}(\mathbf{R}^{3}) \cap \dot{H}^{-1/2 - \epsilon}(\mathbf{R}^{3}).
\end{equation}
\end{lemma}
\emph{Proof:} For any $0 \leq s < 1$,

\begin{equation}\label{8.2}
g \in \dot{H}^{s}(\mathbf{R}^{3}) \Rightarrow |x|^{- 2\epsilon} g \in \dot{H}^{s - 2 \epsilon}(\mathbf{R}^{3}).
\end{equation}
When $s = 0$, $(\ref{8.2})$ follows from Hardy's inequality. When $0 < s < 1$ we adapt the proof of Hardy's inequality (see for example \cite{Tao}). Let $\psi(x) \in C_{0}^{\infty}(\mathbf{R}^{3})$, $\psi(x) = 1$ for $|x| \leq 1$, and $\psi(x) = 0$ for $|x| > 2$. Then let
 
 \begin{equation}\label{5.3}
 \chi_{j}(x) = \psi(2^{-j - 1} x) - \psi(2^{-j} x),
 \end{equation}
 and make the partition of unity, when $x \neq 0$,
 
 \begin{equation}\label{5.4}
 1 = \sum_{j} \chi_{j}(x) = \psi(2^{k_{0}} x) + \sum_{j \geq -k_{0}} \chi_{j}(x).
 \end{equation}
 \noindent $k_{0}$ could be any integer. Combining the Littlewood - Paley partition of unity with the spatial partition of unity,

\begin{equation}
|x|^{-2 \epsilon} g = \sum_{m} \sum_{l} |x|^{-2 \epsilon} \chi_{l}(x) (P_{m} g).
\end{equation}
By Bernstein's inequality,

\begin{equation}
\| |x|^{-2 \epsilon} \chi_{l}(x) (P_{\geq j} g) \|_{L_{x}^{2}(\mathbf{R}^{3})} \lesssim 2^{-2l \epsilon} \sum_{m \geq j} 2^{-ms} \| P_{m} g \|_{\dot{H}^{s}},
\end{equation}
so by Young's inequality,

\begin{equation}
\aligned
\sum_{j} 2^{2j(s - 2 \epsilon)} \| P_{j}((1 - \psi(2^{j} x)) (P_{\geq j} f)) \|_{L^{2}}^{2} \\
\lesssim \sum_{j} 2^{2j(s - 2 \epsilon)} (\sum_{l \geq -j} \sum_{m \geq j} 2^{-2l \epsilon} 2^{-ms} \| P_{m} g \|_{\dot{H}^{s}})^{2} \lesssim \sum_{m} \| P_{m} g \|_{\dot{H}^{s}}^{2} \lesssim \| g \|_{\dot{H}^{s}}^{2}.
\endaligned
\end{equation}
Meanwhile, by H{\"o}lder's inequality, the Sobolev embedding theorem, and Bernstein's inequality,

\begin{equation}
\aligned
\| P_{j} (|x|^{-2 \epsilon} \chi_{l}(x) (P_{\geq j} g)) \|_{L_{x}^{2}(\mathbf{R}^{3})} \lesssim 2^{j} \| |x|^{-2 \epsilon} \chi_{l}(x) (P_{\geq j} g) \|_{L_{x}^{6/5}(\mathbf{R}^{3})} \\ \lesssim 2^{j} 2^{l(1 - 2 \epsilon)} \sum_{m \geq j} 2^{-ms} \| P_{m} g \|_{\dot{H}^{s}},
\endaligned
\end{equation}
so again by Young's inequality,

\begin{equation}
\aligned
\sum_{j} 2^{2j(s - 2 \epsilon)} \| P_{j}(\psi(2^{j} x) (P_{\geq j} f)) \|_{L^{2}}^{2} \\
\lesssim \sum_{j} 2^{2j(s - 2 \epsilon)} (\sum_{l \geq -j} \sum_{m \geq j} 2^{j} 2^{-ms} 2^{l(1 - 2 \epsilon)} \| P_{m} g \|_{\dot{H}^{s}})^{2} \lesssim \| g \|_{\dot{H}^{s}}^{2}.
\endaligned
\end{equation}
Also, by Bernstein's inequality,

\begin{equation}
\aligned
\| P_{j} (|x|^{-2 \epsilon} \chi_{l}(x) (P_{\leq j} g)) \|_{L_{x}^{2}(\mathbf{R}^{3})} \lesssim 2^{-j} \| \nabla (|x|^{-2 \epsilon} \chi_{l}(x) (P_{\leq j} g)) \|_{L_{x}^{2}(\mathbf{R}^{3})} \\
\lesssim 2^{-j} \| \nabla (|x|^{-2 \epsilon} \chi_{l}(x)) \|_{L^{3}} \| P_{\leq j} g \|_{L^{6}} + 2^{-j} \| |x|^{-2 \epsilon} \chi_{l}(x) \|_{L^{\infty}} \| \nabla P_{\leq j} g \|_{L^{2}} \\ \lesssim 2^{-j} 2^{-2 \epsilon l} \sum_{m \leq j} 2^{m(1 - s)} \| P_{m} g \|_{\dot{H}^{s}}.
\endaligned
\end{equation}
Again by Young's inequality this implies

\begin{equation}
\aligned
\sum_{j} 2^{2j(s - 2 \epsilon)} \| P_{j}((1 - \psi(2^{j} x)) (P_{\leq j} f)) \|_{L^{2}}^{2} \\
\lesssim \sum_{j} 2^{2j(s - 2 \epsilon)} (\sum_{l \geq -j} \sum_{m \leq j} 2^{-j} 2^{-2 \epsilon l} 2^{m(1 - s)} \| P_{m} g \|_{\dot{H}^{s}})^{2} \lesssim \| g \|_{\dot{H}^{s}}^{2}.
\endaligned
\end{equation}
Finally, by H{\"o}lder's inequality plus the Sobolev embedding theorem,

\begin{equation}
\| |x|^{-2 \epsilon} \chi_{l}(x) P_{\leq j} g \|_{L^{2}(\mathbf{R}^{3})} \lesssim \| P_{\leq j} g \|_{L^{6}} 2^{l(1 - 2 \epsilon)} \lesssim 2^{l(1 - 2 \epsilon)} \sum_{m \leq j} 2^{m(1 - s)} \| P_{m} g \|_{\dot{H}^{s}},
\end{equation}
and by Young's inequality,

\begin{equation}
\aligned
\sum_{j} 2^{2j(s - 2 \epsilon)} \| P_{j}(\psi(2^{j} x) (P_{\leq j} f)) \|_{L^{2}}^{2} \\
\lesssim \sum_{j} 2^{2j(s - 2 \epsilon)} (\sum_{l \leq -j} \sum_{m \leq j} 2^{m(1 - s)} 2^{l(1 - 2 \epsilon)} \| P_{m} g \|_{\dot{H}^{s}})^{2} \lesssim \| g \|_{\dot{H}^{s}}^{2}.
\endaligned
\end{equation}
This proves $(\ref{8.2})$, which directly implies $u_{0} \in \dot{H}^{\frac{1}{2} - \epsilon}$. By duality this also implies $u_{1} \in \dot{H}^{-1/2 - \epsilon}(\mathbf{R}^{3})$. Indeed, suppose $g \in \dot{H}^{1/2 + \epsilon}(\mathbf{R}^{3})$. Then

\begin{equation}\label{8.3}
\aligned
\int u_{1}(x) g(x) dx = \int (|x|^{2 \epsilon} u_{1}(x)) (|x|^{-2 \epsilon} g(x)) dx \\ \lesssim \| |x|^{2 \epsilon} u_{1} \|_{\dot{H}^{-1/2 + \epsilon}(\mathbf{R}^{3})} \| |x|^{-2 \epsilon} g \|_{\dot{H}^{1/2 - \epsilon}(\mathbf{R}^{3})} \\ \lesssim \| |x|^{2 \epsilon} u_{1} \|_{\dot{H}^{-1/2 + \epsilon}(\mathbf{R}^{3})} \| g \|_{\dot{H}^{1/2 + \epsilon}(\mathbf{R}^{3})}.
\endaligned
\end{equation}
This proves the lemma. $\Box$\vspace{5mm}

\noindent Thus, by $\cite{D1}$ we know that the wave equation $(\ref{1.1})$ with initial data prescribed in $(\ref{1.8})$ and $(\ref{1.9})$ has a global solution. Now let $v$ be the conformal transformation of $u$,

\begin{equation}\label{8.4}
v(s, \tau) = e^{\tau} \frac{\sinh s}{s} u(e^{\tau} \sinh s, t_{0} + e^{\tau} \cosh s).
\end{equation}

\begin{proposition}[Initial data]\label{p8.2}
$v(s, 0) \in \dot{H}^{\frac{1}{2} + \epsilon}(\mathbf{R}^{3})$ and $\partial_{\tau} v(s, 0) \in \dot{H}^{\epsilon - \frac{1}{2}}(\mathbf{R}^{3})$.
\end{proposition}
\emph{Proof:} Choose $s_{0}$ so that $\cosh s_{0} + t_{0} = 0$. Make a linear - nonlinear decomposition of $u$,

\begin{equation}\label{8.8}
u(t, x) = S(t)(u_{0}, u_{1}) + (u(t,x) - S(t)(u_{0}, u_{1})) = u_{l}(t,x) + u_{nl}(t,x),
\end{equation}
and let

\begin{equation}\label{8.9}
v_{1}(s, \tau) = \frac{e^{\tau} \sinh s}{s} u_{l}(e^{\tau} \sinh s, e^{\tau} \cosh s + t_{0}),
\end{equation}
and

\begin{equation}\label{8.10}
v_{2}(s, \tau) = \frac{e^{\tau} \sinh s}{s} u_{nl}(e^{\tau} \sinh s, e^{\tau} \cosh s + t_{0}).
\end{equation}
Recalling the free evolution of the radial wave equation in $\mathbf{R}^{3}$, since $t_{0} + \cosh s_{0} = 0$ and $\sinh s_{0} - (t_{0} + \cosh s_{0}) = -e^{-s_{0}} - t_{0} > -1 - t_{0}$, if $s \geq s_{0}$,

\begin{equation}\label{8.11}
s v_{1}(s, \tau) = \frac{1}{2} [u_{0}(t_{0} + e^{\tau + s})(t_{0} + e^{\tau + s}) + u_{0}(-t_{0} - e^{\tau - s})(-t_{0} - e^{\tau - s})] + \frac{1}{2} \int_{-t_{0} - e^{\tau - s}}^{t_{0} + e^{\tau + s}} u_{1}(r) r dr.
\end{equation}
\noindent We will first show that

\begin{equation}\label{8.14}
\aligned
(1 - \psi(\frac{s}{s_{0}})) \frac{1}{s} (\partial_{s} + \partial_{\tau})(s v_{1})(s, \tau)|_{\tau = 0} \\ = (1 - \psi(\frac{s}{s_{0}})) \frac{1}{s} [u_{0}'(t_{0} + e^{s})(t_{0} + e^{s}) e^{s} + u_{0}(t_{0} + e^{s}) e^{s} + e^{s} u_{1}(t_{0} + e^{s})(t_{0} + e^{s})]
\endaligned
\end{equation}
lies in $\dot{H}^{-1/2 + \epsilon}(\mathbf{R}^{3})$. By duality it suffices to estimate

\begin{equation}
\int_{0}^{\infty} g(s) s (1 - \psi(\frac{s}{s_{0}})) [(\partial_{s} + \partial_{\tau})(s v_{1})(s, \tau)|_{\tau = 0}] ds,
\end{equation}
for some radial $g \in \dot{H}^{1/2 - \epsilon}(\mathbf{R}^{3})$. Now if $g \in \dot{H}^{\sigma}(\mathbf{R}^{3})$ is a radial function, $-1 \leq \sigma \leq 1$, then by direct calculation,

\begin{equation}\label{8.6}
|x| g(|x|) = C \int_{0}^{\infty} \sin(|x| r) \hat{g}(r) r dr,
\end{equation}
which can be extended to an odd function lying in $\dot{H}^{\sigma}(\mathbf{R})$.\vspace{5mm}

\noindent Next we will make use of a technical lemma.

\begin{lemma}\label{l8.3}
Suppose $\chi(x) = \chi(|x|)$ is a radial function, $\chi \in C_{0}^{\infty}(\mathbf{R}^{3})$, and $\chi(x)$ supported on $\frac{1}{2} \leq |x| \leq 2$. Then for any $-1 < \sigma < 1$,

\begin{equation}\label{8.12}
\| \chi(2^{-k} x) f \|_{\dot{H}^{\sigma}(\mathbf{R}^{3})} \lesssim 2^{-2k \epsilon} \| |x|^{2 \epsilon} f \|_{\dot{H}^{\sigma}(\mathbf{R}^{3})}.
\end{equation}
\end{lemma}
\emph{Proof:} First take $0 \leq \sigma \leq 1$. By the product rule in \cite{Taylor} and the Sobolev embedding theorem, if $\frac{1}{p} = \frac{1}{2} - \frac{\sigma}{3}$,

\begin{equation}\label{8.13}
\aligned
\| |\nabla|^{\sigma} \chi(2^{-k} x) f \|_{L^{2}(\mathbf{R}^{3})} \lesssim \| |\nabla|^{\sigma} (|x|^{2 \epsilon} f) \|_{L^{2}(\mathbf{R}^{3})} \| \chi(2^{-k} x) |x|^{-2 \epsilon} \|_{L^{\infty}(\mathbf{R}^{3})} \\
+ \| |x|^{2 \epsilon} f \|_{L^{p}(\mathbf{R}^{3})} \| \nabla (\chi(2^{-k} x) |x|^{-2 \epsilon}) \|_{L^{3}(\mathbf{R}^{3})} \lesssim 2^{-2k \epsilon}  \| |x|^{2 \epsilon} f \|_{\dot{H}^{\sigma}(\mathbf{R}^{3})}.
\endaligned
\end{equation}
\noindent The case when $-1 \leq \sigma < 0$ follows by duality. $\Box$\vspace{5mm}

\noindent Now make the decomposition

\begin{equation}\label{8.15}
1 - \psi(\frac{s}{s_{0}}) = \sum_{k \geq s_{0}}^{\infty} \chi_{k}(x),
\end{equation}
where each $\chi_{k}(x)$, $k \geq s_{0}$ is supported on $k \leq |x| \leq k + 2$.\vspace{5mm}

\noindent Now for $-\frac{1}{2} \leq \sigma \leq \frac{1}{2}$ we have the scaling identity,

\begin{equation}\label{8.17}
\| u(\lambda x) \|_{\dot{H}^{\sigma}(\mathbf{R})} = \lambda^{\sigma - 1/2} \| u \|_{\dot{H}^{\sigma}(\mathbf{R})}.
\end{equation}
Also by the product rule and Hardy's inequality,

\begin{equation}
\| |x|^{2 \epsilon} u_{0}'(x) \|_{\dot{H}^{-1/2 + \epsilon}} \lesssim \| |x|^{2 \epsilon} u_{0} \|_{\dot{H}^{1/2 + \epsilon}} + \| |x|^{2 \epsilon - 1} u_{0} \|_{\dot{H}^{-1/2 + \epsilon}}.
\end{equation}
Therefore, by lemma $\ref{l8.3}$, $(\ref{8.6})$, and $(\ref{8.17})$, 

\begin{equation}\label{8.18}
\aligned
\| \chi_{k}(s) u_{0}'(t_{0} + e^{s})(t_{0} + e^{s}) e^{s} \|_{\dot{H}^{-1/2 + \epsilon}(\mathbf{R})} + 
\| \chi_{k}(s) u_{0}(t_{0} + e^{s}) e^{s} \|_{\dot{H}^{-1/2 + \epsilon}(\mathbf{R})} \\
+ \| \chi_{k}(s) e^{s} u_{1}(t_{0} + e^{s})(t_{0} + e^{s}) \|_{\dot{H}^{-1/2 + \epsilon}(\mathbf{R})} \\ \lesssim e^{-\epsilon k} (\| |x|^{2 \epsilon} u_{0} \|_{\dot{H}^{1/2 + \epsilon}(\mathbf{R}^{3})} + \| |x|^{2 \epsilon} u_{1} \|_{\dot{H}^{-1/2 + \epsilon}(\mathbf{R}^{3})}).
\endaligned
\end{equation}
Therefore, by $(\ref{8.6})$,

\begin{equation}\label{8.19}
\aligned
\int_{0}^{\infty} \psi(\frac{s}{s_{0}}) g(s) [(\partial_{\tau} + \partial_{s}) v_{1}(s, \tau)|_{\tau = 0}] s ds \\ \lesssim \sum_{k \geq s_{0}} 2^{-\epsilon k} \| g \|_{\dot{H}^{1/2 - \epsilon}(\mathbf{R}^{3})} (\| |x|^{2 \epsilon} u_{0} \|_{\dot{H}^{1/2 + \epsilon}(\mathbf{R}^{3})} + \| |x|^{2 \epsilon} u_{1} \|_{\dot{H}^{-1/2 + \epsilon}(\mathbf{R}^{3})}) \\
\lesssim  \| g \|_{\dot{H}^{1/2 - \epsilon}(\mathbf{R}^{3})} (\| |x|^{2 \epsilon} u_{0} \|_{\dot{H}^{1/2 + \epsilon}(\mathbf{R}^{3})} + \| |x|^{2 \epsilon} u_{1} \|_{\dot{H}^{-1/2 + \epsilon}(\mathbf{R}^{3})}).
\endaligned
\end{equation}
Next,

\begin{equation}\label{8.20}
\aligned
(\partial_{s} - \partial_{\tau})(s v_{1}(s, \tau))|_{\tau = 0} = u_{0}'(-t_{0} - e^{ - s})(-t_{0} - e^{ - s}) e^{- s} \\ + u_{0}(-t_{0} - e^{- s}) e^{- s} + u_{1}(-t_{0} - e^{- s})(-t_{0} - e^{ - s}) e^{-s}.
\endaligned
\end{equation}
This time, for any $k \geq s_{0}$, by $(\ref{8.17})$,

\begin{equation}\label{8.22}
\aligned
\| \chi_{k}(s) u_{0}'(-t_{0} - e^{- s})(-t_{0} - e^{- s}) e^{- s} \|_{\dot{H}^{\epsilon - 1/2}(\mathbf{R})} + \| \chi_{k}(s) u_{0}(-t_{0} - e^{- s}) e^{- s} \|_{\dot{H}^{\epsilon - 1/2}(\mathbf{R})} \\ + \| e^{- s} u_{1}(-t_{0} - e^{- s})(-t_{0} - e^{- s}) \|_{\dot{H}^{\epsilon - 1/2}(\mathbf{R})} \lesssim e^{-\epsilon k} (\| u_{0} \|_{\dot{H}^{1/2 + \epsilon}(\mathbf{R}^{3})} + \| u_{1} \|_{\dot{H}^{\epsilon - 1/2}(\mathbf{R}^{3})}).
\endaligned
\end{equation}
Finally, by the radial Sobolev embedding theorem,

\begin{equation}
\aligned
\| |x| u_{1} \|_{L_{t,x}^{\infty}(\mathbf{R} \times \{ x : |x| > |t| + R \})} \lesssim \| u_{1} \|_{L_{t}^{\infty} \dot{H}^{1/2 + \epsilon}(\mathbf{R} \times \mathbf{R}^{3})} + \| u_{1} \|_{L_{t}^{\infty} \dot{H}^{1/2 - \epsilon}(\mathbf{R} \times \mathbf{R}^{3})} \\ \lesssim \| u(0) \|_{\dot{H}^{1/2 + \epsilon}(\mathbf{R}^{3})} + \| u(0) \|_{\dot{H}^{1/2 - \epsilon}(\mathbf{R}^{3})} + \| u_{t}(0) \|_{\dot{H}^{-1/2 + \epsilon}(\mathbf{R}^{3})} + \| u_{t}(0) \|_{\dot{H}^{-1/2 - \epsilon}(\mathbf{R}^{3})} \lesssim 1.
\endaligned
\end{equation}
Therefore, by $(\ref{8.9})$,

\begin{equation}\label{8.23}
|v_{1}(s, 0)| \lesssim \frac{1}{s^{2}},
\end{equation}
and

\begin{equation}\label{8.24}
\int_{s_{0}}^{\infty} |v_{1}(s, 0)|^{\frac{3}{2 - \epsilon}} s^{2} ds < \infty.
\end{equation}
Therefore, we have proved 

\begin{equation}\label{8.26}
\aligned
\| \partial_{s} ((1 - \psi(\frac{s}{s_{0}})) v_{1})(s, \tau) \|_{\dot{H}^{- 1/2 + \epsilon}(\mathbf{R}^{3})} + \| \partial_{\tau} ((1 - \psi(\frac{s}{s_{0}})) v_{1})(s, \tau) \|_{\dot{H}^{- 1/2 + \epsilon}(\mathbf{R}^{3})} \\ \lesssim \| u_{0} \|_{\dot{H}^{1/2 + \epsilon}(\mathbf{R}^{3})} + \| u_{1} \|_{\dot{H}^{-1/2 + \epsilon}(\mathbf{R}^{3})} + \| |x|^{2 \epsilon} u_{0} \|_{\dot{H}^{1/2 + \epsilon}(\mathbf{R}^{3})} + \| |x|^{2 \epsilon} u_{1} \|_{\dot{H}^{-1/2 + \epsilon}(\mathbf{R}^{3})}.
\endaligned
\end{equation}
\noindent Next let us consider $(1 - \psi(\frac{s}{s_{0}})) v_{2}(s, \tau)$. Now for $s \geq s_{0}$, $t_{0} + \cosh s > -1 - t_{0} > 0$, and $\sinh s - (t_{0} + \cosh s) \geq -t_{0} - 1 > R$.\vspace{5mm}

\noindent By $(\ref{6.25})$, small data arguments, finite propagation speed, and the radial Sobolev embedding theorem,

\begin{equation}
\| |x| u \|_{L_{t,x}^{\infty}(\mathbf{R} \times \{ x : |x| > |t| + R \})} \lesssim \| u_{0} \|_{\dot{H}^{1/2 + \epsilon}} + \| u_{0} \|_{\dot{H}^{1/2 - \epsilon}} + \| u_{1} \|_{\dot{H}^{-1/2 + \epsilon}} + \| u_{1} \|_{\dot{H}^{-1/2 - \epsilon}}.
\end{equation}

%\begin{equation}\label{8.35}
%\aligned
%\| |x| u_{1} \|_{L_{t,x}^{\infty}(\mathbf{R} \times \{ x : |x| > |t| + R \})} \lesssim \| u_{1} \|_{L_{t}^{\infty} \dot{H}^{1/2 + \epsilon}(\mathbf{R} \times \mathbf{R}^{3})} + \| u_{1} \|_{L_{t}^{\infty} \dot{H}^{1/2 - \epsilon}(\mathbf{R} \times \mathbf{R}^{3})} \\ \lesssim \| u(0) \|_{\dot{H}^{1/2 + \epsilon}(\mathbf{R}^{3})} + \| u(0) \|_{\dot{H}^{1/2 - \epsilon}(\mathbf{R}^{3})} + \| u_{t}(0) \|_{\dot{H}^{-1/2 + \epsilon}(\mathbf{R}^{3})} + \| u_{t}(0) \|_{\dot{H}^{-1/2 - \epsilon}(\mathbf{R}^{3})} \lesssim 1.
%\endaligned
%\end{equation}
\noindent Therefore,

\begin{equation}\label{8.36}
\| |x|^{1/2 - \epsilon/2} u^{3} \|_{L_{t}^{1} L_{x}^{2}(\mathbf{R} \times \{ x : |x| > R + |t| \})} \lesssim \int_{\mathbf{R}} (\int_{R + |t|} s^{1 - \epsilon} \frac{s^{2}}{s^{6}} ds )^{1/2} dt \lesssim \int_{\mathbf{R}} \frac{1}{(R + |t|)^{1 + \epsilon/2}} dt < \infty.
\end{equation}
By direct computation

\begin{equation}\label{8.37}
(\partial_{\tau} + \partial_{s}) \int_{-t_{0} - e^{\tau - s} + t}^{t_{0} + e^{s + \tau} - t} u^{3}(r) r dr|_{\tau = 0} = (t_{0} + e^{s} - t) u^{3}(t_{0} + e^{s} - t) e^{s}.
\end{equation}
By the Sobolev embedding theorem and duality, $L^{\frac{1}{1 - \epsilon}}(\mathbf{R}) \subset \dot{H}^{-\frac{1}{2} + \epsilon}(\mathbf{R})$, so

\begin{equation}\label{8.38}
 \int_{0}^{\infty} (\int_{s_{0}}^{\infty} |(t_{0} + e^{s} - t) u^{3}(t_{0} + e^{s} - t) e^{s}|^{\frac{1}{1 - \epsilon}} ds)^{1 - \epsilon} dt < \infty,
\end{equation}
which implies that $(\partial_{\tau} + \partial_{s}) v_{2}|_{\tau = 0} \in \dot{H}^{-1/2 + \epsilon}(\mathbf{R}^{3})$. Also,

\begin{equation}\label{8.39}
(\partial_{\tau} - \partial_{s}) \int_{-t_{0} - e^{\tau - s} + t}^{t_{0} + e^{\tau + s} - t} u^{3}(r) r dr|_{\tau = 0} = (- t_{0} - e^{- s} + t) u^{3}(-t_{0} - e^{- s} + t) e^{- s},
\end{equation}
and

\begin{equation}
\int_{0}^{\infty} (\int_{s_{0}}^{\infty} |(- t_{0} - e^{- s} + t) u^{3}(-t_{0} - e^{- s} + t) e^{- s}|^{\frac{1}{1 - \epsilon}} ds)^{1 - \epsilon} dt < \infty.
\end{equation}
Then again by the finite propagation speed, $|v_{2}(s, 0)| \lesssim \frac{1}{s}$ when $s \geq s_{0}$, which proves $\frac{1}{s} v_{2}(s, 0) \lesssim \frac{1}{s^{2}}$, and therefore

\begin{equation}
\| (1 - \psi(\frac{s}{s_{0}})) v(s, \tau) \|_{\dot{H}^{1/2 + \epsilon}(\mathbf{R}^{3})} + \| \partial_{\tau} (1 - (\psi(\frac{s}{s_{0}})) v(s, \tau))| \|_{\dot{H}^{-1/2 + \epsilon}(\mathbf{R}^{3})} \lesssim 1.
\end{equation}
\noindent Now we turn to $\psi(\frac{s}{s_{0}}) v(s, \tau)|_{\tau = 0}$. Making a different linear - nonlinear decomposition of $u$,

\begin{equation}\label{8.5}
\aligned
u(t,x) = S(t - 1 - t_{0})(u(-t_{0} - 1,x), u_{t}(-t_{0} - 1,x)) \\ + [u(t,x) - S(t - 1 - t_{0})(u(-t_{0} - 1, x), u_{t}(-t_{0} - 1,x))] = u_{l}(t,x)  + u_{nl}(t,x).
\endaligned
\end{equation}
Again by the fundamental solution to the linear wave equation for radial data, if $v_{1}$ is given by $(\ref{8.9})$ under the new $u_{l}(t,x)$,
 
\begin{equation}\label{8.7}
\aligned
s v_{1}(s, \tau) = \frac{1}{2} u(-1 - t_{0}, e^{\tau + s} - 1)(e^{\tau + s} - 1) + \frac{1}{2} u(-1 - t_{0}, 1 - e^{\tau - s}) (1 - e^{\tau - s}) \\ + \frac{1}{2} \int_{1 - e^{\tau - s}}^{e^{\tau + s} - 1} u_{t}(-1 - t_{0}, r) r dr.
\endaligned
\end{equation}
By the global well - posedness result of \cite{D1}, 

\begin{equation}\label{estimate}
\| u(-1 + t_{0}, x) \|_{\dot{H}^{1/2 + \epsilon}} + \| u_{t}(-1 + t_{0}, x) \|_{\dot{H}^{-1/2 + \epsilon}} \lesssim 1.
\end{equation}
Therefore, by direct computation,

\begin{equation}\label{8.16}
\aligned
(\partial_{\tau} + \partial_{s}) (s v_{1}(s, \tau))|_{\tau = 0} = \frac{1}{2} u'(-1 - t_{0}, e^{s} - 1) (e^{s} - 1) e^{s} \\ + \frac{1}{2} u(-1 - t_{0}, e^{s} - 1) e^{s}
+ \frac{1}{2} u_{t}(-1 - t_{0}, e^{s} - 1) (e^{s} - 1) e^{s},
\endaligned
\end{equation}
and

\begin{equation}\label{8.21}
\aligned
(\partial_{\tau} - \partial_{s}) (s v_{1}(s, \tau))|_{\tau = 0} = -\frac{1}{2} u'(-1 - t_{0}, 1 - e^{-s}) (1 - e^{-s}) e^{-s} \\ - \frac{1}{2} u(-1 - t_{0}, 1 - e^{-s}) e^{-s}
- \frac{1}{2} u_{t}(-1 - t_{0}, 1 - e^{-s}) (1 - e^{-s}) e^{-s}.
\endaligned
\end{equation}
Since $\psi(\frac{s}{s_{0}}) (\frac{e^{s} - 1}{s}) e^{s}$ and $ \psi(\frac{s}{s_{0}})(\frac{1 - e^{-s}}{s}) e^{-s}$ and all their derivatives are uniformly bounded above and below, by $(\ref{estimate})$ and Hardy's inequality,

\begin{equation}\label{8.25}
\| \frac{1}{s} (\partial_{\tau} + \partial_{s}) (s v_{1}(s, \tau))|_{\tau = 0} \|_{\dot{H}^{\epsilon - 1/2}(\mathbf{R}^{3})} + \| \frac{1}{s} (\partial_{\tau} - \partial_{s}) (s v_{1}(s, \tau))|_{\tau = 0} \|_{\dot{H}^{\epsilon - 1/2}(\mathbf{R}^{3})} \lesssim 1.
\end{equation}
Also,

\begin{equation}\label{8.25.1}
\aligned
\frac{1}{s} v_{1}(s, \tau)|_{\tau = 0} =  \frac{1}{2s^{2}} u(-1 - t_{0}, e^{s} - 1)(e^{s} - 1) \\ + \frac{1}{2s^{2}} u(-1 - t_{0}, 1 - e^{- s}) (1 - e^{- s}) + \frac{1}{2s^{2}} \int_{1 - e^{- s}}^{e^{s} - 1} u_{t}(-1 - t_{0}, r) r dr.
\endaligned
\end{equation}
By Hardy's inequality,

\begin{equation}\label{8.25.2}
\aligned
\| \frac{1}{2s} \psi(\frac{s}{s_{0}}) u(-1 - t_{0}, e^{s} - 1) (\frac{e^{s} - 1}{s}) \|_{\dot{H}^{-1/2 + \epsilon}(\mathbf{R}^{3})} \\ + \| \frac{1}{2s} (1 - \psi(\frac{s}{s_{0}})) u(-1 - t_{0}, 1 - e^{- s}) (\frac{1 - e^{- s}}{s}) \|_{\dot{H}^{-1/2 + \epsilon}(\mathbf{R}^{3})} \lesssim 1.
\endaligned
\end{equation}
Finally, for any $0 < \theta < 1$ we have a uniform bound

\begin{equation}\label{8.25.3}
\| u_{t}(-1 - t_{0}, \theta (e^{s} - 1) + (1 - \theta)(1 - e^{-s})) \frac{ \theta (e^{s} - 1) + (1 - \theta)(1 - e^{-s})}{s} \|_{\dot{H}^{-1/2 + \epsilon}(\mathbf{R}^{3})} \lesssim 1,
\end{equation}
and since $\frac{\sinh s}{s} \psi(\frac{s}{s_{0}})$ and all its derivatives are uniformly bounded, making a change of variables,

\begin{equation}\label{8.25.5}
\aligned
\| \frac{1}{s^{2}} \psi(\frac{s}{s_{0}}) \int_{1 - e^{-s}}^{e^{s} - 1} u_{t}(-1 - t_{0}, r) r dr \|_{\dot{H}^{-1/2 + \epsilon}(\mathbf{R}^{3})} \\
= \| \frac{\sinh s}{s} \psi(\frac{s}{s_{0}})  \int_{0}^{1} u_{t}(-1 - t_{0}, \theta(e^{s} - 1) + (1 - \theta)(1 - e^{-s})) \\ \times \frac{\theta (e^{s} - 1) + (1 - \theta)(1 - e^{-s})}{s} d\theta \|_{\dot{H}^{-1/2 + \epsilon}} \lesssim 1.
\endaligned
\end{equation}
Therefore,

\begin{equation}\label{8.25.4}
\| \partial_{\tau} (\psi(\frac{s}{s_{0}}) v_{1}(\tau, s))|_{\tau = 0} \|_{\dot{H}^{-1/2 + \epsilon}(\mathbf{R}^{3})} + \| \partial_{s} (\psi(\frac{s}{s_{0}}) v_{1}(\tau, s))|_{\tau = 0} \|_{\dot{H}^{-1/2 + \epsilon}(\mathbf{R}^{3})} \lesssim 1.
\end{equation}

\noindent It only remains to evaluate $\psi(\frac{s}{s_{0}}) v_{2}(s, \tau)$, where $v_{2}$ is given by $(\ref{8.10})$ with $u_{nl}$ given by $(\ref{8.5})$.\vspace{5mm}

\noindent Again, by theorem $\ref{t7.6}$ and global well - posedness (see \cite{D1}), $\| u \|_{L_{t,x}^{4}([t_{0}, 0] \times \mathbf{R}^{3})} < \infty$, so for any $R > 0$,

\begin{equation}\label{8.27}
\| P_{j} u \|_{L_{t,x}^{2}([t_{0}, 0] \times B_{R})} \lesssim 2^{-j(1/2 + \epsilon)} R^{1/2},
\end{equation}

\begin{equation}\label{8.28}
\| P_{j} u \|_{L_{t}^{2} L_{x}^{\infty}([t_{0}, 0] \times \mathbf{R}^{3})} \lesssim 2^{j(1/2 - \epsilon)},
\end{equation}
and by the radial Sobolev embedding theorem,

\begin{equation}\label{8.29}
\| P_{j} u \|_{L_{t,x}^{\infty}([t_{0}, 0] \times A_{R})} \lesssim \inf(R^{-1 + 4 \epsilon} 2^{3j \epsilon} \| u \|_{L_{t}^{\infty} \dot{H}^{1/2 + \epsilon}([t_{0}, 0] \times \mathbf{R}^{3})}, R^{-1 + 2 \epsilon} 2^{j \epsilon} \| u \|_{L_{t}^{\infty} \dot{H}^{1/2 + \epsilon}([t_{0}, 0] \times \mathbf{R}^{3})}).
\end{equation}
Therefore,

\begin{equation}\label{8.30}
\| |x|^{1/2 - 2 \epsilon} u^{3} \|_{L_{t}^{1} L_{x}^{2}([t_{0}, 0] \times \{ x : |x| \leq 2e^{s_{0}} + |t_{0}| \})} \lesssim 1.
\end{equation}
Now for any $t_{0} \leq t \leq e^{s_{0}}$,

\begin{equation}\label{8.31}
(\partial_{\tau} + \partial_{s}) \int_{t_{0} + e^{s + \tau} + (-1 + t_{0} + t)}^{-t_{0} - e^{-s + \tau} + (t - t_{0} + 1)} r u^{3}(t, r) dr|_{\tau = 0} = e^{s} u^{3}(t_{0} + e^{s} + (-1 + t_{0} + t), t) (t_{0} + e^{s} + (-1 + t_{0} + t)).
\end{equation}
Then by $(\ref{8.30})$,

\begin{equation}\label{8.32}
\int_{t_{0}}^{e^{s_{0}}} (\int_{0}^{2s_{0}} |e^{s} u^{3}(t_{0} + e^{s} + (-1 + t_{0} + t), t) (t_{0} + e^{s} + (-1 + t_{0} + t))|^{\frac{1}{1 - \epsilon}} ds)^{1 - \epsilon} dt \lesssim 1.
\end{equation}
Also,

\begin{equation}\label{8.33}
(\partial_{\tau} - \partial_{s}) \int_{t_{0} + e^{s + \tau} - (t_{0} - 1 - t)}^{-t_{0} - e^{-s + \tau} + (t_{0} - 1 - t)} r u^{3}(t, r) dr|_{\tau = 0} = e^{-s} u^{3}(-t_{0} + e^{-s} + (t_{0} - 1 - t), t) (-t_{0} + e^{-s} + (t_{0} - 1 - t)),
\end{equation}
so again by $(\ref{8.30})$,

\begin{equation}\label{8.34}
\int_{t_{0}}^{e^{s_{0}}} (\int_{0}^{2s_{0}} |e^{-s} u^{3}(-t_{0} + e^{-s} + (t_{0} - 1 - t), t) (-t_{0} + e^{-s} + (t_{0} - 1 - t))|^{\frac{1}{1 - \epsilon}} ds)^{1 - \epsilon} dt \lesssim 1.
\end{equation}
Finally, as in $(\ref{8.25.3})$ - $(\ref{8.25.4})$, since for any $0 < \theta < 1$,

\begin{equation}
\aligned
\int_{t_{0}}^{e^{s_{0}}} (\int_{0}^{2s_{0}} |e^{s} u^{3}(\theta(t_{0} + e^{s} + (-1 + t_{0} + t)) + (1 - \theta)(-t_{0} + e^{-s} + (t_{0} - 1 - t)) , t) \\ \times (\theta(t_{0} + e^{s} + (-1 + t_{0} + t)) + (1 - \theta)(-t_{0} + e^{-s} + (t_{0} - 1 - t))|^{\frac{1}{1 - \epsilon}} ds)^{1 - \epsilon} dt \lesssim 1,
\endaligned
\end{equation}
for any $0 < \theta < 1$, we have proved

\begin{equation}
\| \partial_{\tau} (\psi(\frac{s}{s_{0}}) v_{2}(s, \tau))|_{\tau = 0} \|_{\dot{H}^{-1/2 + \epsilon}(\mathbf{R}^{3})} + \| \partial_{s} (\psi(\frac{s}{s_{0}}) v_{2}(s, \tau))|_{\tau = 0} \|_{\dot{H}^{-1/2 + \epsilon}(\mathbf{R}^{3})} \lesssim 1.
\end{equation}
This finally proves proposition $\ref{p8.2}$. $\Box$

\section{Multi - linear estimates}
The proofs of propositions $\ref{p6.2}$ - $\ref{p6.4}$ will utilize several multi - linear estimates, which will be proved in this section. All implicit constants in this section may depend on the $\epsilon > 0$ in theorem $\ref{t1.4}$. To simplify notation let $s = \frac{1}{2} + \epsilon$. $I$ is the operator defined in $(\ref{6.6})$.

\begin{theorem}\label{t5.1}
Suppose $|M(x)| \lesssim e^{-c |x|}$ for some $c > 0$, and $f$, $g$, $h$ are radial functions. Then

\begin{equation}\label{5.2}
\aligned
\| I(M(x) fg) (P_{\leq k_{0}} h) \|_{L_{t,x}^{1}(J \times \mathbf{R}^{3})} + \| M(x) fg (P_{\leq k_{0}} h) \|_{L_{t,x}^{1}(J \times \mathbf{R}^{3})} \\ \lesssim 2^{k_{0}/2} |k_{0}| \| f \|_{l_{k_{0}}^{\infty} L_{t,x}^{2}(J \times \mathbf{R}^{3})} \| g \|_{l_{k_{0}}^{\infty} L_{t,x}^{2}(J \times \mathbf{R}^{3})} \| h \|_{L_{t}^{\infty} L_{x}^{2}(J \times \mathbf{R}^{3})}.
\endaligned
\end{equation}
Also,

\begin{equation}\label{5.2.1}
\aligned
\| I(M(x) fg) \frac{1}{|x|} \|_{L_{t,x}^{1}(J \times \{ x : |x| > 2^{k_{0} - 1} \})} + \| M(x) fg \frac{1}{|x|} \|_{L_{t,x}^{1}(J \times \{ x : |x| > 2^{k_{0} - 1} \})} \\ \lesssim |k_{0}| \| f \|_{l_{k_{0}}^{\infty} L_{t,x}^{2}(J \times \mathbf{R}^{3})} \| g \|_{l_{k_{0}}^{\infty} L_{t,x}^{2}(J \times \mathbf{R}^{3})}.
\endaligned
\end{equation}

\end{theorem}
 \emph{Proof:} Recall the spatial partition of unity in $(\ref{5.3})$ and $(\ref{5.4})$. By H{\"o}lder's inequality and the Sobolev embedding theorem,

\begin{equation}\label{5.8}
\aligned
\| \psi(2^{k_{0}} x) M(x) fg(P_{\leq k_{0}} h) \|_{L_{t,x}^{1}(J \times \mathbf{R}^{3})}  \lesssim \| f \|_{L_{t,x}^{2}(J \times B_{1 - k_{0}})} \| g \|_{L_{t,x}^{2}(J \times B_{1 - k_{0}})} \| P_{\leq k_{0}} h \|_{L_{t,x}^{\infty}(J \times \mathbf{R}^{3})} \\
\lesssim 2^{k_{0}/2} \| f \|_{l_{k_{0}}^{\infty} L_{t,x}^{2}(J \times \mathbf{R}^{3})} \| g \|_{l_{k_{0}}^{\infty} L_{t,x}^{2}(J \times \mathbf{R}^{3})} \| h \|_{L_{t}^{\infty} L_{x}^{2} (J \times \mathbf{R}^{3})}.
\endaligned
\end{equation}
Also, by the radial Sobolev embedding theorem,

\begin{equation}\label{5.9}
\aligned
\sup_{j \geq -k_{0}} \| \chi_{j}(x) fg (P_{\leq k_{0}} h) \|_{L_{t,x}^{1}(J \times \mathbf{R}^{3})}
\lesssim \sup_{j \geq -k_{0}} \| f \|_{L_{t,x}^{2}(J \times A_{j})} \| g \|_{L_{t,x}^{2}(J \times A_{j})} \| P_{\leq k_{0}} h \|_{L_{t,x}^{\infty}(J \times A_{j})}
 \\ \lesssim 2^{k_{0}/2} \| f \|_{l_{k_{0}}^{\infty} L_{t,x}^{2}(J \times \mathbf{R}^{3})} \| g \|_{l_{k_{0}}^{\infty} L_{t,x}^{2}(J \times \mathbf{R}^{3})} \| h \|_{L_{t}^{\infty} L_{x}^{2}(J \times \mathbf{R}^{3})}.
%\lesssim (\sup_{R > \frac{1}{N}} R^{-1/2} \| f \|_{L_{t,x}^{2}(J \times \{ x : |x| \leq R \})}) (\sup_{R > \frac{1}{N}} R^{-1/2} \| g \|_{L_{t,x}^{2}(J \times \{ x : |x| \leq R \})}) \| h \|_{L_{t}^{\infty} \dot{H}^{1/2} (J \times \mathbf{R}^{3})}.
\endaligned
\end{equation} 
 Therefore, since $|M(x)| \lesssim e^{-c |x|}$,
 
\begin{equation}\label{5.9.1}
\aligned
\sum_{j \geq -k_{0}} \| \chi_{j}(x) M(x) fg (P_{\leq k_{0}} h) \|_{L_{t,x}^{1}(J \times \mathbf{R}^{3})} \\ \lesssim |k_{0}| 2^{k_{0}/2} \| f \|_{l_{k_{0}}^{\infty} L_{t,x}^{2}(J \times \mathbf{R}^{3})} \| g \|_{l_{k_{0}}^{\infty} L_{t,x}^{2}(J \times \mathbf{R}^{3})} \| h \|_{L_{t}^{\infty} L_{x}^{2}(J \times \mathbf{R}^{3})}.
\endaligned
\end{equation}
$\| I(M(x) fg) (P_{\leq k_{0}} h) \|_{L_{t,x}^{1}}$ may be estimated with a similar argument combined with the fact that the Littlewood - Paley kernels are rapidly decreasing. Indeed, by $(\ref{7.4})$, if $I(x)$ is the kernel of $I$, then for any $M > 0$,

% The kernels of $P_{j}$ and $P_{\leq j}$ are rapidly decreasing outside $|x| \geq 2^{-j}$, that is, for any $k$, if $K_{j}(x)$ is the kernel of $P_{j}$ or $P_{\leq j}$,
 
 %\begin{equation}\label{5.5}
% K_{j}(x) \lesssim_{k} \frac{2^{3j}}{(1 + 2^{j} |x|)^{k}}.
% \end{equation}
% Therefore, if $I(x)$ is the kernel of $I$, then
 
 \begin{equation}\label{5.6}
 I(x) \lesssim_{M} \frac{1}{(2^{k_{0}} |x|)^{M}}.
 \end{equation}
 So if $j, l \geq -k_{0}$, then by the radial Sobolev embedding theorem,
 
 \begin{equation}\label{5.7}
 \aligned
 \| \chi_{j}(x) I(\chi_{l}(x) M(x) fg) (P_{\leq k_{0}} h) \|_{L_{t,x}^{1}(J \times \mathbf{R}^{3})} \\ \lesssim 2^{-10 |j - l|} \| \chi_{l}(x) M(x) fg \|_{L_{t,x}^{1}(J \times \mathbf{R}^{3})} \| \chi_{j}(x) (P_{\leq k_{0}} h) \|_{L_{t,x}^{\infty}(J \times \mathbf{R}^{3})} \\
 \lesssim 2^{k_{0}/2} 2^{-10 |j - l|} 2^{l} 2^{-j} e^{-c 2^{l}} \| f \|_{l_{k_{0}}^{\infty} L_{t,x}^{2}(J \times \mathbf{R}^{3})} \| g \|_{l_{k_{0}} L_{t,x}^{2}(J \times \mathbf{R}^{3})} \| h \|_{L_{t}^{\infty} L_{x}^{2}(J \times \mathbf{R}^{3})}.
 \endaligned
 \end{equation}
 
 %Then making a partition of unity,
 
 %\begin{equation}\label{5.7}
 %\aligned
 %I(Mfg) (P_{\leq N} h) = I(\psi(Nx) M f g) (P_{\leq N} h) + \sum_{j > 0} I(\chi_{j}(x) M fg)(P_{\leq N} h) \\
% = \psi(Nx) I(\psi(Nx) Mf g)(P_{\leq N} h) + \sum_{l > 0} \chi_{l}(x) I(\psi(Nx) Mfg)(P_{\leq N} h) \\
% + \sum_{j > 0} \psi(Nx) I(\chi_{j}(x) Mf g)(P_{\leq N} h) + \sum_{j, l > 0} \chi_{l}(x) I(\chi_{j}(x) Mfg)(P_{\leq N} h).
 %\endaligned
 %\end{equation}
 
 \noindent Therefore, by $(\ref{5.8})$, $(\ref{5.9.1})$, $(\ref{5.7})$, and Young's inequality, the proof of $(\ref{5.2})$ is complete.\vspace{5mm}

\noindent The proof of $(\ref{5.2.1})$ is similar. Because $|M(x)| \lesssim e^{-c |x|}$,

\begin{equation}\label{5.13}
 \aligned
\| (1 - \psi(2^{k_{0}} x)) \frac{1}{|x|} M(x) fg \|_{L_{t,x}^{1}(J \times \mathbf{R}^{3})} \leq \sum_{j \geq -k_{0}} 2^{-j} \| \chi_{j}(x) M(x) fg \|_{L_{t,x}^{1}(J \times \mathbf{R}^{3})} \\ \lesssim |k_{0}| \cdot \sup_{j \geq -k_{0}} 2^{-j} \| f \|_{L_{t,x}^{2}(J \times A_{j})} \| g \|_{L_{t,x}^{2})(J \times A_{j})} \lesssim |k_{0}| \| f \|_{l_{k_{0}}^{\infty} L_{t,x}^{2}(J \times \mathbf{R}^{3})} \| g \|_{l_{k_{0}}^{\infty} L_{t,x}^{2}(J \times \mathbf{R}^{3})}.
\endaligned
 \end{equation}
 Also, since $\frac{1}{|x|} \leq 2^{k_{0}}$ when $|x| \geq 2^{-k_{0}}$,

\begin{equation}\label{5.14}
\aligned
\| (1 - \psi(2^{k_{0}} x)) \frac{1}{|x|} I(\psi(2^{k_{0}} x) M(x) fg) \|_{L_{t,x}^{1}(J \times \mathbf{R}^{3})} \lesssim 2^{k_{0}} \| f \|_{L_{t,x}^{2}(J \times B_{1 - k_{0}})} \| g \|_{L_{t,x}^{2}(J \times B_{1 - k_{0}})} \\ \lesssim \| f \|_{l_{k_{0}}^{\infty} L_{t,x}^{2}(J \times \mathbf{R}^{3})} \| g \|_{l_{k_{0}}^{\infty} L_{t,x}^{2}(J \times \mathbf{R}^{3})}. 
\endaligned
\end{equation}
Also, by $(\ref{5.6})$, $(\ref{5.13})$, and Young's inequality,

\begin{equation}\label{5.15}
\aligned
\sum_{j, l \geq -k_{0}} \| \chi_{l}(x) \frac{1}{|x|} I(\chi_{j}(x) M(x) fg) \|_{L_{t,x}^{1}(J \times \mathbf{R}^{3})}
\lesssim \sum_{j, l \geq -k_{0}} 2^{-10 |j - l|} 2^{-l} \| \chi_{j}(x) M(x) fg \|_{L_{t,x}^{1}(J \times \mathbf{R}^{3})} \\
\lesssim |k_{0}| \| f \|_{l_{k_{0}}^{\infty} L_{t,x}^{2}(J \times \mathbf{R}^{3})} \| g \|_{l_{k_{0}}^{\infty} L_{t,x}^{2}(J \times \mathbf{R}^{3})}.
\endaligned
\end{equation}

\noindent This concludes the proof of $(\ref{5.2.1})$. $\Box$

\begin{theorem}\label{t5.2}
Suppose that $M(x)$ is some function satisfying $|M(x)| \lesssim e^{-c |x|}$ for some $c > 0$. Then

\begin{equation}\label{5.16}
\aligned
\| I(M(x) g h) f \|_{L_{t,x}^{1}(J \times \mathbf{R}^{3})} \lesssim |k_{0}| \| g \|_{L_{t,x}^{2}(J \times \mathbf{R}^{3})} \| f \|_{l_{k_{0}}^{\infty} L_{t,x}^{2}(J \times \mathbf{R}^{3})} \\ \times  (2^{-k_{0}/2} \| h \|_{L_{t,x}^{\infty}(J \times \mathbf{R}^{3})} + \| |x|^{1/2} h \|_{L_{t,x}^{\infty}(J \times \mathbf{R}^{3})}).
\endaligned
\end{equation}
\end{theorem}
\emph{Proof:} Again by $(\ref{5.4})$,

\begin{equation}\label{5.17}
\aligned
 I(M(x) gh) f = \psi(2^{k_{0}} x) I(\psi(2^{k_{0}} x) M(x) g h) f + \sum_{l \geq -k_{0}} \chi_{l}(x) I(\psi(2^{k_{0}} x) M(x) gh) f \\
 + \sum_{j \geq -k_{0}} \psi(2^{k_{0}} x) I(\chi_{j}(x) M(x) gh)f + \sum_{j, l \geq -k_{0}} \chi_{l}(x) I(\chi_{j}(x) M(x) gh) f.
 \endaligned
\end{equation}
By H{\"o}lder's inequality and the support properties of $\psi(Nx)$,

\begin{equation}\label{5.18}
\aligned
\| \psi(Nx) I(\psi(2^{k_{0}} x) M(x) gh) f \|_{L_{t,x}^{1}(J \times \mathbf{R}^{3})} \lesssim \| f \|_{L_{t,x}^{2}(J \times B_{1 - k_{0}})} \| g \|_{L_{t,x}^{2}(J \times \mathbf{R}^{3})} \| h \|_{L_{t,x}^{\infty}(J \times \mathbf{R}^{3})} \\
\lesssim 2^{-k_{0}/2} \| f \|_{l_{k_{0}}^{\infty} L_{t,x}^{2}(J \times \mathbf{R}^{3})} \| g \|_{L_{t,x}^{2}(J \times \mathbf{R}^{3})} \| h \|_{L_{t,x}^{\infty}(J \times \mathbf{R}^{3})}.
\endaligned
\end{equation}
Also,

\begin{equation}\label{5.19}
\aligned
\sum_{j \geq -k_{0}} \| \psi(2^{k_{0}} x) I(\chi_{j}(x) M(x) gh) f \|_{L_{t,x}^{1}(J \times \mathbf{R}^{3})} \\ \lesssim \sum_{j \geq -k_{0}} \| f \|_{L_{t,x}^{2}(J \times B_{1 - k_{0}} )} \| g \|_{L_{t,x}^{2}(J \times \mathbf{R}^{3})} \| \chi_{j}(x) h \|_{L_{t,x}^{\infty}(J \times \mathbf{R}^{3})} \\
\lesssim 2^{-k_{0}/2} \| f \|_{l_{k_{0}}^{\infty} L_{t,x}^{2}(J \times \mathbf{R}^{3})} \| g \|_{L_{t,x}^{2}(J \times \mathbf{R}^{3})} \sum_{j \geq -k_{0}} \| \chi_{j}(x) h \|_{L_{t,x}^{\infty}(J \times \mathbf{R}^{3})} \\
\lesssim \| f \|_{l_{k_{0}}^{\infty} L_{t,x}^{2}(J \times \mathbf{R}^{3})} \| g \|_{L_{t,x}^{2}(J \times \mathbf{R}^{3})} \| |x|^{1/2} h \|_{L_{t,x}^{\infty}(J \times \mathbf{R}^{3})}.
\endaligned
\end{equation}
Next, by $(\ref{5.6})$,

\begin{equation}\label{5.20}
\aligned
\sum_{l \geq -k_{0}} \| \chi_{l}(x) I(\psi(2^{k_{0}} x) M(x) gh) f \|_{L_{t,x}^{1}(J \times \mathbf{R}^{3})} \\ \lesssim \sum_{l \geq -k_{0}} 2^{-10 (l + k_{0})} \| f \|_{L_{t,x}^{2}(J \times A_{l})} \| g \|_{L_{t,x}^{2}(J \times \mathbf{R}^{3})} \| h \|_{L_{t,x}^{\infty}(J \times \mathbf{R}^{3})} \\
\lesssim \sum_{l \geq -k_{0}} 2^{-10(l + k_{0})} 2^{l/2}  \| f \|_{l_{k_{0}}^{\infty} L_{t,x}^{2}(J \times \mathbf{R}^{3})} \| g \|_{L_{t,x}^{2}(J \times \mathbf{R}^{3})} \| h \|_{L_{t,x}^{\infty}(J \times \mathbf{R}^{3})} \\
\lesssim 2^{-k_{0}/2} \| f \|_{l_{k_{0}}^{\infty} L_{t,x}^{2}(J \times \mathbf{R}^{3})} \| g \|_{L_{t,x}^{2}(J \times \mathbf{R}^{3})} \| h \|_{L_{t,x}^{\infty}(J \times \mathbf{R}^{3})}.
\endaligned
\end{equation}
Finally by $(\ref{5.6})$ and $|M(x)| \lesssim e^{-c |x|}$,

\begin{equation}\label{5.21}
\aligned
\sum_{j, l \geq -k_{0}} \| \chi_{l}(x) I(\chi_{j}(x) M(x) gh) f \|_{L_{t,x}^{1}(J \times \mathbf{R}^{3})} \\
\lesssim \sum_{j, l > 0} 2^{-10 |j - l|} \| f \|_{L_{t,x}^{2}(J \times A_{l})} \| g \|_{L_{t,x}^{2}(J \times \mathbf{R}^{3})} \| \chi_{j}(x) M(x) h \|_{L_{t,x}^{\infty}(J \times \mathbf{R}^{3})} \\
\lesssim \sum_{j, l > 0} 2^{-10 |j - l|} 2^{l/2 - j/2} e^{-c 2^{j}} \| f \|_{l_{k_{0}}^{\infty} L_{t,x}^{2}(J \times \mathbf{R}^{3})} \| g \|_{L_{t,x}^{2}(J \times \mathbf{R}^{3})} \| |x|^{1/2} h \|_{L_{t,x}^{\infty}(J \times \mathbf{R}^{3})} \\
\lesssim |k_{0}| \| f \|_{l_{k_{0}}^{\infty} L_{t,x}^{2}(J \times \mathbf{R}^{3})} \| g \|_{L_{t,x}^{2}(J \times \mathbf{R}^{3})} \| |x|^{1/2} h \|_{L_{t,x}^{\infty}(J \times \mathbf{R}^{3})}.
\endaligned
\end{equation}
This proves the theorem. $\Box$

\begin{theorem}\label{t5.3}
If $M(x)$ is a function satisfying $|M(x)| \lesssim e^{-c |x|}$, then

\begin{equation}\label{5.22}
\aligned
\| I (M(x) (P_{\geq k_{0} - 7} f)^{3}) \|_{L_{t}^{1} L_{x}^{2}(J \times \mathbf{R}^{3})} \lesssim 2^{k_{0}(1 - s - \epsilon)} (\sup_{j > k_{0} - 7} 2^{-j(1 - s)} \mathcal S(j, f))^{3} \\ + 2^{k_{0}(1 - s - \epsilon)} (\sup_{j > k_{0} - 7} 2^{-j(1 - s)} \mathcal S(j, f))^{2} \| f \|_{L_{t}^{\infty} \dot{H}^{s}(J \times \mathbf{R}^{3})} .
% \\ \lesssim N^{1 - s} N^{-\epsilon} (\sup_{N_{1} \geq \frac{N}{100}} N_{1}^{\epsilon} \| P_{N_{1}} f \|_{L_{t,x}^{4}(J \times \mathbf{R}^{3})})^{2}  (\sup_{N_{1} \geq \frac{N}{100}} N_{1}^{s - 1} \| P_{N_{1}} f \|_{L_{t}^{2} L_{x}^{\infty}(J \times \mathbf{R}^{3})}) \\
%+ N^{1 - s} N^{-\epsilon} (\sup_{N_{1} \geq \frac{N}{100}} N_{1}^{s} (\sup_{R > 0} R^{-1/2} \| P_{N_{1}} f \|_{L_{t,x}^{2}(J \times \{ x : |x| \leq R \})})^{2} \| f \|_{L_{t}^{\infty} \dot{H}^{s}(J \times \mathbf{R}^{3})}.
\endaligned
\end{equation}
\end{theorem}
\emph{Proof:} By the Sobolev embedding theorem,

\begin{equation}\label{5.23}
\| P_{\leq k_{0}} (M(x) (P_{\geq k_{0} - 7} f)^{3}) \|_{L_{t}^{1} L_{x}^{2}(J \times \mathbf{R}^{3})} \lesssim 2^{3k_{0}/2} \| M(x) (P_{\geq k_{0} - 7}f)^{3} \|_{L_{t,x}^{1}(J \times \mathbf{R}^{3})}
\end{equation}

\begin{equation}\label{5.24}
\lesssim 2^{3k_{0}/2} \sum_{k_{0} - 7 \leq j_{1} \leq j_{2} \leq j_{3}} \| M(x) (P_{j_{1}} f)(P_{j_{2}} f)(P_{j_{3}} f) \|_{L_{t,x}^{1}(J \times \mathbf{R}^{3})}.
\end{equation}
Now, by the radial Sobolev embedding theorem $\| |x|^{1 - \epsilon} f \|_{L^{\infty}} \lesssim \| f \|_{\dot{H}^{s}}$ and Bernstein's inequality, for any $l \in \mathbf{Z}$, 

\begin{equation}\label{5.25}
\aligned
\| (P_{j_{1}} f)(P_{j_{2}} f)(P_{j_{3}} f) \|_{L_{t,x}^{1}(J \times A_{l})} \lesssim \| P_{j_{1}} f \|_{L_{t,x}^{\infty}(J \times A_{l})} \| P_{j_{2}} f \|_{L_{t,x}^{2}(J \times A_{l})} \| P_{j_{3}} f \|_{L_{t,x}^{2}(J \times A_{l})} \\
\lesssim 2^{l \epsilon} \| P_{j_{1}} f \|_{L_{t}^{\infty} \dot{H}^{s}(J \times \mathbf{R}^{3})} \| P_{j_{2}} f \|_{l^{\infty} L_{t,x}^{2}(J \times \mathbf{R}^{3})}) \| P_{j_{3}} f \|_{l^{\infty} L_{t,x}^{2}(J \times \mathbf{R}^{3})}.
\endaligned
\end{equation}
Now by H{\"o}lder's inequality, for any $l \in \mathbf{Z}$,

\begin{equation}
2^{-l/2} \| f \|_{L_{t,x}^{2}(J \times B_{l})} \lesssim 2^{l} \| f \|_{L_{t}^{2} L_{x}^{\infty}(J \times \mathbf{R}^{3})},
\end{equation}
so in particular, when $l \leq -j$,

\begin{equation}\label{L2estimate}
2^{-l/2} \| P_{j} f \|_{L_{t,x}^{2}(J \times B_{l})} \lesssim 2^{-j} \mathcal S(j, f).
\end{equation}
Combining this fact with $(\ref{5.25})$ and the fact that $|M(x)| \lesssim e^{-cR/2}$ for $\frac{R}{2} \leq |x| \leq R$,

\begin{equation}\label{5.26}
\aligned
(\ref{5.24}) \lesssim 2^{3k_{0}/2} \sum_{k_{0} - 7 \leq j_{1} \leq j_{2} \leq j_{3}} \| P_{j_{1}} f \|_{L_{t}^{\infty} \dot{H}^{s}(J \times \mathbf{R}^{3})} \| P_{j_{2}} f \|_{l^{\infty} L_{t,x}^{2}(J \times \mathbf{R}^{3})}  \| P_{j_{3}} f \|_{l^{\infty} L_{t,x}^{2}(J \times \mathbf{R}^{3})} \\
\lesssim (\sup_{j > k_{0} - 7} 2^{-j(1 - s)} \mathcal S(j, f))^{2} \| f \|_{L_{t}^{\infty} \dot{H}^{s}(J \times \mathbf{R}^{3})} \cdot 2^{3k_{0}/2} \sum_{k_{0} - 7 \leq j_{1} \leq j_{2} \leq j_{3}} 2^{-s j_{2}} 2^{-s j_{3}} \\
\lesssim 2^{k_{0}(1 - s - \epsilon)} (\sup_{j > k_{0} - 7} 2^{-j(1 - s)} \mathcal S(j, f))^{2} \| f \|_{L_{t}^{\infty} \dot{H}^{s}(J \times \mathbf{R}^{3})}.
\endaligned
\end{equation}
When estimating

\begin{equation}\label{5.27}
\sum_{j > k_{0}} 2^{k_{0}(1 - s)} 2^{-j(1 - s)} \| P_{j} (M(x) (P_{\geq k_{0} - 7} f)^{3}) \|_{L_{t}^{1} L_{x}^{2}(J \times \mathbf{R}^{3})},
\end{equation}
the terms in which $j_{2} \leq j$ and the terms in which $j_{2} \geq j$ will be analyzed in two different manners, where once again $j_{1} \leq j_{2} \leq j_{3}$.\vspace{5mm}

\noindent Again by the Sobolev embedding theorem, $(\ref{5.25})$, $(\ref{L2estimate})$, and $|M(x)| \lesssim e^{-c |x|}$,

\begin{equation}\label{5.28}
\aligned
\sum_{k_{0} < j \leq j_{2}} \sum_{k_{0} - 7 \leq j_{1} \leq j_{2} \leq j_{3}} 2^{k_{0}(1 - s)} 2^{-j(1 - s)} \| P_{j} (M(x) (P_{j_{1}} f)(P_{j_{2}} f)(P_{j_{3}} f)) \|_{L_{t}^{1} L_{x}^{2}(J \times \mathbf{R}^{3})} \\
\lesssim \sum_{k_{0} < j \leq j_{2}} \sum_{k_{0} - 7 \leq j_{1} \leq j_{2} \leq j_{3}} 2^{k_{0}(1 - s)} 2^{j(1/2 + s)} \| P_{j_{1}} f \|_{L_{t}^{\infty} \dot{H}^{s}(J \times \mathbf{R}^{3})} \| P_{j_{2}} f \|_{l^{\infty} L_{t,x}^{2}(J \times \mathbf{R}^{3})} \| P_{j_{3}} f \|_{l^{\infty} L_{t,x}^{2}(J \times \mathbf{R}^{3})} \\ \lesssim 2^{k_{0}(1 - s - \epsilon)} (\sup_{j > k_{0} - 7} 2^{-j(1 - s)} \mathcal S(j, f))^{2} \| f \|_{L_{t}^{\infty} \dot{H}^{s}(J \times \mathbf{R}^{3})}.
\endaligned
\end{equation}
Meanwhile,

\begin{equation}\label{5.29}
\aligned
\sum_{k_{0} - 7 \leq j_{2} \leq j} \sum_{k_{0} - 7 \leq j_{1} \leq j_{2} \leq j_{3}} 2^{k_{0}(1 - s)} 2^{-j(1 - s)} \| P_{j_{1}} f \|_{L_{t}^{2} L_{x}^{\infty}(J \times \mathbf{R}^{3})} \| P_{j_{2}} f \|_{L_{t,x}^{4}(J \times \mathbf{R}^{3})} \| P_{j_{3}} f \|_{L_{t,x}^{4}(J \times \mathbf{R}^{3})} \\
\lesssim  (\sup_{j_{2} \geq k_{0} - 7} 2^{\epsilon j_{2}} \| P_{j_{2}} f \|_{L_{t,x}^{4}(J \times \mathbf{R}^{3})})^{2} \\ \times \sum_{k_{0} - 7 \leq j_{2} \leq j} \sum_{k_{0} - 7 \leq j_{1} \leq j_{2}} 2^{-j(1 - s)} 2^{-2 \epsilon j_{2}} 2^{j_{1}(1 - s)} (\sup_{j_{1} > k_{0} - 7} 2^{-j_{1}(1 - s)} \| P_{j_{1}} f \|_{L_{t}^{2} L_{x}^{\infty}(J \times \mathbf{R}^{3})}) \\
\lesssim 2^{k_{0}(1 - s - 2\epsilon)} (\sup_{j > k_{0} - 7} 2^{-j(1 - s)} \mathcal S(j, f))^{3}.
\endaligned
\end{equation}
This proves the theorem. $\Box$

\begin{proposition}[Bilinear estimate]\label{p5.4}
Suppose $|M(x)| \lesssim e^{-c |x|}$. Then

\begin{equation}\label{5.30}
\| M(x) (P_{j} f) g \|_{L_{t,x}^{2}(J \times \mathbf{R}^{3})} \lesssim |k_{0}| 2^{-j} \mathcal S(j, f) (N^{-1/2} \| g \|_{L_{t,x}^{\infty}(J \times \mathbf{R}^{3})} + \| g \|_{L_{t}^{\infty} \dot{H}^{1}(J \times \mathbf{R}^{3})}).
\end{equation}
\end{proposition}
\emph{Proof:} First, by H{\"o}lder's inequality,

\begin{equation}\label{5.31}
\aligned
 \| (P_{j} f) g \|_{L_{t,x}^{2}(J \times B_{-j})} \lesssim 2^{-j} \| P_{j} f \|_{L_{t}^{2} L_{x}^{\infty}(J \times \mathbf{R}^{3})} \| g \|_{L_{t}^{\infty} L^{6}(J \times \mathbf{R}^{3})} \lesssim 2^{-j} \mathcal S(j, f) \| g \|_{L_{t}^{\infty} \dot{H}^{1}(J \times \mathbf{R}^{3})}.
 \endaligned
\end{equation}
Next, for any $-j \leq l \leq 1 - k_{0}$,

\begin{equation}
\| (P_{j} f) g \|_{L_{t,x}^{2}(J \times A_{l})} \lesssim 2^{l/2} \| P_{j} f \|_{l_{j}^{\infty} L_{t,x}^{2}(J \times \mathbf{R}^{3})} \| g \|_{L_{t,x}^{\infty}(J \times \mathbf{R}^{3})} \lesssim 2^{-j} 2^{l/2} \mathcal S(j, f) \| g \|_{L_{t,x}^{\infty}(J \times \mathbf{R}^{3})}.
\end{equation}
Summing up,

\begin{equation}
\sum_{-j \leq l \leq 1 - k_{0}} 2^{-j} 2^{l/2} \mathcal S(j, f) \| g \|_{L_{t,x}^{\infty}(J \times \mathbf{R}^{3})} \lesssim 2^{-k_{0}/2} 2^{-j} \mathcal S(j, f) \| g \|_{L_{t,x}^{\infty}(J \times \mathbf{R}^{3})}.
\end{equation}
Now by the radial Sobolev embedding theorem,

\begin{equation}\label{5.32}
\aligned
\sum_{-k_{0} \leq j \leq 0} \| fg \|_{L_{t,x}^{2}(J \times A_{j})}
\lesssim |k_{0}| \sup_{-k_{0} \leq j \leq 0} (2^{-j/2} \| f \|_{L_{t,x}^{2}(J \times A_{j})}) (2^{j/2} \| g \|_{L_{t,x}^{\infty}(J \times A_{j})}) \\ \lesssim |k_{0}| \| f \|_{l^{\infty} L_{t,x}^{2}(J \times \mathbf{R}^{3})} \| g \|_{L_{t}^{\infty} \dot{H}^{1}(J \times \mathbf{R}^{3})}.
\endaligned
\end{equation}
Finally, by the radial Sobolev embedding combined with the fact that $|M(x)| \lesssim e^{-c|x|}$, for some constant $c > 0$,

\begin{equation}\label{5.33}
\aligned
\sum_{j \geq 0} \| M(x) fg \|_{L_{t,x}^{2}(J \times A_{j})} \lesssim  \sum_{j \geq 0} e^{-c2^{j}} ( \sup_{j} 2^{-j/2} \| f \|_{L_{t,x}^{2}(J \times A_{j})}) (2^{j/2} \| g \|_{L_{t,x}^{\infty}(J \times A_{j})}) \\ \lesssim \| f \|_{l^{\infty} L_{t,x}^{2}(J \times \mathbf{R}^{3})} \| g \|_{L_{t}^{\infty} \dot{H}^{1}(J \times \mathbf{R}^{3})}.
\endaligned
\end{equation}
This proves the proposition. $\Box$\vspace{5mm}

\noindent We conclude this section with a tri - linear estimate close to the origin.

\begin{theorem}\label{t5.5}
If $|M(x)| \lesssim e^{-c |x|}$ for some $c > 0$ and $\psi(2^{k_{0}} x)$ is as in $(\ref{5.4})$,

\begin{equation}\label{5.34}
\aligned
\| I(M(x) fgh) \psi(2^{k_{0}} x) \frac{1}{|x|^{3/4}} \|_{L_{t,x}^{4/3}(J \times \mathbf{R}^{3})} \lesssim 2^{-k_{0} s} (\sup_{j} 2^{-j(1 - s)} \mathcal S(j, f)) \| g^{2} \|_{l_{k_{0}}^{\infty} L_{t,x}^{2}(J \times \mathbf{R}^{3})}^{1/2} \| h \|_{L_{t,x}^{\infty}(J \times \mathbf{R}^{3})}.
%\\ \lesssim N^{- s} (\sup_{M_{1}} M_{1}^{s - 1} \| P_{M_{1}} f \|_{L_{t}^{2} L_{x}^{\infty}(J \times \mathbf{R}^{3})}) (\sup_{R > \frac{1}{N}} R^{-1/4} \| g \|_{L_{t,x}^{4}(J \times \{ x : |x| \leq R \})}) \| h \|_{L_{t,x}^{\infty}(J \times \mathbf{R}^{3})} \\
%+ N^{-s} (\sup_{M_{1}} M_{1}^{s} (\sup_{R > 0} R^{-1/2} \| P_{M_{1}} f \|_{L_{t,x}^{2}(J \times \{ x : |x| \leq R \})})) (\sup_{R > \frac{1}{N}} R^{-1/4} \| g \|_{L_{t,x}^{4}(J \times \{ x : |x| \leq R \})}) \| h \|_{L_{t,x}^{\infty}(J \times \mathbf{R}^{3})}.
\endaligned
\end{equation}
\end{theorem}
\emph{Proof:} Combining the Littlewood - Paley decomposition, $(\ref{6.6})$, and $(\ref{5.4})$,

\begin{equation}\label{5.35}
\aligned
\| I(M(x) fgh) \frac{1}{|x|^{3/4}} \psi(2^{k_{0}} x) \|_{L_{t,x}^{4/3}(J \times \mathbf{R}^{3})} \lesssim \| P_{\leq k_{0}}(M(x) fgh) \frac{1}{|x|^{3/4}} \psi(2^{k_{0}} x) \|_{L_{t,x}^{4/3}(J \times \mathbf{R}^{3})} \\ + \sum_{j > k_{0}} 2^{k_{0}(1 - s)} 2^{-j(1 - s)} \| P_{j}(M(x) fgh) \frac{1}{|x|^{3/4}} \psi(2^{k_{0}} x) \|_{L_{t,x}^{4/3}(J \times \mathbf{R}^{3})}
\endaligned
\end{equation}

\begin{equation}\label{5.36}
\aligned
\lesssim \| P_{\leq k_{0}} (M(x) \psi(2^{k_{0}} x) fgh) \frac{1}{|x|^{3/4}} \psi(2^{k_{0}} x) \|_{L_{t,x}^{4/3}(J \times \mathbf{R}^{3})} \\ + \sum_{l \geq -k_{0}} \| P_{\leq k_{0}} (\chi_{l}(x) M(x) fgh) \frac{1}{|x|^{3/4}} \psi(2^{k_{0}} x) \|_{L_{t,x}^{4/3}(J \times \mathbf{R}^{3})} \\
+ \sum_{j > k_{0}} 2^{k_{0}(1 - s)} 2^{-j(1 - s)} \| P_{j}(M(x) \psi(2^{k_{0}} x) fgh) \frac{1}{|x|^{3/4}} \psi(2^{k_{0}} x) \|_{L_{t,x}^{4/3}(J \times \mathbf{R}^{3})} \\ + \sum_{j > k_{0}} \sum_{l \geq -k_{0}} 2^{k_{0}(1 - s)} 2^{-j(1 - s)} \| P_{j}(\chi_{l}(x) M(x) fgh) \frac{1}{|x|^{3/4}} \psi(2^{k_{0}} x) \|_{L_{t,x}^{4/3}(J \times \mathbf{R}^{3})}
\endaligned
\end{equation}

\begin{equation}\label{5.37}
\aligned
\lesssim \sum_{j_{2} \leq k_{0}} \| P_{\leq k_{0}} (M(x) \psi(2^{k_{0}} x) (P_{j_{2}} f) gh) \frac{1}{|x|^{3/4}} \psi(2^{k_{0}} x) \|_{L_{t,x}^{4/3}(J \times \mathbf{R}^{3})} \\
+ \sum_{j_{2} > k_{0}} \| P_{\leq k_{0}} (M(x) \psi(2^{k_{0}} x) (P_{j_{2}} f) gh) \frac{1}{|x|^{3/4}} \psi(2^{k_{0}} x) \|_{L_{t,x}^{4/3}(J \times \mathbf{R}^{3})}
\endaligned
\end{equation}

\begin{equation}\label{5.38}
\aligned
+ \sum_{l \geq -k_{0}} \sum_{j_{2} \leq k_{0}} \| P_{\leq k_{0}} (\chi_{l}(x) M(x) (P_{j_{2}} f) gh) \frac{1}{|x|^{3/4}} \psi(2^{k_{0}} x) \|_{L_{t,x}^{4/3}(J \times \mathbf{R}^{3})} \\ + \sum_{l \geq -k_{0}} \sum_{j_{2} > k_{0}} \| P_{\leq k_{0}} (\chi_{l}(x) M(x) (P_{j_{2}} f) gh) \frac{1}{|x|^{3/4}} \psi(2^{k_{0}} x) \|_{L_{t,x}^{4/3}(J \times \mathbf{R}^{3})}
\endaligned
\end{equation}

\begin{equation}\label{5.39}
\aligned
+ \sum_{j > k_{0}} \sum_{j_{2} \leq j} 2^{k_{0}(1 - s)} 2^{-j(1 - s)} \| P_{j}(M(x) \psi(2^{k_{0}} x) (P_{j_{2}} f)gh) \frac{1}{|x|^{3/4}} \psi(2^{k_{0}} x) \|_{L_{t,x}^{4/3}(J \times \mathbf{R}^{3})} \\
\sum_{j > k_{0}} \sum_{j_{2} > j} 2^{k_{0}(1 - s)} 2^{-j(1 - s)} \| P_{j}(M(x) \psi(2^{k_{0}} x) (P_{j_{2}} f)gh) \frac{1}{|x|^{3/4}} \psi(2^{k_{0}} x) \|_{L_{t,x}^{4/3}(J \times \mathbf{R}^{3})}
\endaligned
\end{equation}

\begin{equation}\label{5.40}
\aligned
+ \sum_{j > k_{0}} \sum_{j_{2} \leq j} \sum_{l \geq -k_{0}} 2^{k_{0}(1 - s)} 2^{-j(1 - s)} \| P_{j}(\chi_{l}(x) M(x) (P_{j_{2}} f)gh) \frac{1}{|x|^{3/4}} \psi(2^{k_{0}} x) \|_{L_{t,x}^{4/3}(J \times \mathbf{R}^{3})} \\
+ \sum_{j > k_{0}} \sum_{j_{2} > j} \sum_{l \geq -k_{0}} 2^{k_{0}(1 - s)} 2^{-j(1 - s)} \| P_{j}(\chi_{l}(x) M(x) (P_{j_{2}} f)gh) \frac{1}{|x|^{3/4}} \psi(2^{k_{0}} x) \|_{L_{t,x}^{4/3}(J \times \mathbf{R}^{3})}.
\endaligned
\end{equation}
Now by H{\"o}lder's inequality,

\begin{equation}\label{5.41}
\aligned
\| P_{\leq k_{0}} (M(x) \psi(2^{k_{0}} x) (P_{j_{2}} f) gh) \psi(2^{k_{0}} x) \frac{1}{|x|^{3/4}} \|_{L_{t,x}^{4/3}(J \times \mathbf{R}^{3})} \\ \lesssim 2^{-3k_{0}/4} \| P_{j_{2}} f \|_{L_{t}^{2} L_{x}^{\infty}(J \times \mathbf{R}^{3})} \| g \|_{L_{t,x}^{4}(J \times B_{1 - k_{0}})} \| h \|_{L_{t,x}^{\infty}(J \times \mathbf{R}^{3})} \\ \lesssim 2^{j_{2} (1 - s)} 2^{-k_{0}} (\sup_{j_{2}} 2^{j_{2} (s - 1)} \| P_{j_{2}} f \|_{L_{t}^{2} L_{x}^{\infty}(J \times \mathbf{R}^{3})}) \| g^{2} \|_{l_{k_{0}}^{\infty} L_{t,x}^{2}(J \times \mathbf{R}^{3})}^{1/2} \| h \|_{L_{t,x}^{\infty}(J \times \mathbf{R}^{3})},
\endaligned
\end{equation}
so then

\begin{equation}\label{5.42}
\sum_{j_{2} \leq k_{0}} (\ref{5.41}) \lesssim 2^{-k_{0} s} (\sup_{j_{2}} 2^{j_{2} (s - 1)} \mathcal S(j_{2}, f)) \| g^{2} \|_{l_{k_{0}}^{\infty} L_{t,x}^{2}(J \times \mathbf{R}^{3})}^{1/2} \| h \|_{L_{t,x}^{\infty}(J \times \mathbf{R}^{3})}.
\end{equation}
Next, by H{\"o}lder's inequality and the Sobolev embedding theorem,

\begin{equation}\label{5.43}
\aligned
\| P_{\leq k_{0}}(M(x) \psi(2^{k_{0}} x) (P_{j_{2}} f) gh) \frac{1}{|x|^{3/4}} \psi(2^{k_{0}} x) \|_{L_{t,x}^{4/3}(J \times \mathbf{R}^{3})} \\
\lesssim 2^{3k_{0}/4} \| P_{j_{2}} f \|_{L_{t,x}^{2}(J \times B_{\frac{2}{N}})} \| g \|_{L_{t,x}^{4}(J \times B_{\frac{2}{N}})} \| h \|_{L_{t,x}^{\infty}(J \times \mathbf{R}^{3})} \\
\lesssim \| P_{j_{2}} f \|_{l_{k_{0}}^{\infty} L_{t,x}^{2}(J \times \mathbf{R}^{3})} \| g^{2} \|_{l_{k_{0}}^{\infty} L_{t,x}^{2}(J \times \mathbf{R}^{3})}^{1/2} \| h \|_{L_{t,x}^{\infty}(J \times \mathbf{R}^{3})},
\endaligned
\end{equation}
so

\begin{equation}\label{5.44}
\aligned
\sum_{j_{2} > k_{0}} (\ref{5.43}) \lesssim 2^{-k_{0} s} (\sup_{j_{2} > k_{0}} 2^{j_{2} (s - 1)} \mathcal S(j_{2}, f)) \| g^{2} \|_{l_{k_{0}}^{\infty} L_{t,x}^{2}(J \times \mathbf{R}^{3})}^{1/2} \| h \|_{L_{t,x}^{\infty}(J \times \mathbf{R}^{3})}.
\endaligned
\end{equation}
Therefore, $(\ref{5.37})$ is bounded by the right hand side of $(\ref{5.34})$.\vspace{5mm}

\noindent Next, by H{\"o}lder's inequality and $(\ref{5.4})$,

\begin{equation}\label{5.45}
\aligned
\| P_{\leq k_{0}} (\chi_{l}(x) M(x) (P_{j_{2}} f) gh) \frac{1}{|x|^{3/4}} \psi(2^{k_{0}} x) \|_{L_{t,x}^{4/3}(J \times \mathbf{R}^{3})} \\
\lesssim 2^{-10 (l + k_{0})} 2^{-3k_{0}/4} \| P_{j_{2}} f \|_{L_{t}^{2} L_{x}^{\infty}(J \times \mathbf{R}^{3})} \| g \|_{L_{t,x}^{4}(J \times B_{l})} \| h \|_{L_{t,x}^{\infty}(J \times \mathbf{R}^{3})} \\
\lesssim 2^{-10(l + k_{0})} 2^{l/4} 2^{-3k_{0}/4} 2^{j_{2}(1 - s)} (2^{j_{2} (s - 1)} \mathcal S(j_{2}, f)) \| g^{2} \|_{l_{k_{0}}^{\infty} L_{t,x}^{2}(J \times \mathbf{R}^{3})}^{1/2} \| h \|_{L_{t,x}^{\infty}(J \times \mathbf{R}^{3})},
\endaligned
\end{equation}
so

\begin{equation}\label{5.46}
\sum_{l \geq -k_{0}} \sum_{j_{2} \leq k_{0}} (\ref{5.45}) \lesssim 2^{-k_{0} s} (\sup_{j_{2}} 2^{j_{2} (s - 1)} \mathcal S(j_{2}, f)) \| g^{2} \|_{l_{k_{0}}^{\infty} L_{t,x}^{2}(J \times \mathbf{R}^{3})}^{1/2} \| h \|_{L_{t,x}^{\infty}(J \times \mathbf{R}^{3})}.
\end{equation}
Also, by H{\"o}lder's inequality, the Sobolev embedding theorem, and $(\ref{5.4})$,

\begin{equation}\label{5.47}
\aligned
\| P_{\leq k_{0}} (\chi_{l}(x) M(x) (P_{j_{2}} f) gh) \frac{1}{|x|^{3/4}} \psi(2^{k_{0}} x) \|_{L_{t,x}^{4/3}(J \times \mathbf{R}^{3})} \\
\lesssim 2^{-10(l + k_{0})} 2^{3k_{0}/4} \| P_{j_{2}} f \|_{L_{t,x}^{2}(J \times B_{l})} \| g \|_{L_{t,x}^{4}(J \times B_{l})} \| h \|_{L_{t,x}^{\infty}(J \times \mathbf{R}^{3})} \\
\lesssim 2^{-10(l + k_{0})} 2^{3l/4} 2^{3k_{0}/4} \| P_{j_{2}} f \|_{l_{k_{0}}^{\infty} L_{t,x}^{2}(J \times \mathbf{R}^{3})} \| g^{2} \|_{l_{k_{0}}^{\infty} L_{t,x}^{2}(J \times \mathbf{R}^{3})}^{1/2} \| h \|_{L_{t,x}^{\infty}(J \times \mathbf{R}^{3})},
\endaligned
\end{equation}
so

\begin{equation}\label{5.48}
\aligned
\sum_{l \geq -k_{0}} \sum_{j_{2} > k_{0}} (\ref{5.47}) \lesssim 2^{-k_{0} s} (\sup_{j_{2}} 2^{j_{2} (s - 1)} \mathcal S(j_{2}, f)) \| g^{2} \|_{l_{k_{0}}^{\infty} L_{t,x}^{2}(J \times \mathbf{R}^{3})}^{1/2} \| h \|_{L_{t,x}^{\infty}(J \times \mathbf{R}^{3})}.
\endaligned
\end{equation}
Thus $(\ref{5.38})$ is bounded by the right hand side of $(\ref{5.34})$.\vspace{5mm}

\noindent Next, by H{\"o}lder's inequality and the Sobolev embedding theorem,

\begin{equation}\label{5.49}
\aligned
2^{k_{0}(1 - s)} 2^{-j(1 - s)} \| P_{j}(M(x) \psi(2^{k_{0}} x) (P_{j_{2}} f) gh) \frac{1}{|x|^{3/4}} \psi(2^{k_{0}} x) \|_{L_{t,x}^{4/3}(J \times \mathbf{R}^{3})} \\
\lesssim 2^{k_{0}(\frac{5}{8} - s)} 2^{-j(1 - s)} \| P_{j_{2}} f \|_{L_{t}^{2} L_{x}^{\infty}(J \times \mathbf{R}^{3})}^{3/4} \| P_{j_{2}} f \|_{L_{t,x}^{2}(J \times B_{1 - k_{0}})}^{1/4} \| g \|_{L_{t,x}^{4}(J \times B_{1 - k_{0}})} \| h \|_{L_{t,x}^{\infty}(J \times \mathbf{R}^{3})} \\
\lesssim 2^{k_{0}(\frac{1}{4} - s)} 2^{-j(1 - s)}  \| P_{j_{2}} f \|_{L_{t}^{2} L_{x}^{\infty}(J \times \mathbf{R}^{3})}^{3/4} \| P_{j_{2}} f \|_{l_{k_{0}}^{\infty} L_{t,x}^{2}(J \times \mathbf{R}^{3})}^{1/4} \| g^{2} \|_{l_{k_{0}}^{\infty} L_{t,x}^{2}(J \times \mathbf{R}^{3})}^{1/2} \| h \|_{L_{t,x}^{\infty}(J \times \mathbf{R}^{3})} \\
\lesssim 2^{k_{0}(\frac{1}{4} - s)} 2^{(j_{2} - j)(1 - s)} 2^{-j_{2}/4} (\sup_{j_{2}} 2^{j_{2} (s - 1)} \mathcal S(j_{2}, f)) \| g^{2} \|_{l_{k_{0}}^{\infty} L_{t,x}^{2}(J \times \mathbf{R}^{3})}^{1/2} \| h \|_{L_{t,x}^{\infty}(J \times \mathbf{R}^{3})},
\endaligned
\end{equation}
and therefore,

\begin{equation}\label{5.50}
\aligned
\sum_{j > k_{0}} \sum_{j_{2} \leq j} (\ref{5.49}) \lesssim 2^{-k_{0} s} (\sup_{j_{2}} 2^{j_{2} (s - 1)} \mathcal S(j_{2}, f)) \| g^{2} \|_{l_{k_{0}}^{\infty} L_{t,x}^{2}(J \times \mathbf{R}^{3})}^{1/2} \| h \|_{L_{t,x}^{\infty}(J \times \mathbf{R}^{3})}.
\endaligned
\end{equation}
Also by the Sobolev embedding theorem and H{\"o}lder's inequality, when $j_{2} \geq j \geq k_{0}$,

\begin{equation}\label{5.51}
\aligned
2^{k_{0}(1 - s)} 2^{-j(1 - s)} \| P_{j} (M(x) \psi(2^{k_{0}} x) (P_{j_{2}} f) gh) \frac{1}{|x|^{3/4}} \psi(2^{k_{0}} x) \|_{L_{t,x}^{4/3}(J \times \mathbf{R}^{3})} \\
\lesssim 2^{k_{0}(\frac{11}{12} - s)} 2^{-j(1 - s)} 2^{5j/6} \| P_{j_{2}} f \|_{L_{t,x}^{2}(J \times B_{1 - k_{0}})} \| g \|_{L_{t,x}^{4}(J \times B_{1 - k_{0}})} \| h \|_{L_{t,x}^{\infty}(J \times \mathbf{R}^{3})} \\
\lesssim 2^{k_{0}(\frac{1}{6} - s)} 2^{-j(1 - s)} 2^{5j/6} \| P_{j_{2}} f \|_{l_{k_{0}}^{\infty} L_{t,x}^{2}(J \times \mathbf{R}^{3})} \| g^{2} \|_{l_{k_{0}}^{\infty} L_{t,x}^{2}(J \times \mathbf{R}^{3})}^{1/2} \| h \|_{L_{t,x}^{\infty}(J \times \mathbf{R}^{3})} \\
\lesssim 2^{k_{0}(\frac{1}{6} - s)} 2^{-j(1 - s)} 2^{5j/6} 2^{-j_{2} s} (\sup_{j_{2}} 2^{-j_{2}(1 - s)} \mathcal S(j_{2}, f)) \| g^{2} \|_{l_{k_{0}}^{\infty} L_{t,x}^{2}(J \times \mathbf{R}^{3})}^{1/2} \| h \|_{L_{t,x}^{\infty}(J \times \mathbf{R}^{3})},
\endaligned
\end{equation}
so then

\begin{equation}\label{5.52}
\aligned
\sum_{j > k_{0}} \sum_{j_{2} > j} (\ref{5.51}) \lesssim 2^{-k_{0} s} (\sup_{j_{2}} 2^{-j_{2}(1 - s)} \mathcal S(j_{2}, f)) \| g^{2} \|_{l_{k_{0}}^{\infty} L_{t,x}^{2}(J \times \mathbf{R}^{3})}^{1/2} \| h \|_{L_{t,x}^{\infty}(J \times \mathbf{R}^{3})}.
\endaligned
\end{equation}
Therefore, $(\ref{5.39})$ is bounded by the right hand side of $(\ref{5.34})$.\vspace{5mm}

\noindent Finally, by H{\"o}lder's inequality and $(\ref{5.4})$,

\begin{equation}\label{5.53}
\aligned
2^{k_{0}(1 - s)} 2^{-j(1 - s)} \| P_{j}(M(x) \chi_{l}(x) (P_{j_{2}} f) gh) \frac{1}{|x|^{3/4}} \psi(2^{k_{0}} x) \|_{L_{t,x}^{4/3}(J \times \mathbf{R}^{3})} \\
\lesssim 2^{-10(l + k_{0})} 2^{k_{0}(\frac{5}{8} - s)} 2^{-j(1 - s)} \| P_{j_{2}} f \|_{L_{t}^{2} L_{x}^{\infty}(J \times \mathbf{R}^{3})}^{3/4} \| P_{j_{2}} f \|_{L_{t,x}^{2}(J \times B_{l})}^{1/4} \| g \|_{L_{t,x}^{4}(J \times B_{l})} \| h \|_{L_{t,x}^{\infty}(J \times \mathbf{R}^{3})} \\
\lesssim 2^{-10(k_{0} + l)} 2^{3l/8} 2^{k_{0}(\frac{5}{8} - s)} 2^{(j_{2} - j)(1 - s)} 2^{-j_{2}/4} (\sup_{j_{2}} 2^{-j_{2}(1 - s)} \mathcal S(j_{2}, f))   \| g^{2} \|_{l_{k_{0}}^{\infty} L_{t,x}^{2}(J \times \mathbf{R}^{3})}^{1/2} \| h \|_{L_{t,x}^{\infty}(J \times \mathbf{R}^{3})},
\endaligned
\end{equation}
and therefore,

\begin{equation}\label{5.54}
\aligned
\sum_{j > k_{0}} \sum_{j_{2} \leq j} \sum_{l \geq -k_{0}} (\ref{5.53}) \lesssim 2^{-k_{0} s} (\sup_{j_{2}} 2^{-j_{2}(1 - s)} \mathcal S(j_{2}, f))  \| g^{2} \|_{l_{k_{0}}^{\infty} L_{t,x}^{2}(J \times \mathbf{R}^{3})}^{1/2} \| h \|_{L_{t,x}^{\infty}(J \times \mathbf{R}^{3})}.
\endaligned
\end{equation}
Also by the Sobolev embedding theorem, H{\"o}lder's inequality, $j_{2} \geq j \geq k_{0}$, and $(\ref{5.4})$,

\begin{equation}\label{5.55}
\aligned
2^{k_{0}(1 - s)} 2^{-j(1 - s)} \| P_{j} (M(x) \chi_{l}(x) (P_{j_{2}} f) gh) \frac{1}{|x|^{3/4}} \psi(2^{k_{0}} x) \|_{L_{t,x}^{4/3}(J \times \mathbf{R}^{3})} \\
\lesssim 2^{-10(k_{0} + l)} 2^{k_{0}(\frac{11}{12} - s)} 2^{-j(1 - s)} 2^{5j/6} \| P_{j_{2}} f \|_{L_{t,x}^{2}(J \times B_{l})} \| g \|_{L_{t,x}^{4}(J \times B_{l})} \| h \|_{L_{t,x}^{\infty}(J \times \mathbf{R}^{3})} \\
\lesssim 2^{-10(k_{0} + l)} 2^{3l/4} 2^{k_{0}(\frac{11}{12} - s)} 2^{-j(1 - s)} 2^{5j/6} \| P_{j_{2}} f \|_{l_{k_{0}}^{\infty} L_{t,x}^{2}(J \times \mathbf{R}^{3})} \| g^{2} \|_{l_{k_{0}}^{\infty} L_{t,x}^{2}(J \times \mathbf{R}^{3})}^{1/2} \| h \|_{L_{t,x}^{\infty}(J \times \mathbf{R}^{3})} \\
\lesssim 2^{-10(k_{0} + l)} 2^{3l/4} 2^{k_{0}(\frac{11}{12} - s)} 2^{-j(\frac{1}{6} - s)} 2^{-j_{2} s} (\sup_{j_{2}} 2^{j_{2} (s - 1)} \mathcal S(j_{2}, f)) \| g^{2} \|_{l_{k_{0}}^{\infty} L_{t,x}^{2}(J \times \mathbf{R}^{3})}^{1/2} \| h \|_{L_{t,x}^{\infty}(J \times \mathbf{R}^{3})}
\endaligned
\end{equation}
so then

\begin{equation}\label{5.56}
\aligned
\sum_{j > k_{0}} \sum_{j_{2} > j} \sum_{l \geq -k_{0}} (\ref{5.55}) \lesssim 2^{-k_{0} s} (\sup_{j_{2}} 2^{j_{2} (s - 1)} \mathcal S(j_{2}, f)) \| g^{2} \|_{l_{k_{0}}^{\infty} L_{t,x}^{2}(J \times \mathbf{R}^{3})}^{1/2} \| h \|_{L_{t,x}^{\infty}(J \times \mathbf{R}^{3})}.
\endaligned
\end{equation}
Therefore, $(\ref{5.40})$ is bounded by the right hand side of $(\ref{5.34})$. This completes the proof of the theorem. $\Box$

\section{Long time Strichartz estimates}
In this section we prove proposition $\ref{p6.2}$. Recalling the proposition,

\begin{proposition}[Long time Strichartz estimates]\label{p2.1}
Suppose $N = 2^{k_{0}}$, $s = \frac{1}{2} + \epsilon$, $\mathcal J$ is an interval on which

\begin{equation}\label{2.1}
\int_{\mathcal J} \int \phi(x) (Iv(x, t))^{4} dx dt \leq C_{1} N^{2(1 - s)},
\end{equation}
and $E(Iv)(t) \leq 2 C_{0}(A, \epsilon) N^{2(1 - s)}$ for all $t \in \mathcal J$. Then for $k_{0}(C_{1}, A, \epsilon)$ sufficiently large,

\begin{equation}\label{2.2}
\sup_{j \geq k_{0} - 7} 2^{-j(1 - s)} \mathcal S(j, v) \lesssim C_{0}(A, \epsilon).
\end{equation}
\end{proposition}
%\textbf{Remark:} The implicit constant in $(\ref{2.2})$ does not depend on $C$. All implicit constants depend only on $A$ and $\epsilon > 0$.\vspace{5mm}

\noindent \emph{Proof:} To prove this we use an induction on frequency estimate combined with a bootstrap argument. Let $J \subset \mathcal J$ be an interval such that

\begin{equation}\label{2.2.1}
\| P_{> k_{0}} Iv \|_{L_{t,x}^{4}(J \times \mathbf{R}^{3})} \lesssim N^{-\epsilon/2}.
\end{equation}
By the dominated convergence theorem combined with corollary $\ref{c9.2}$, such an interval exists. We then make a bootstrap argument and show that this implies

\begin{equation}\label{2.2.2}
\| P_{> k_{0}} Iv \|_{L_{t,x}^{4}(J \times \mathbf{R}^{3})} \lesssim N^{-\epsilon},
\end{equation}
which then by the dominated convergence theorem shows that $J \subset \mathcal J$ can be extended to all of $\mathcal J$.\vspace{5mm}

%\noindent Let

%\begin{equation}\label{2.2.1}
%\mathcal S(M, J) = \sup_{M_{1} \geq M} [M_{1} (\sup_{R > 0} \| P_{M_{1}} Iv \|_{L_{t,x}^{2}(J \times \{ x : |x| \leq R \})}) + M_{1}^{1/2} \| P_{M_{1}} Iv \|_{L_{t,x}^{4}(J \times \mathbf{R}^{3})}].
%\end{equation}

\noindent By theorem $\ref{t7.6}$,

\begin{equation}\label{2.3}
\aligned
\mathcal S(j, v) \lesssim 2^{j} \| P_{j} v_{0} \|_{L^{2}(\mathbf{R}^{3})} + \| P_{j} v_{1} \|_{L_{x}^{2}(\mathbf{R}^{3})} +  \| \phi(x) (P_{> j - 7} v)(P_{\leq k_{0}} v)(P_{\leq j} v) \|_{l^{1} L_{t,x}^{2}(J \times \mathbf{R}^{3})} \\ + 2^{j} \| \phi(x) (P_{> j - 10} v)(P_{\leq k_{0}} v)(P_{> j} v) \|_{L_{t}^{2} L_{x}^{1}(J \times \mathbf{R}^{3})} + 2^{j/2} \| \phi(x)(P_{> j - 7} v)(P_{> k_{0}} v)^{2} \|_{L_{t,x}^{4/3}(J \times \mathbf{R}^{3})}  \\ + \| P_{j} (\phi(x) (P_{\leq j - 7} v)^{3}) \|_{L_{t}^{1} L_{x}^{2}(J \times \mathbf{R}^{3})}.
\endaligned
\end{equation}
Since $E(Iv)(t) \leq 2 C_{0}(A, \epsilon) N^{2(1 - s)}$ for all $t \in J$,

\begin{equation}\label{2.4}
2^{j} \| P_{j} v_{0} \|_{L^{2}(\mathbf{R}^{3})} + \| P_{j} v_{1} \|_{L_{x}^{2}(\mathbf{R}^{3})} \lesssim C_{0}(A, \epsilon) (2^{j(1 - s)} + N^{1 - s}).
\end{equation}
Next, by proposition $\ref{p5.4}$, the Sobolev embedding theorem, and the fact that $\phi(x)^{1/2} = (\frac{|x|}{\sinh |x|})$,

\begin{equation}\label{2.5}
\aligned
\| \phi(x)^{1/2} (P_{> j - 7} v)(P_{\leq k_{0}} v) \|_{L_{t,x}^{2}(J \times \mathbf{R}^{3})} \lesssim \sum_{l > j - 7} 2^{-l} \mathcal S(l, v) \| P_{\leq k_{0}} v \|_{L_{t}^{\infty} \dot{H}^{1}(J \times \mathbf{R}^{3})} \\
\lesssim |k_{0}| 2^{k_{0}(1 - s)} 2^{-js} (\sup_{l > j - 7} 2^{-l(1 - s)} \mathcal S(l, v)) \cdot C_{0}(A, \epsilon) \\ = \ln(N) N^{1 - s} 2^{-js} (\sup_{l > j - 7} 2^{-l(1 - s)} \mathcal S(l, v)) \cdot C_{0}(A, \epsilon).
\endaligned
\end{equation}
%while by the radial Sobolev embedding

%\begin{equation}\label{2.6}
%\aligned
%\sum_{\frac{1}{N} \leq 2^{j} \leq 1} \| (P_{> \frac{M}{100}} v)(P_{\leq N} v) \|_{L_{t,x}^{2}(J \times \{ x : |x| \sim 2^{j} \})} \\
%\lesssim \sum_{\frac{1}{N} \leq 2^{j} \leq 1} (\frac{M}{100} \sup_{j} 2^{-j/2} \| P_{> \frac{M}{100}} v \|_{L_{t,x}^{2}(J \times \{ x : |x| \sim 2^{j} \})}) (2^{j/2} \| P_{\leq N} v \|_{L_{t,x}^{\infty}(J \times \{ x : |x| \sim 2^{j} \})}) \\ \lesssim \ln(N) (\frac{N^{1 - s}}{M} + \frac{1}{M^{s}}) \mathcal S(J, \frac{M}{100}).
%\endaligned
%\end{equation}
%Finally, by the radial Sobolev embedding combined with the fact that $\phi(x)^{1/2} = \frac{|x|}{\sinh |x|} \lesssim e^{-c|x|}$, for some constant $c > 0$,

%\begin{equation}\label{2.7}
%\aligned
%\sum_{j \geq 0} \| \phi(x)^{1/2} (P_{> \frac{M}{100}} v)(P_{\leq N} v) \|_{L_{t,x}^{2}(J \times \{ x : |x| \sim 2^{j} \})} \\
%\lesssim \frac{1}{M}  \sum_{j \geq 0} e^{-cj} (\frac{M}{100} \sup_{j} 2^{-j/2} \| P_{> \frac{M}{100}} v \|_{L_{t,x}^{2}(J \times \{ x : |x| \sim 2^{j} \})}) (2^{j/2} \| P_{\leq N} v \|_{L_{t,x}^{\infty}(J \times \{ x : |x| \sim 2^{j} \})}) \\ \lesssim (\frac{N^{1 - s}}{M} + \frac{1}{M^{s}}) \mathcal S(J, \frac{M}{100}).
%\endaligned
%\end{equation}
%Therefore, combining $(\ref{2.5})$ - $(\ref{2.7})$,

%\begin{equation}\label{2.8}
%\| \phi(x)^{1/2} (P_{> \frac{M}{100}} v)(P_{\leq N} v) \|_{L_{t,x}^{2}(J \times \mathbf{R}^{3})} \lesssim \ln(N) (\frac{N^{1 - s}}{M} + \frac{1}{M^{s}}) \mathcal S(J, \frac{M}{100}).
%\end{equation}
\noindent Next, by Bernstein's inequality, $(\ref{2.4})$, $(\ref{2.5})$, and $\phi \in L^{\infty}$,

\begin{equation}\label{2.13}
\aligned
2^{j} \| \phi(x) (P_{> j - 7} v)(P_{\leq k_{0}} v)(P_{> j} v) \|_{L_{t}^{2} L_{x}^{1}(J \times \mathbf{R}^{3})} \\ 
\lesssim 2^{j} \| \phi(x)^{1/2} (P_{> j} v)(P_{\leq k_{0}} v) \|_{L_{t,x}^{2}(J \times \mathbf{R}^{3})} \| P_{> j} v \|_{L_{t}^{\infty} L_{x}^{2}(J \times \mathbf{R}^{3})} \\
\lesssim 2^{j(1 - s)} \ln(N) N^{1 - s} (\sup_{l > j - 7} 2^{-l(1 - s)} \mathcal S(l, v)) \cdot C_{0}(A, \epsilon) 
 \| P_{> j} v \|_{L_{t}^{\infty} L_{x}^{2}(J \times \mathbf{R}^{3})} \\ \lesssim 2^{j(1 - s)} \ln(N) N^{1 - s} (\sup_{l > j - 7} 2^{-l(1 - s)} \mathcal S(l, v)) \cdot C_{0}(A, \epsilon)^{2}  (2^{-j} N^{1 - s} +  2^{-js}).
 \endaligned
\end{equation}

\noindent Next, by the Sobolev embedding theorem, the radial Sobolev embedding theorem, the fact that $\phi(x)^{1/2} = (\frac{|x|}{\sinh |x|})$, and bounds on the energy $E(Iv(t))$,

\begin{equation}\label{2.9}
\aligned
2^{-j/2} \| P_{\leq j} v \|_{L_{t,x}^{\infty}(J \times B_{-j})} + \sum_{l > -j} 2^{l/2} \| \phi(x)^{1/2} (P_{\leq j} v) \|_{L_{t,x}^{\infty}(J \times A_{l})} \\\lesssim (|j| + 1) \| P_{\leq j} v \|_{L_{t}^{\infty} \dot{H}^{1}(J \times \mathbf{R}^{3})} \lesssim (|j| + 1) C_{0}(A, \epsilon) (N^{1 - s} + 2^{j(1 - s)}),
\endaligned
\end{equation}
so by H{\"o}lder's inequality and $(\ref{2.5})$,

%\begin{equation}\label{2.10}
%\aligned
%\sum_{\frac{1}{M} \leq 2^{j} \leq 1} 2^{j/2} \| P_{M} I(\phi(x) (P_{> \frac{M}{100}} v)(P_{\leq N} v)(P_{\leq M} v)) \|_{L_{t,x}^{2}(J \times \{ x : |x| \sim 2^{j} \})} \\ \lesssim \ln(M) \| \phi(x)^{1/2} (P_{> \frac{M}{100}} v)(P_{\leq N} v) \|_{L_{t,x}^{2}(J \times \mathbf{R}^{3})} \| |x|^{1/2} (P_{\leq M} v) \|_{L_{t,x}^{\infty}(J \times \mathbf{R}^{3})} \\
%\lesssim \ln(M)  \ln(N) (\frac{N^{1 - s}}{M} + \frac{1}{M^{s}})(N^{1 - s} + M^{1 - s}) \mathcal S(J, \frac{M}{100}).
%\endaligned
%\end{equation}
%Again since $\phi(x) \lesssim e^{-c |x|}$,

%\begin{equation}\label{2.11}
%\aligned
%\sum_{j \geq 0} 2^{j/2} \| P_{> M} I(\phi(x) (P_{> \frac{M}{100}} v)(P_{\leq N} v)(P_{\leq M} v)) \|_{L_{t,x}^{2}(J \times \{ x : |x| \sim 2^{j} \})} \\ \lesssim \| \phi(x)^{1/2} (P_{> \frac{M}{100}} v)(P_{\leq N} v) \|_{L_{t,x}^{2}(J \times \mathbf{R}^{3})} \sum_{j \geq 0} \| \phi(x)^{1/2} |x|^{1/2} (P_{\leq M} v) \|_{L_{t,x}^{\infty}(J \times \mathbf{R}^{3})} \\
%\lesssim \ln(M)  (\frac{N^{1 - s}}{M} + \frac{1}{M^{s}})(N^{1 - s} + M^{1 - s}) \mathcal S(J, \frac{M}{100}).
%\endaligned
%\end{equation}
%Therefore,

\begin{equation}\label{2.12}
\aligned
\| \phi(x) (P_{> j - 7} v)(P_{\leq k_{0}} v)(P_{\leq j} v) \|_{l^{1} L_{t, x}^{2}(J \times \mathbf{R}^{3})} \lesssim (\ref{2.9}) \cdot \| \phi(x)^{1/2} (P_{> j} v)(P_{\leq k_{0}} v) \|_{l^{\infty} L_{t,x}^{2}(J \times \mathbf{R}^{3})} \\ \lesssim (\sup_{l > j - 7} 2^{-l(1 - s)} \mathcal S(l, v)) \cdot C_{0}(A, \epsilon)^{2} 2^{j(1 - s)} \ln(N) (|j| + 1) N^{1 - s} (2^{-j} N^{1 - s} +  2^{-js}).
\endaligned
\end{equation}
Next, by the bootstrap assumption $\| P_{> k_{0}} v \|_{L_{t,x}^{4}(J \times \mathbf{R}^{3})} \lesssim N^{-\epsilon/2}$,

\begin{equation}\label{2.14}
\aligned
2^{j/2} \| \phi(x)(P_{> j - 7} v)(P_{> k_{0}} v)^{2} \|_{L_{t,x}^{4/3}(J \times \mathbf{R}^{3})} \\ \lesssim 2^{j/2} \| P_{> j - 7} v \|_{L_{t,x}^{4}(J \times \mathbf{R}^{3})} \| P_{> k_{0}} v \|_{L_{t,x}^{4}(J \times \mathbf{R}^{3})}^{2} \lesssim 2^{j(1 - s)} N^{-\epsilon} (\sup_{l > j - 7} 2^{-l(1 - s)} \mathcal S(l, v)).
\endaligned
\end{equation}
Finally, decompose

\begin{equation}\label{2.15}
\aligned
\phi(x) = ((P_{> j - 7} \phi(x)^{1/2})(P_{> j - 7} \phi(x)^{1/2})) + ((P_{\leq j - 7} \phi(x)^{1/2})(P_{\leq j - 7} \phi(x)^{1/2})) \\
+ 2 ((P_{> j - 7} \phi(x)^{1/2})\phi(x)^{1/2}) - 2 ((P_{> j - 7} \phi(x)^{1/2})(P_{> j - 7} \phi(x)^{1/2})).
\endaligned
\end{equation}
Because $(\frac{|x|}{\sinh |x|})$ and all its derivatives are smooth and rapidly decreasing, for any $x$, $k$, $j > 0$,

\begin{equation}\label{2.16}
|P_{> j - 7} \phi(x)^{1/2}| \lesssim_{k} 2^{-j} (1 + |x|)^{-k},
\end{equation}
so by the radial Sobolev embedding theorem and a crude summation of local bounds (in particular corollary $\ref{c9.2}$),

\begin{equation}\label{2.17}
\aligned
\| (P_{> j - 7} \phi(x)^{1/2})^{2} (P_{\leq j - 7} v)^{3} \|_{L_{t}^{1} L_{x}^{2}(J \times \mathbf{R}^{3})} \\ \lesssim 2^{-2j} \| P_{\leq j - 7} v \|_{L_{t}^{2} L_{x}^{\infty}(J \times \mathbf{R}^{3})}^{2} \| |x|^{1/2} P_{\leq j - 7} v \|_{L_{t,x}^{\infty}(J \times \mathbf{R}^{3})} \\ \lesssim
2^{-2j} (1 + \frac{2^{j(1 - s)}}{N^{1 - s}})^{3}\| Iv \|_{L_{t}^{2} L_{x}^{\infty}(J \times \mathbf{R}^{3})}^{2} \| Iv \|_{L_{t}^{\infty} \dot{H}^{1}(J \times \mathbf{R}^{3})} \\ \lesssim 2^{-2j} C_{1}^{2} C_{0}(A, \epsilon)^{3} N^{5(1 - s)} (1 + \frac{2^{j(1 - s)}}{N^{1 - s}})^{3} .
\endaligned
\end{equation}
Next, using $(\ref{2.1})$, corollary $\ref{c9.2}$, and splitting $P_{\leq j - 7} v = Iv - P_{> j - 7} Iv + P_{\leq j - 7} (1 - I) v$,

\begin{equation}\label{2.18}
\aligned
\| (P_{> j - 7} \phi(x)^{1/2}) \phi(x)^{1/2} (Iv)^{2}(P_{\leq j - 7} v) \|_{L_{t}^{1} L_{x}^{2}(J \times \mathbf{R}^{3})} \\ \lesssim \| \phi(x)^{1/2} (Iv)^{2} \|_{L_{t,x}^{2}(J \times \mathbf{R}^{3})} \| P_{\leq j - 7} v \|_{L_{t}^{2} L_{x}^{\infty}(J \times \mathbf{R}^{3})} \| P_{> j - 7} \phi(x)^{1/2} \|_{L_{x}^{\infty}(\mathbf{R}^{3})} \\ \lesssim 2^{-j} C_{1}^{2} C_{0}(A, \epsilon) N^{3(1 - s)} (1 + \frac{2^{j(1 - s)}}{N^{1 - s}}),
\endaligned
\end{equation}
and by Bernstein's inequality,

\begin{equation}\label{2.19}
\aligned
\| (P_{> j - 7} \phi(x)^{1/2}) \phi(x)^{1/2} (P_{> j - 7} Iv)^{2}(P_{\leq j - 7} v) \|_{L_{t}^{1} L_{x}^{2}(J \times \mathbf{R}^{3})} \\
\lesssim 2^{-j} \| Iv \|_{L_{t}^{2} L_{x}^{\infty}(J \times \mathbf{R}^{3})} \| P_{\leq j - 7} v \|_{L_{t}^{2} L_{x}^{\infty}(J \times \mathbf{R}^{3})} \| P_{> j - 7} Iv \|_{L_{t}^{\infty} L_{x}^{2}(J \times \mathbf{R}^{3})} \\ \lesssim 2^{-2j} C_{1}^{2} C_{0}(A, \epsilon)^{3} N^{5(1 - s)} (1 + \frac{2^{j(1 - s)}}{N^{1 - s}}).
\endaligned
\end{equation}
Also, by definition of the $I$ - operator and the bootstrap assumption,

\begin{equation}\label{2.20}
\aligned
\| \phi(x)^{1/2} (P_{> j - 7} \phi(x)^{1/2})(P_{\leq j - 7} (1 - I)v)^{3} \|_{L_{t}^{1} L_{x}^{2}(J \times \mathbf{R}^{3})} \\ \lesssim  \| P_{> j - 7} \phi(x)^{1/2} \|_{L_{x}^{\infty}(\mathbf{R}^{3})} \| P_{\leq j - 7} (1 - I)v \|_{L_{t}^{2} L_{x}^{\infty}(J \times \mathbf{R}^{3})} \| P_{> k_{0}} v \|_{L_{t,x}^{4}(J \times \mathbf{R}^{3})}^{2} \\ \lesssim 2^{-j} C_{1}(A, \epsilon) C_{0}(A, \epsilon) (1 + \frac{2^{j(1 - s)}}{N^{1 - s}}) N^{2(1 - s)} N^{-\epsilon}.
\endaligned
\end{equation}
Since

\begin{equation}\label{2.21}
P_{j} ((P_{\leq j - 7} v)^{3} (P_{\leq j - 7} \phi(x)^{1/2})^{2} = 0,
\end{equation}
then by $(\ref{2.3})$, $(\ref{2.4})$, $(\ref{2.13})$, $(\ref{2.12})$, $(\ref{2.14})$, $(\ref{2.17})$, $(\ref{2.18})$,

\begin{equation}\label{2.22}
\aligned
2^{-j(1 - s)} \mathcal S(j, v) \lesssim C_{0}(A, \epsilon) N^{1 - s} + N^{-\epsilon} (\sup_{l > j - 7} 2^{-l(1 - s)} \mathcal S(j, v)) \\ + \ln(N) (|j| + 1) C_{0}(A, \epsilon)^{2} (2^{-j} N^{2(1 - s)} + 2^{-js} N^{1 - s}) (\sup_{l > j - 7} 2^{-l(1 - s)} \mathcal S(l, v)) \\ 
+ 2^{-2j} 2^{-j(1 - s)} C_{1}^{2} C_{0}(A, \epsilon)^{3} N^{5(1 - s)} (1 + \frac{2^{j(1 - s)}}{N^{1 - s}})^{3} \\ + 2^{-j} 2^{-j(1 - s)} C_{1}^{2} C_{0}(A, \epsilon) N^{3(1 - s)} (1 + \frac{2^{j(1 - s)}}{N^{1 - s}}).
\endaligned
\end{equation}

%\noindent Now then, for any $M$,

%\begin{equation}\label{2.23}
%M \sum_{M_{1} \geq M} (\sup_{R > 0} R^{-1/2} \| P_{M_{1}} Iv \|_{L_{t,x}^{2}(J \times \{ |x| \leq R \})}) + M^{1/2} \sum_{M_{1} \geq M} \| P_{M_{1}} Iv \|_{L_{t,x}^{4}(J \times \mathbf{R}^{3})}
%\end{equation}

%\begin{equation}\label{2.24}
%\lesssim \sup_{M_{1} \geq M} M_{1} (\sup_{R > 0} R^{-1/2} \| P_{M_{1}} Iv \|_{L_{t,x}^{2}(J \times \{ |x| \leq R \})}) + \sup_{M_{1} \geq M} M_{1}^{1/2} \| P_{M_{1}} Iv \|_{L_{t,x}^{4}(J \times \mathbf{R}^{3})},
%\end{equation}
%and therefore,

%\begin{equation}\label{2.25}
%\aligned
%M (\sup_{R > 0} R^{-1/2} \| P_{> M} Iv \|_{L_{t,x}^{2}(J \times \{ |x| \leq R \})}) + M^{1/2} \| P_{> M} Iv \|_{L_{t,x}^{4}(J \times \mathbf{R}^{3})} \\ \lesssim \ln(N) \ln(M) \sup(\frac{N^{2(1 - s)}}{M}, \frac{N^{1 - s}}{M^{s}}) (\frac{M}{100} \sup_{j} 2^{-j/2} \| P_{> \frac{M}{100}} v \|_{L_{t,x}^{2}(J \times \{ x : |x| \sim 2^{j} \})}) \\ + \| P_{> N} v \|_{L_{t,x}^{4}(J \times \mathbf{R}^{3})}^{2} (\frac{M^{1/2}}{10} \| P_{> \frac{M}{100}} Iv \|_{L_{t,x}^{4}(J \times \mathbf{R}^{3})}) + \frac{N^{1 - s}}{M^{s}} \| P_{> N} v \|_{L_{t,x}^{4}(J \times \mathbf{R}^{3})}^{2} \\ + \frac{1}{M^{2}} N^{5(1 - s)} + \frac{1}{M} N^{3(1 - s)} + N^{1 - s}.
%\endaligned
%\end{equation}
\noindent Choosing $N(C_{1}, A, \epsilon)$ sufficiently large, then for $j$ such that $2^{j} \geq N^{1 - \epsilon}$,

\begin{equation}\label{2.26}
\ln(N) (|j| + 1) C_{0}(A, \epsilon)^{2} (2^{-j} N^{2(1 - s)} + 2^{-js} N^{1 - s}) + N^{-\epsilon} \leq N^{-\epsilon/10},
\end{equation}
and

\begin{equation}\label{2.27}
2^{-2j} C_{1}^{2} C_{0}(A, \epsilon)^{3} N^{5(1 - s)} (1 + \frac{2^{j(1 - s)}}{N^{1 - s}})^{3} + 2^{-j} C_{1}^{2} C_{0}(A, \epsilon) N^{3(1 - s)} (1 + \frac{2^{j(1 - s)}}{N^{1 - s}}) << 1,
\end{equation}
so

\begin{equation}\label{2.28}
\aligned
(\sup_{l \geq j} 2^{-l(1 - s)} \mathcal S(l, v)) \lesssim C_{0}(A, \epsilon) (1 + 2^{-j(1 - s)} N^{1 - s}) + N^{-\epsilon/10} (\sup_{l \geq j - 7} 2^{-l(1 - s)} \mathcal S(l, v)).
\endaligned
\end{equation}
Then by induction on $j$, starting with $j \geq k_{0}(1 - \epsilon)$, and thus $2^{j} \geq N^{1 - \epsilon}$, and corollary $\ref{c9.2}$, there exists some $c > 0$ such that

\begin{equation}\label{2.29}
(\sup_{l \geq k_{0} - 7} 2^{-l(1 - s)} \mathcal S(l, v)) \lesssim C_{0}(A, \epsilon) + N^{-c \epsilon \ln(N)} C_{0}(A, \epsilon) C_{1} N^{2(1 - s)},
\end{equation}
which proves the theorem for $J \subset \mathcal J$. But then this implies $(\ref{2.2.2})$ holds for $J \subset \mathcal J$, and thus by local well - posedness (lemma $\ref{l9.0}$) there exists a larger open interval $J \subset J_{1} \subset \mathcal J$ for which $(\ref{2.2.1})$ holds on then closure of $J_{1}$. By the usual bootstrap arguments, this proves the theorem for all of $\mathcal J$. $\Box$

\begin{corollary}\label{c2.2}
For $k_{0}(C_{1}, A, \epsilon)$ sufficiently large, $N = 2^{k_{0}}$,

%\begin{equation}\label{2.29}
%\sup_{R > 0} R^{-1} \int_{\mathcal J} \int_{|x| \leq R} |\nabla I P_{\geq \frac{N}{100}} v(x,t)|^{2} dx dt + \sup_{R > 0} R^{-1} \int_{\mathcal J} \int_{|x| \leq R} |P_{\geq \frac{N}{100}} Iv_{t}(x,t)|^{2} dx dt \lesssim N^{2(1 - s)}.
%\end{equation}
%Also,

\begin{equation}\label{2.29.1}
\| I P_{\geq k_{0} - 7} v \|_{L_{t}^{2} L_{x}^{\infty}(\mathcal J \times \mathbf{R}^{3})} + \| \nabla I P_{\geq k_{0} - 7} v \|_{l^{\infty} L_{t,x}^{2}(J \times \mathbf{R}^{3})} + \| I P_{\geq k_{0} - 7} v_{t} \|_{l^{\infty} L_{t,x}^{2}(J \times \mathbf{R}^{3})} \lesssim C_{0}(A, \epsilon) N^{1 - s}.
\end{equation}

%\begin{equation}\label{2.29.2}
%\| I P_{\geq \frac{N}{100}} \nabla Iv \|_{L_{t}^{q} L_{x}^{p}(\mathcal J \times \mathbf{R}^{3})} + \| I P_{\geq \frac{N}{100}} Iv_{t} \|_{L_{t}^{q} L_{x}^{p}(\mathcal J \times \mathbf{R}^{3})} \lesssim N^{2(1 - s)},
%\end{equation}

%\begin{equation}\label{2.29.2.1}
%\| P_{> \frac{N}{100}} v \|_{L_{t}^{p} L_{x}^{q}(\mathcal J \times \mathbf{R}^{3})} \lesssim 1,
%\end{equation}
%where

\end{corollary}
\emph{Proof:} Let $(p, q)$ be an $\dot{H}^{\frac{1}{2} + \epsilon}$ - admissible pair

\begin{equation}\label{2.29.3}
\frac{1}{p} = \frac{1}{4} + \frac{\epsilon}{2}, \hspace{5mm} \frac{1}{q} = \frac{1}{4} - \frac{\epsilon}{2}. 
\end{equation}
By interpolation and proposition $\ref{p2.1}$,

\begin{equation}\label{2.30}
\sup_{j \geq k_{0} - 14} \| P_{j} v \|_{L_{t}^{p} L_{x}^{q}(\mathcal J \times \mathbf{R}^{3})} \lesssim C_{0}(A, \epsilon).
\end{equation}
%Therefore,

%\begin{equation}\label{2.31}
%\sup_{M \geq N} \| P_{M} v \|_{L_{t}^{p} L_{x}^{q}(\mathcal J \times \mathbf{R}^{3})} \lesssim 1.
%\end{equation}
\noindent Then let

\begin{equation}\label{2.32}
\frac{1}{\tilde{p}} = \frac{1}{p} + \frac{1}{2}, \hspace{5mm} \frac{1}{\tilde{q}} = \frac{1}{q} + \frac{1}{2}.
\end{equation}
Combining proposition $\ref{p2.1}$ with $(\ref{2.30})$,

\begin{equation}\label{2.33}
\aligned
\| (P_{> k_{0} - 14} v)^{3} \|_{L_{t}^{\tilde{p}} L_{x}^{\tilde{q}}(\mathcal J \times \mathbf{R}^{3})}  &\lesssim \sum_{k_{0} - 14 \leq j_{1} \leq j_{2} \leq j_{3}} \| P_{j_{1}} v \|_{L_{t}^{p} L_{x}^{q}(\mathcal J \times \mathbf{R}^{3})} \| P_{j_{2}} v \|_{L_{t,x}^{4}(\mathcal J \times \mathbf{R}^{3})} \| P_{j_{3}} v \|_{L_{t,x}^{4}(\mathcal J \times \mathbf{R}^{3})} \\ &\lesssim \sum_{k_{0} - 14 \leq j_{1} \leq j_{2} \leq j_{3}} 2^{-\epsilon j_{2}} 2^{-\epsilon j_{3}} \lesssim N^{-2 \epsilon}.
\endaligned
\end{equation}
Meanwhile, by $(\ref{2.12})$ and the conclusion of proposition $\ref{p2.1}$,

\begin{equation}\label{2.34}
\aligned
\sum_{j > k_{0} - 7} \| \phi(x) (P_{> j - 7} v)(P_{\leq k_{0}} v)^{2} \|_{l_{j}^{1} L_{t,x}^{2}(\mathcal J \times \mathbf{R}^{3})} \lesssim C_{0}(A, \epsilon)^{3} \sum_{j > k_{0} - 7} 2^{-\frac{\epsilon}{2} j} N^{1 - s} \lesssim C_{0}(A, \epsilon) N^{1 - s}.
\endaligned
\end{equation}
Finally by $(\ref{2.17})$ - $(\ref{2.20})$,

\begin{equation}\label{2.35}
\aligned
\| I P_{> k_{0} - 7} (\phi(x) (P_{\leq k_{0} - 14} v)^{3}) \|_{L_{t}^{1} L_{x}^{2}(\mathcal J \times \mathbf{R}^{3})} \\ \lesssim \sum_{j > k_{0} - 7} \| \tilde{P}_{j} (\phi(x) (P_{\leq k_{0} - 14} v)^{3}) \|_{L_{t}^{1} L_{x}^{2}(\mathcal J \times \mathbf{R}^{3})} \lesssim C_{0}(A, \epsilon) N^{1 - s}.
\endaligned
\end{equation}
Therefore, since 

\begin{equation}
P_{> k_{0} - 7} I \sim N^{1 - s} |\nabla|^{s - 1},
\end{equation}
by theorem $\ref{t7.7}$, $(\ref{2.33})$ - $(\ref{2.35})$,

\begin{equation}\label{2.36}
\| \nabla I P_{> k_{0} - 7} v \|_{l^{\infty} L_{t,x}^{2}(\mathcal J \times \mathbf{R}^{3})} + \| I P_{> k_{0} - 7} v_{t} \|_{l^{\infty} L_{t,x}^{2}(\mathcal J \times \mathbf{R}^{3})} + \| I P_{> k_{0} - 7} Iv \|_{L_{t}^{2} L_{x}^{\infty}(\mathcal J \times \mathbf{R}^{3})} \lesssim C_{0}(A, \epsilon) N^{1 - s}.
\end{equation}
This proves the theorem. $\Box$

\section{Change of energy}
Next recall proposition $\ref{p6.3}$.
\begin{theorem}\label{t3.1}
Suppose $\mathcal J$ is an interval on which

\begin{equation}\label{3.0}
\int_{\mathcal J} \int (\frac{\cosh |x|}{\sinh |x|}) \phi(x) (Iv(x,t))^{4} dx dt \leq C_{1} N^{1 - 2 \epsilon} = C_{1} N^{2(1 - s)},
\end{equation}
and $\sup_{t \in \mathcal J} E(Iv(t)) \leq 2 C_{0}(A, \epsilon) N^{2(1 - s)}$. Then for $k_{0}(C_{1}, A, \epsilon)$ sufficiently large,

\begin{equation}\label{3.1}
\sup_{t \in \mathcal J} E(Iv(t)) \leq \frac{3}{2} C_{0}(A, \epsilon) N^{2(1 - s)}.
\end{equation}
\end{theorem}
\emph{Proof:} By proposition $\ref{p2.1}$ and corollary $\ref{c2.2}$, we can choose $k_{0}(C_{1}, A, \epsilon)$ sufficiently large so that

\begin{equation}\label{3.2}
\sup_{j \geq k_{0} - 7} 2^{-j(1 - s)} \mathcal S(j, v) \lesssim C_{0}(A, \epsilon),
\end{equation}
and

\begin{equation}\label{3.3}
\aligned
\| \nabla I P_{\geq k_{0} - 7} v \|_{l^{\infty} L_{t,x}^{2}(\mathcal J \times \mathbf{R}^{3})} + \| I P_{\geq k_{0} - 7} v_{t} \|_{l^{\infty} L_{t,x}^{2}(\mathcal J \times \mathbf{R}^{3})} \\ + \| I P_{\geq k_{0} - 7} v \|_{L_{t}^{2} L_{x}^{\infty}(\mathcal J \times \mathbf{R}^{3})} \lesssim C_{0}(A, \epsilon) N^{1 - s}.
\endaligned
\end{equation}
The change of the modified energy

\begin{equation}\label{3.4}
E(Iv(t)) = \frac{1}{2} \int |\nabla Iv(x,t)|^{2} dx + \frac{1}{2} \int (Iv_{t}(x,t))^{2} dx + \frac{1}{4} \int \phi(x) (Iv(x,t))^{4} dx,
\end{equation}
is given by

\begin{equation}\label{3.5}
\frac{d}{dt} E(Iv(t)) = \int (Iv_{t}(x,t)) \cdot [\phi(x) (Iv(x,t))^{3} - I(\phi(x) v(x,t)^{3})] dx.
\end{equation}

\noindent \textbf{Remark:} When $I = 1$ then the energy is clearly conserved.\vspace{5mm}

\noindent First observe that $(\ref{3.0})$, $(\ref{3.2})$, and the fact that $\frac{1}{|x|} \lesssim \frac{\cosh |x|}{\sinh |x|}$ imply

\begin{equation}\label{3.16}
\aligned
\| \phi(x)^{1/2} (P_{\leq k_{0}} v)^{2} \|_{l_{k_{0}}^{\infty} L_{t,x}^{2}(J \times \mathbf{R}^{3})} \lesssim \| \frac{\phi(x)^{1/4}}{|x|^{1/4}} (Iv) \|_{L_{t,x}^{4}(J \times \mathbf{R}^{3})}^{2} + N^{1/2} \| P_{> k_{0}} Iv \|_{L_{t,x}^{4}(J \times \mathbf{R}^{3})}^{2} \\ \lesssim C_{1}^{1/2} N^{1 - s} + C_{0}(A, \epsilon) N^{1/2 - 2 \epsilon},
\endaligned
\end{equation}
and

\begin{equation}\label{3.16.1}
\| (P_{> k_{0}} v)^{2} \|_{l_{k_{0}}^{\infty} L_{t,x}^{2}(J \times \mathbf{R}^{3})} + \| (P_{> k_{0}} Iv)^{2} \|_{l_{k_{0}}^{\infty} L_{t,x}^{2}(J \times \mathbf{R}^{3})} \lesssim C_{0}(A, \epsilon) N^{1/2 - 2 \epsilon}.
\end{equation}
Therefore, by theorem $\ref{t5.1}$, proposition $\ref{p2.1}$, the bootstrap assumption $(\ref{3.0})$, and the decay of $\phi(x)^{1/2}$,

\begin{equation}\label{3.23}
\aligned
\int_{J} \int \phi(x) (P_{> k_{0} - 7} Iv)(Iv)^{2} (I P_{\leq k_{0}} v_{t}) dx dt \\
\lesssim \ln(N) N^{1/2} \| \phi(x)^{1/2} (Iv)^{2} \|_{l_{k_{0}}^{\infty} L_{t,x}^{2}(J \times \mathbf{R}^{3})} \| P_{> k_{0} - 7} Iv \|_{l_{k_{0}}^{\infty} L_{t,x}^{2}(J \times \mathbf{R}^{3})} \| P_{\leq k_{0} - 7} Iv_{t} \|_{L_{t}^{\infty} L_{x}^{2}(J \times \mathbf{R}^{3})} \\
\lesssim \ln(N) N^{1/2 - s} C_{0}(A, \epsilon)^{2} C_{1}^{1/2} N^{2(1 - s)}.
\endaligned
\end{equation}

\noindent Also, by theorem $\ref{t5.1}$, proposition $\ref{p2.1}$, and $(\ref{3.0})$,

\begin{equation}\label{3.14}
\aligned
\int_{J} \int I(\phi(x) (P_{> k_{0} - 7} v)v^{2}) (I P_{\leq k_{0}} v_{t}) dx dt \\ 
\lesssim \ln(N) N^{1/2} \| P_{\leq k_{0}} Iv_{t} \|_{L_{t}^{\infty} L_{x}^{2}(J \times \mathbf{R}^{3})} \| \phi(x)^{1/2} v^{2} \|_{l_{k_{0}}^{\infty} L_{t,x}^{2}(J \times \mathbf{R}^{3})} \| P_{\geq k_{0} - 7} v \|_{l_{k_{0}}^{\infty} L_{t,x}^{2}(J \times \mathbf{R}^{3})} \\
\lesssim \ln(N) N^{1/2} \| I v_{t} \|_{L_{t}^{\infty} L_{x}^{2}(J \times \mathbf{R}^{3})} (C_{1}^{1/2} N^{1 - s} + C_{0}(A, \epsilon) N^{1/2 - 2 \epsilon})  \sum_{j > k_{0} - 7} 2^{-js} \mathcal S(j, v) \\
\lesssim \ln(N) N^{1/2 - s} C_{0}(A, \epsilon)^{2} N^{1 - s} (C_{1}^{1/2} N^{1 - s} + C_{0}(A, \epsilon) N^{1/2 - 2 \epsilon}).
\endaligned
\end{equation}

\noindent Also, by $(\ref{5.2.1})$, proposition $\ref{p2.1}$, corollary $\ref{c2.2}$, $(\ref{3.0})$, the radial Sobolev embedding theorem, and Bernstein's inequality

\begin{equation}\label{3.22}
\aligned
\int_{J} \int_{|x| \geq \frac{2}{N}} \phi(x) (P_{> k_{0} - 7} Iv)(Iv)^{2} (I P_{> k_{0}} v_{t}) dx dt \\ \lesssim \ln(N) \| \phi(x)^{1/2} (Iv)^{2} \|_{l_{k_{0}}^{\infty} L_{t,x}^{2}(J \times \mathbf{R}^{3})} \| I P_{> k_{0}} Iv_{t} \|_{l_{k_{0}}^{\infty} L_{t,x}^{2}(J \times \mathbf{R}^{3})} \| |x| (I P_{> k_{0}} v) \|_{L_{t,x}^{\infty}(J \times \mathbf{R}^{3})} \\
\lesssim \ln(N) C_{0}(A, \epsilon) C_{1}^{1/2} N^{2(1 - s)}  \sum_{j > k_{0} - 7} 2^{-j/2} \| P_{j} Iv \|_{L_{t}^{\infty} \dot{H}^{1}(J \times \mathbf{R}^{3})} \\
\lesssim \ln(N) C_{0}(A, \epsilon)^{2} C_{1}^{1/2} N^{2(1 - s)}  N^{1/2 - s},
\endaligned
\end{equation}
while

\begin{equation}\label{3.22.1}
\aligned
\int_{J} \int_{|x| \leq \frac{2}{N}} \phi(x) (P_{> k_{0} - 7} Iv) (Iv)^{2} (I P_{> k_{0}} v_{t}) dx dt \\
\lesssim \| Iv \|_{L_{t,x}^{4}(J \times B_{\frac{2}{N}})}^{2} \| P_{> k_{0} - 7} Iv \|_{L_{t}^{2} L_{x}^{\infty}(J \times \mathbf{R}^{3})} \| IP_{> k_{0}} v_{t} \|_{L_{t}^{\infty} L_{x}^{2}(J \times \mathbf{R}^{3})} \\ \lesssim C_{0}(A, \epsilon)^{2} N^{2(1 - s)} N^{-1/2} \| (Iv)^{2} \|_{l_{k_{0}}^{\infty} L_{t,x}^{2}(J \times \mathbf{R}^{3})} \\ \lesssim C_{0}(A, \epsilon)^{2} C_{1}^{1/2} N^{\frac{3}{2} - 2s} N^{1 - s}.
\endaligned
\end{equation}
Meanwhile, by the Sobolev embedding theorem and the bounds on the energy,

\begin{equation}
\| P_{\leq k_{0} - 7} v \|_{L_{t}^{\infty} \dot{H}^{1}(J \times \mathbf{R}^{3})} + N^{-1/2} \| P_{\leq k_{0} - 7} v \|_{L_{t,x}^{\infty}(J \times \mathbf{R}^{3})} \lesssim \sup_{t \in \mathcal J} E(Iv(t))^{1/2} \lesssim C_{0}(A, \epsilon) N^{1 - s},
\end{equation}
so by $(\ref{5.16})$, proposition $\ref{p2.1}$, corollary $\ref{c2.2}$, and proposition $\ref{p5.4}$,
%Next, by $(\ref{2.7})$ and corollary $\ref{c2.2}$,

%\begin{equation}\label{3.16}
%\aligned
%\int_{\mathcal J} \int (Iv_{t}(x,t)) \phi(x) (I P_{\geq \frac{N}{100}} v)^{2} (P_{\leq \frac{N}{100}} v) dx dt \\ \lesssim \| \phi(x)^{1/2} (IP_{\geq \frac{N}{100}} v)(P_{\leq \frac{N}{100}} v) \|_{L_{t,x}^{2}(\mathcal J \times \mathbf{R}^{3})} \| I P_{\geq \frac{N}{100}} v \|_{L_{t}^{2} L_{x}^{\infty}(\mathcal J \times \mathbf{R}^{3})} \| Iv_{t} \|_{L_{t}^{\infty} L_{x}^{2}(\mathcal J \times \mathbf{R}^{3})} \\
%\lesssim \ln(N) N^{-2 \epsilon} N^{2(1 - s)}.
%\endaligned
%\end{equation}
%Also by corollary $\ref{c2.2}$ and $(\ref{2.7})$,

\begin{equation}\label{3.17}
\aligned
\int_{J} \int I(\phi(x) v(x,t) (P_{> k_{0} - 7} v)(P_{\leq k_{0}} v))(x,t)  (I P_{> k_{0}} Iv_{t})(x,t)  dx dt \\
\lesssim \ln(N) C_{0}(A, \epsilon) N^{1 - s} \| I P_{\geq k_{0}} v_{t} \|_{l_{k_{0}}^{\infty} L_{t,x}^{2}(J \times \mathbf{R}^{3})} \| \phi(x)^{1/2} (P_{> k_{0} - 7} v) v \|_{L_{t,x}^{2}(J \times \mathbf{R}^{3})} \\
\lesssim \ln(N) N^{2(1 - s)} C_{0}(A, \epsilon)^{2} (\sum_{j > k_{0} - 7} 2^{-j} \mathcal S(j, v) (\| P_{\leq k_{0}} v \|_{L_{t}^{\infty} \dot{H}^{1}} \\ + N^{-1/2} \| P_{\leq k_{0}} v \|_{L_{t,x}^{\infty}}) + \| P_{> k_{0} - 7} v \|_{L_{t,x}^{4}(J \times \mathbf{R}^{3})}^{2}) \\
\lesssim \ln(N) C_{0}(A, \epsilon)^{4} N^{2(1 - s)} (\ln(N) N^{1 - 2s} + N^{-2 \epsilon}).
\endaligned
\end{equation}

\noindent Also, by proposition $\ref{p2.1}$ and theorem $\ref{t5.3}$,

\begin{equation}\label{3.11}
\aligned
\| I(\phi(x) (P_{> k_{0} - 7} v)^{3}) (I P_{\geq k_{0}} v_{t}) \|_{L_{t, x}^{1}(J \times \mathbf{R}^{3})} \\ \lesssim \| I(\phi(x) (P_{> k_{0} - 7} v)^{3}) \|_{L_{t}^{1} L_{x}^{2}(J \times \mathbf{R}^{3})} \| I v_{t} \|_{L_{t}^{\infty} L_{x}^{2}(J \times \mathbf{R}^{3})} \\ \lesssim C_{0}(A, \epsilon)^{3} N^{1 - s - \epsilon} (\sup_{t \in J} E(Iv(t)))^{1/2} \lesssim C_{0}(A, \epsilon)^{4} N^{2(1 - s)} N^{-\epsilon}.
\endaligned
\end{equation}

%\noindent Also, for any $j > 0$, since $\phi(x) \lesssim e^{-2c |x|}$ for some $c > 0$, by the radial Sobolev embedding, 

%\begin{equation}\label{3.23.1}
%\aligned
%\int_{J} \int \phi(x) \chi_{j}(x) (Iv)^{2} (IP_{\geq \frac{N}{100}} v) Iv_{t}(x,t) dx dt \\
%\lesssim e^{-c \frac{2^{j}}{N}} \| \frac{\phi(x)^{1/4}}{|x|^{1/4}} Iv \|_{L_{t,x}^{4}(J \times \mathbf{R}^{3})}^{2} \| P_{\geq \frac{N}{100}} Iv \|_{L_{t,x}^{2}(J \times \{ x : |x| \leq \frac{2^{j + 1}}{N} \})} \| \chi_{j}(x) |x|^{1/2} P_{\leq \frac{N}{100}} Iv_{t} \|_{L_{t,x}^{\infty}(J \times \mathbf{R}^{3})} \\
%\lesssim e^{-c \frac{2^{j}}{N}} N^{1 - s} N^{1/2 - s} \| \frac{\phi(x)^{1/4}}{|x|^{1/4}} Iv \|_{L_{t,x}^{4}(J \times \mathbf{R}^{3})}^{2}.
%\endaligned
%\end{equation}
%Summing in $j$,

%\begin{equation}\label{3.23.2}
%\sum_{j > 0} (\ref{3.23.1}) \lesssim \ln(N) N^{1 - s} N^{-\epsilon} \| \frac{\phi(x)^{1/4}}{|x|^{1/4}} Iv \|_{L_{t,x}^{4}(J \times \mathbf{R}^{3})}^{2} << N^{2(1 - s)}.
%\end{equation}

\noindent Now then,

\begin{equation}\label{3.6}
\phi(x) (IP_{\leq k_{0} - 7} v(x,t))^{3} - I(\phi(x) P_{\leq k_{0} - 7} v(x,t)^{3}) = (1 - I) (\phi(x) (P_{\leq k_{0} - 7} v(x,t))^{3}).
\end{equation}
Recall from $(\ref{2.17})$ - $(\ref{2.20})$ that

\begin{equation}\label{3.7}
\aligned
\| (1 - I) (\phi(x) (P_{\leq k_{0} - 7} v)^{3}) \|_{L_{t}^{1} L_{x}^{2}(J \times \mathbf{R}^{3})} \\ \lesssim C_{1}^{2} C_{0}(A, \epsilon)^{3} \sum_{j > k_{0}} 2^{-2j} N^{5(1 - s)} + 2^{-j} N^{3(1 - s)} + 2^{-2j} 2^{3j(1 - s)} N^{2(1 - s)} + 2^{-js} N^{2(1 - s)} \\
\lesssim C_{1}^{2} C_{0}(A, \epsilon)^{3} (N^{3 - 5s} + N^{2 - 3s}).
\endaligned
\end{equation}
Since

\begin{equation}
\| Iv_{t} \|_{L_{t}^{\infty} L_{x}^{2}(J \times \mathbf{R}^{3})} \lesssim C_{0}(A, \epsilon) N^{1 - s},
\end{equation}

\noindent Therefore we have proved that for $k_{0}(C_{1}, A, \epsilon)$ sufficiently large, $N = 2^{k_{0}}$,

\begin{equation}\label{3.21}
\aligned
\int_{\mathcal J} \int |\frac{d}{dt} E(Iv(t))| dt \lesssim N^{2(1 - s)} N^{-\epsilon/2} C_{1}^{2} C_{0}(A, \epsilon)^{4} << N^{2(1 - s)}.
\endaligned
\end{equation}

\noindent This finally proves

\begin{equation}\label{3.24}
\int_{\mathcal J} |\frac{d}{dt} E(Iv(t))| dt << N^{2(1 - s)},
\end{equation}
and therefore for all $t \in \mathcal J$,

\begin{equation}\label{3.25}
|E(Iv(t)) - E(Iv(0))| << N^{2(1 - s)},
\end{equation}
and therefore since $E(Iv(0)) \leq C_{0}(A, \epsilon) N^{2(1 - s)}$, $\sup_{t \in \mathcal J} E(Iv(t)) \leq \frac{3}{2} C_{0}(A, \epsilon) N^{2(1 - s)}$. $\Box$

\section{Morawetz estimates}
Finally we prove proposition $\ref{p6.4}$.

\begin{proposition}[Morawetz estimates]\label{t4.1}
Suppose $v$ solves the conformal wave equation on $\mathcal J$ with $\sup_{t \in \mathcal J} E(Iv(t)) \leq 2 C_{0}(A, \epsilon) N^{2(1 - s)}$ and

\begin{equation}
\int_{\mathcal J} \int \phi(x) (\frac{\cosh |x|}{\sinh |x|}) |Iv(x,t)|^{4} dx dt \leq C_{1} N^{2(1 - s)}.
\end{equation} 
Then for $k_{0}(C_{1}, A, \epsilon)$ sufficiently large, $N = 2^{k_{0}}$, if $\phi(x) = (\frac{|x|}{\sinh |x|})^{2}$,

\begin{equation}\label{4.1}
\int_{\mathcal J} \int \phi(x) (\frac{\cosh |x|}{\sinh |x|}) |Iv(x,t)|^{4} dx dt \lesssim C_{0}(A, \epsilon) N^{2(1 - s)}.
\end{equation}
\end{proposition}
\emph{Proof:} Let $a(x) = |x|$ and let

\begin{equation}\label{4.2}
M(t) = \int Iv_{t}(x,t) (\nabla Iv(x,t) \cdot \nabla a(x) + \frac{1}{2} \Delta a(x) v(x,t)) dx.
\end{equation}
%Then let $J \subset \mathcal J$ be an interval such that

%\begin{equation}
%\int_{J} \int \phi(x) (\frac{\cosh |x|}{\sinh |x|}) (Iv(t,x))^{4} dx dt \leq 10 \sup_{t \in \mathcal J} |M(t)|.
%\end{equation}
\noindent By Hardy's inequality,

\begin{equation}\label{4.4}
\sup_{t \in \mathcal J} |M(t)| \lesssim \| \nabla Iv \|_{L_{t}^{\infty} L_{x}^{2}(\mathcal J \times \mathbf{R}^{3})} \| Iv_{t} \|_{L_{t}^{\infty} L_{x}^{2}(\mathcal J \times \mathbf{R}^{3})} \lesssim E(Iv(t)) \lesssim C_{0}(A, \epsilon) N^{2(1 - s)}.
\end{equation}
Then following the computations of \cite{Shen}, 

\begin{equation}\label{4.5}
\int_{\mathcal J} \int \phi(x) (\frac{\cosh |x|}{\sinh |x|}) (Iv(x,t))^{4} dx dt \lesssim |\int_{\mathcal J} \frac{d}{dt} M(t) dt| + |\int_{\mathcal J} \mathcal E(t) dt|.
\end{equation}
with error terms are given by

\begin{equation}\label{4.3}
\mathcal E(t) = \int (I(\phi(x) v^{3}(x,t)) - \phi(x) (I v)^{3}(x,t)) (\nabla Iv(x,t) \cdot \nabla a(x) + \frac{1}{2} \Delta a(x) v(x,t)) dx,
\end{equation}
\noindent By the fundamental theorem of calculus,

\begin{equation}\label{4.6}
|\int_{J} \frac{d}{dt} M(t) dt| \lesssim \sup_{t \in J} |M(t)| \lesssim C_{0}(A, \epsilon)^{2} N^{2(1 - s)}.
\end{equation}
Therefore it remains to estimate

\begin{equation}\label{4.7}
|\int_{\mathcal J} \mathcal E(t) dt|.
\end{equation}
Split the error into two terms,

\begin{equation}
\aligned
\mathcal E_{1}(t) = \int (I(\phi(x) v^{3}(x,t)) - \phi(x) (I v)^{3}(x,t)) (\nabla Iv(x,t) \cdot \nabla a(x)) dx, \\
\mathcal E_{2}(t) = \int \int (I(\phi(x) v^{3}(x,t)) - \phi(x) (I v)^{3}(x,t)) \frac{1}{2} \Delta a(x) v(x,t)) dx.
\endaligned
\end{equation}
By direct calculation, $\nabla a(x) = \frac{x}{|x|}$, and also by definition of the energy and corollary $\ref{c2.2}$, $\nabla Iv$ and $Iv_{t}$ have the same estimates, and thus the terms in $\mathcal E_{1}(t)$ may be estimated in a manner which is exactly analogous to the corresponding terms in the previous section.\vspace{5mm}

\noindent Indeed, as in $(\ref{3.23})$,

\begin{equation}\label{4.21}
\aligned
\int_{\mathcal J} \int \phi(x) (P_{> k_{0} - 7} Iv)(Iv)^{2} (\nabla a(x) \cdot \nabla I P_{\leq k_{0}} v) dx dt \\
\lesssim \ln(N) N^{1/2} \| \phi(x)^{1/2} (Iv)^{2} \|_{l_{k_{0}}^{\infty} L_{t,x}^{2}(\mathcal J \times \mathbf{R}^{3})} \| P_{> k_{0} - 7} Iv \|_{l_{k_{0}}^{\infty} L_{t,x}^{2}(\mathcal J \times \mathbf{R}^{3})} \| P_{\leq k_{0}} \nabla Iv \|_{L_{t}^{\infty} L_{x}^{2}(\mathcal J \times \mathbf{R}^{3})} \\
\lesssim \ln(N) N^{1/2 - s} C_{0}(A, \epsilon)^{2} C_{1}^{1/2} N^{2(1 - s)}.
\endaligned
\end{equation}
Next, as in $(\ref{3.14})$,

\begin{equation}\label{4.20}
\aligned
\int_{\mathcal J} \int I(\phi(x) (P_{> k_{0} - 7} v)v^{2}) (\nabla a(x) \cdot \nabla I P_{\leq k_{0}} v) dx dt \\ 
\lesssim \ln(N) N^{1/2} \| \nabla P_{\leq k_{0}} Iv \|_{L_{t}^{\infty} L_{x}^{2}(\mathcal J \times \mathbf{R}^{3})} \| \phi(x)^{1/2} v^{2} \|_{l_{k_{0}}^{\infty} L_{t,x}^{2}(\mathcal J \times \mathbf{R}^{3})} \| P_{\geq k_{0} - 7} v \|_{l_{k_{0}}^{\infty} L_{t,x}^{2}(\mathcal J \times \mathbf{R}^{3})} \\
\lesssim \ln(N) N^{1/2} (C_{1}^{1/2} N^{1 - s} + C_{0}(A, \epsilon) N^{1/2 - 2 \epsilon}) \| \nabla I v \|_{L_{t}^{\infty} L_{x}^{2}(\mathcal J \times \mathbf{R}^{3})} \sum_{j > k_{0} - 7} 2^{-js} \mathcal S(j, v) \\
\lesssim \ln(N) N^{1/2 - s} C_{0}(A, \epsilon)^{2} N^{1 - s} (C_{1}^{1/2} N^{1 - s} + C_{0}(A, \epsilon) N^{1/2 - 2 \epsilon}).
\endaligned
\end{equation}
As in $(\ref{3.22.1})$,

\begin{equation}\label{4.22}
\aligned
\int_{\mathcal J} \int_{|x| \leq \frac{2}{N}} \phi(x) (P_{> k_{0} - 7} Iv) (Iv)^{2} (\nabla a(x) \cdot \nabla I P_{> k_{0}} v) dx dt \\
\lesssim \| Iv \|_{L_{t,x}^{4}(\mathcal J \times B_{\frac{2}{N}})}^{2} \| P_{> k_{0} - 7} Iv \|_{L_{t}^{2} L_{x}^{\infty}(\mathcal J \times \mathbf{R}^{3})} \| \nabla IP_{> k_{0}} v \|_{L_{t}^{\infty} L_{x}^{2}(\mathcal J \times \mathbf{R}^{3})} \\ \lesssim C_{0}(A, \epsilon)^{2} N^{2(1 - s)} N^{-1/2} \| (Iv)^{2} \|_{l_{k_{0}}^{\infty} L_{t,x}^{2}(\mathcal J \times \mathbf{R}^{3})} \\ \lesssim C_{0}(A, \epsilon)^{2} C_{1}^{1/2} N^{\frac{3}{2} - 2s} N^{1 - s}.
\endaligned
\end{equation}
Also, as in $(\ref{3.22})$,

\begin{equation}\label{4.23}
\aligned
\int_{\mathcal J} \int_{|x| \geq \frac{2}{N}} \phi(x) (P_{> k_{0} - 7} Iv)(Iv)^{2} (\nabla a(x) \cdot \nabla I P_{> k_{0}} v) dx dt \\ \lesssim \ln(N) \| \phi(x)^{1/2} (Iv)^{2} \|_{l_{k_{0}}^{\infty} L_{t,x}^{2}(\mathcal J \times \mathbf{R}^{3})} \| \nabla I P_{> k_{0}} v \|_{l_{k_{0}}^{\infty} L_{t,x}^{2}(\mathcal J \times \mathbf{R}^{3})} \| |x| (I P_{> k_{0}} v) \|_{L_{t,x}^{\infty}(\mathcal J \times \mathbf{R}^{3})} \\
\lesssim \ln(N) C_{0}(A, \epsilon) C_{1}^{1/2} N^{2(1 - s)} \sum_{j > k_{0} - 7} 2^{-j/2} \| P_{j} Iv \|_{L_{t}^{\infty} \dot{H}^{1}(\mathcal J \times \mathbf{R}^{3})} \\
\lesssim \ln(N) C_{0}(A, \epsilon)^{2} C_{1}^{1/2} N^{2(1 - s)} N^{1/2 - s}.
\endaligned
\end{equation}
%Summing in $j$,

%\noindent By Hardy's inequality,

%\begin{equation}\label{4.24}
%\sum_{j > 0} (\ref{4.23}) \lesssim \ln(N) N^{1 - s} N^{-\epsilon} \| \frac{\phi(x)^{1/4}}{|x|^{1/4}} Iv \|_{L_{t,x}^{4}(J \times \mathbf{R}^{3})}^{2}.
%\end{equation}

\noindent As in $(\ref{3.11})$,

\begin{equation}\label{4.31}
\aligned
\int_{\mathcal J} \int I(\phi(x) (P_{> k_{0} - 7} v)^{3}) (\nabla a(x) \cdot \nabla I P_{\geq k_{0}} v) dx dt \\ \lesssim \| I(\phi(x) (P_{> k_{0} - 7} v)^{3}) \|_{L_{t}^{1} L_{x}^{2}(\mathcal J \times \mathbf{R}^{3})} \| \nabla I v \|_{L_{t}^{\infty} L_{x}^{2}(\mathcal J \times \mathbf{R}^{3})} \\ \lesssim C_{0}(A, \epsilon)^{3} N^{1 - s - \epsilon} (\sup_{t \in \mathcal J} E(Iv(t)))^{1/2} \lesssim C_{0}(A, \epsilon)^{4} N^{2(1 - s)} N^{-\epsilon}.
\endaligned
\end{equation}

\noindent Next, following $(\ref{3.17})$,

\begin{equation}\label{4.32}
\aligned
\int_{\mathcal J} \int I(\phi(x) v(x,t) (P_{> k_{0} - 7} v)(P_{\leq k_{0}} v))(x,t)  (\nabla I P_{> k_{0}} Iv)(x,t)  dx dt \\
\lesssim \ln(N) C_{0}(A, \epsilon) N^{1 - s} \| \nabla I P_{\geq k_{0}} v \|_{l_{k_{0}}^{\infty} L_{t,x}^{2}(\mathcal J \times \mathbf{R}^{3})} \| \phi(x)^{1/2} (P_{> k_{0} - 7} v) v \|_{L_{t,x}^{2}(\mathcal J \times \mathbf{R}^{3})} \\
\lesssim \ln(N) N^{2(1 - s)} C_{0}(A, \epsilon)^{2} (\ln(N) N^{1 - s} \| P_{> k_{0} - 7} v \|_{l_{k_{0}}^{\infty} L_{t,x}^{2}(\mathcal J \times \mathbf{R}^{3})} + \| P_{> k_{0} - 7} v \|_{L_{t,x}^{4}(\mathcal J \times \mathbf{R}^{3})}^{2}) \\
\lesssim \ln(N) C_{0}(A, \epsilon)^{4} N^{2(1 - s)} (\ln(N) N^{1 - 2s} + N^{-2 \epsilon}).
\endaligned
\end{equation}

\noindent Finally, by $(\ref{3.7})$,

\begin{equation}\label{4.11}
\aligned
\int_{\mathcal J} \int (I - 1) (\phi(x) (P_{\leq k_{0} - 7} v)^{3}) (\nabla a(x) \cdot \nabla I v(x,t)) dx dt \\
\lesssim \| (I - 1) (\phi(x) (P_{\leq k_{0} - 7} v)^{3}) \|_{L_{t}^{1} L_{x}^{2}(\mathcal J \times \mathbf{R}^{3})} \| \nabla Iv \|_{L_{t}^{\infty} L_{x}^{2}(\mathcal J \times \mathbf{R}^{3})} \lesssim C_{1}^{2} C_{0}^{4}(A, \epsilon) (N^{4 - 6s} + N^{3 - 4s}).
\endaligned
\end{equation}
Therefore, we have proved that

\begin{equation}
\int_{\mathcal J} |\mathcal E_{1}(t)| dt \lesssim N^{2(1 - s)} N^{-\epsilon/2} C_{1}^{2} C_{0}^{4}(A, \epsilon).
\end{equation}
Also by $(\ref{3.7})$, $\Delta a(x) \lesssim \frac{1}{|x|}$, and Hardy's inequality,

%\begin{equation}\label{4.12}
%(I - 1) ((P_{\leq \frac{N}{100}} \phi(x)^{1/2})^{2} (P_{\leq \frac{N}{100}} v)^{3}) = 0,
%\end{equation}

\begin{equation}\label{4.13}
\aligned
\int_{\mathcal J} \int (I - 1) (\phi(x) (P_{\leq k_{0} - 7} v)^{3}) [\frac{1}{2} \Delta a(x) Iv(x,t)] dx dt \\ \lesssim  \| (I - 1) (\phi(x) (P_{\leq k_{0} - 7} v)^{3}) \|_{L_{t}^{1} L_{x}^{2}(\mathcal J \times \mathbf{R}^{3})} \| \frac{1}{|x|} Iv \|_{L_{t}^{\infty} L_{x}^{2}(\mathcal J \times \mathbf{R}^{3})} \\ \lesssim C_{1}^{2} C_{0}^{4}(A, \epsilon) (N^{4 - 6s} + N^{3 - 4s}).
\endaligned
\end{equation}
Therefore, it only remains to estimate

\begin{equation}
\int_{\mathcal J} \int \phi(x) (Iv)^{2} (I P_{> k_{0} - 7} v) \frac{1}{|x|} (Iv) dx dt,
\end{equation}
and

\begin{equation}
\int_{\mathcal J} \int I(\phi(x) v^{2} (P_{> k_{0} - 7} v)) \frac{1}{|x|} (Iv) dx dt.
\end{equation}
First, by $(\ref{5.2.1})$, proposition $\ref{p2.1}$, $(\ref{3.16})$, $(\ref{3.16.1})$, and the radial Sobolev embedding theorem,

\begin{equation}\label{4.19}
\aligned
 \int_{\mathcal J} \int_{|x| > \frac{2}{N}} I(\phi(x) (P_{> k_{0} - 7} v) v^{2}) \frac{1}{|x|} (Iv) dx dt \\
\lesssim N^{1/2} \ln(N) \| P_{> k_{0} - 7} v \|_{l_{k_{0}}^{\infty} L_{t,x}^{2}(\mathcal J \times \mathbf{R}^{3})} \| \phi(x)^{1/2} v^{2} \|_{l_{k_{0}}^{\infty} L_{t,x}^{2}(\mathcal J \times \mathbf{R}^{3})} \| |x|^{1/2} Iv \|_{L_{t,x}^{\infty}(\mathcal J \times \mathbf{R}^{3})} \\
\lesssim N^{3/2 - 2s} \ln(N) C_{0}(A, \epsilon)^{2} (C_{1}^{1/2} N^{1 - s} + C_{0}(A, \epsilon)^{2} N^{1/2 - 2 \epsilon}).
\endaligned
\end{equation}
%Now by proposition $\ref{p2.1}$, theorem $\ref{t5.3}$, and the fact that $\phi(x) \lesssim e^{-c |x|}$ for some $c > 0$,

\noindent Next, following $(\ref{3.11})$, by Hardy's inequality,

\begin{equation}\label{4.18.1}
\aligned
\int_{\mathcal J} \int_{|x| \leq \frac{2}{N}} I(\phi(x) (P_{> k_{0} - 7} v)^{3})(x,t) \frac{1}{|x|} (Iv)(x,t) dx dt \\
\lesssim \| I(\phi(x) (P_{> k_{0} - 7} v)^{3}) \|_{L_{t}^{1} L_{x}^{2}(\mathcal J \times \mathbf{R}^{3})} \| \frac{1}{|x|} Iv \|_{L_{t}^{\infty} L_{x}^{2}(\mathcal J \times \mathbf{R}^{3})} \lesssim C_{0}(A, \epsilon)^{4} N^{-\epsilon} N^{2(1 - s)}.
\endaligned
\end{equation}
Also, by theorem $\ref{t5.5}$, $(\ref{3.16})$, $(\ref{3.16.1})$, the fact that $\Delta a(x) \lesssim \frac{1}{|x|}$, the support of $\psi(2^{k_{0}} x)$, and proposition $\ref{p2.1}$,

\begin{equation}\label{4.18.2}
\aligned
\int_{J} \int I(\phi(x) (P_{> k_{0} - 7} v)(P_{\leq k_{0} - 7} v) v) \frac{1}{|x|} \psi(2^{k_{0}} x) (Iv)(x,t) dx dt \\
\lesssim \| I(\phi(x) (P_{> k_{0} - 7} v)(P_{\leq k_{0} - 7} v) v) \psi(2^{k_{0}} x) \frac{1}{|x|^{3/4}} \|_{L_{t,x}^{4/3}(J \times \mathbf{R}^{3})} \| \frac{\phi(x)^{1/4}}{|x|^{1/4}} Iv \|_{L_{t,x}^{4}(J \times \mathbf{R}^{3})}
\endaligned
\end{equation}

\begin{equation}\label{4.18.3}
\aligned
\lesssim N^{-s} C_{0}(A, \epsilon) C_{1}^{1/4} N^{\frac{1 - s}{2}} \| v^{2} \|_{l_{k_{0}}^{\infty} L_{t,x}^{2}(\mathcal J \times \mathbf{R}^{3})}^{1/2} \| P_{\leq k_{0} - 7} v \|_{L_{t,x}^{\infty}(\mathcal J \times \mathbf{R}^{3})} \\
\lesssim C_{0}(A, \epsilon)^{2} C_{1}^{1/4} N^{\frac{1 - s}{2}} N^{1/2 - s} N^{1 - s} (C_{1}^{1/4} N^{\frac{1 - s}{2}} + C_{0}(A, \epsilon) N^{1/4 - \epsilon}).
\endaligned
\end{equation}

\noindent Finally, by $\Delta a(x) \lesssim \frac{1}{|x|}$, $\frac{1}{|x|} \lesssim (\frac{\cosh |x|}{\sinh |x|})$, Hardy's inequality and proposition $\ref{p2.1}$,

\begin{equation}\label{4.16}
\aligned
\int_{J} \int \phi(x) (Iv(x,t))^{2} (I P_{> k_{0} - 7} v)(x,t) (\Delta a(x) Iv(x,t)) dx dt \\ \lesssim \| \frac{\phi(x)^{1/4}}{|x|^{1/4}} Iv \|_{L_{t,x}^{4}(J \times \mathbf{R}^{3})}^{3} \| \frac{1}{|x|^{1/2}} (IP_{\geq k_{0} - 7} v) \|_{L_{t}^{\infty} L_{x}^{2}(J \times \mathbf{R}^{3})}^{1/2} \| I P_{\geq k_{0} - 7} v \|_{L_{t}^{2} L_{x}^{\infty}(J \times \mathbf{R}^{3})}^{1/2} \\
\lesssim N^{-\epsilon/2} N^{\frac{1 - s}{2}} \| \frac{\phi(x)^{1/4}}{|x|^{1/4}} Iv \|_{L_{t,x}^{4}(J \times \mathbf{R}^{3})}^{3} \lesssim C_{1}^{3/4} N^{-\epsilon} N^{2(1 - s)}.
\endaligned
\end{equation}
Therefore, by $(\ref{4.19})$ - $(\ref{4.16})$,

\begin{equation}\label{4.25}
\int_{\mathcal J} |\mathcal E_{2}(t)| dt << N^{2(1 - s)},
\end{equation}

\noindent which proves the proposition. $\Box$

\nocite*
\bibliographystyle{plain}

\begin{thebibliography}{[00]}
\bibitem{AM}
	\newblock C. Antonini and F. Merle,
	\newblock ``Optimal bounds on positive blow - up solutions for a semi linear wave equation",
	\newblock \textit{International Mathematics Research Notices} (2001) no. 21, 1141 -- 1167.

\bibitem{B}
	\newblock J. Bourgain,
	\newblock ``Refinements of {S}trichartz' inequality and applications to {$2$}{D} - {NLS} with critical nonlinearity",
	\newblock \textit{International Mathematics Research Notices} \textbf{5} (1998) 253 -- 283.		

\bibitem{CCT}
	\newblock J. Colliander, M. Christ, and T. Tao,
	\newblock ``Asymptotics, frequency modulation, and low regularity ill - posedness for canonical defocusing equations",
	\newblock \textit{American Journal of Mathematics} \textbf{125} 6 (2003) 1235 -- 1293.
	
\bibitem{CKSTT}
	\newblock C. Colliander, M. Keel, G. Staffilani, H. Takaoka, and T. Tao
	\newblock ``Almost conservation laws and global rough solutions to a nonlinear {S}chrodinger equation",
	\newblock \textit{Mathematics Research Letters} \textbf{9} 5 -- 6 (2002) 659 -- 682.
	
	
\bibitem{D1}
	\newblock B. Dodson,
	\newblock ``Global well - posedness for the defocusing, cubic, nonlinear wave equation in three dimensions for radial initial data in $\dot{H}^{s} \times \dot{H}^{s - 1}$, $s > \frac{1}{2}$",
	\newblock Preprint, arXiv:1506.06239.
	
\bibitem{DL}
	\newblock B. Dodson and A. Lawrie,
	\newblock ``Scattering for the radial 3d cubic wave equation",
	\newblock \textit{Analysis and PDE} \textbf{8} (2015) 467 -- 497.		

\bibitem{GV}
	\newblock J. Ginibre and G. Velo
	\newblock ``Generalized {S}trichartz inequalities for the wave equation",
	\newblock \textit{Journal of Functional Analysis} \textbf{133} 1 (1995) 50 -- 68.	
	
\bibitem{Gril}
	\newblock M. Grillakis,
	\newblock ``Regularity and asymptotic behaviour of the wave equation with critical nonlinearity",
	\newblock \textit{Annals of Mathematics} 132 (1990) 485 -- 509.	

\bibitem{KTwm}
	\newblock M. Keel and T. Tao
	\newblock ``Local and global well posedness of wave maps on $\mathbf{R}^{1 + 1}$ for rough data'',
	\newblock \textit{International Mathematics Research Notices}, \textbf{21} (1998) 1117 -- 1156.	

\bibitem{KPV}
	\newblock C. Kenig, G. Ponce, and L. Vega,
	\newblock ``Global well - posedness for semi - linear wave equations",
	\newblock \textit{Communications in Partial Differential Equations} \textbf{25} 9 - 10 (2000) 1741 -- 1752.	


\bibitem{KM}
	\newblock S. Klainerman and M. Machedon,
	\newblock ``Space - time estimates for null forms and the local existence theorem",
	\newblock \textit{Communications on Pure and Applied Mathematics} \textbf{46} 9 (1993) 1221 -- 1268.	


\bibitem{LS}
	\newblock H. Lindblad and C. Sogge,
	\newblock ``On existence and scattering with minimal regularity for semilinear wave equations",
	\newblock \textit{Journal of Functional Analysis} \textbf{130} (1995) 357 -- 426.

\bibitem{Roy}
	\newblock T. Roy,
	\newblock ``Global well - posedness for the radial defocusing cubic wave equation on {$\Bbb R^{3}$} and for rough data",
	\newblock \textit{Electronic Journal of Differential Equations} No. 166 (2007).
	
\bibitem{Shen}
	\newblock R. Shen,
	\newblock ``Scattering of solutions to the defocusing energy sub - critical semi - linear wave equation in 3{D}",
	\newblock Preprint, arXiv:1512.00705.
	
\bibitem{Stri}
	\newblock R. S. Strichartz,
	\newblock ``Restrictions of {F}ourier transforms to quadratic surfaces and decay of solutions of wave equations'',
	\newblock \textit{Duke Mathematical Journal} \textbf{44} no. 3 (1977) 705 - 714.

\bibitem{Tao} 
     	\newblock T. Tao,
     	\newblock \textit{Nonlinear Dispersive Equations. Local and Global Analysis},
     	\newblock CBMS Regional Conference Series in Mathematics \textbf{104} Published for the Conference Board of the Mathematical Sciences, Washington, DC, 2006.

\bibitem{Taylor}
	\newblock M. Taylor,
	\newblock ``Tools for {PDE} : Pseudodifferential operators, paradifferential operators, and layer potentials'',
	\newblock \textit{American Mathematical Society} Providence, RI, 2000.
	

\end{thebibliography}

\end{document}